\documentclass[10pt]{article}
\usepackage{hyperref}
\usepackage{amssymb}
\usepackage{amsmath}
\usepackage{tikz,tikz-cd}
\usepackage{authblk}
\usepackage{color}

\input{liemacs10.sty} 
\newcommand{\sV}{{\tt V}}

\newcommand{\bD}{\mathbb D}

\newcommand{\Out}{\mathop{{\rm Out}}\nolimits}

\renewcommand{\phi}{\varphi}

\newcommand{\sech}{\mathop{{\rm sech}}\nolimits}

\renewcommand{\phi}{\varphi}
\newcommand{\bracket}[2]{\langle #1, #2\rangle}
\renewcommand{\mlabel}{\label} 

\begin{document} 

\title{Reflection positivity  and its relation to \\ 
  disc, half plane and the strip}

\author[1]{Maria Stella Adamo}
\author[2]{Karl-Hermann Neeb}
\author[3]{Jonas Schober}
\date{}

\affil[1,2]{Department Mathematik, Friedrich-Alexander-Universit\"at Erlangen-N\"urnberg \authorcr
Cauerstrasse 11, 91058 Erlangen, Germany}
\affil[3]{Instituto de Investigaciones en Matem\'{a}ticas Aplicadas y en Sistemas,
Universidad Nacional Aut\'{o}noma de M\'{e}xico, Mexico City, Mexico}

\maketitle

\abstract{We develop a novel perspective
  on reflection positivity (RP) on the strip 
  by systematically developing the analogies with
  the unit disc and the upper half plane in the complex plane.
  These domains correspond to the three conjugacy classes
  of one-parameter groups in the M\"obius group
  (elliptic for the disc, parabolic for the upper half plane
  and hyperbolic for the strip). In all cases, reflection
  positive functions correspond to positive functionals
  on $H^\infty$ for a suitable involution.
  For the strip, reflection positivity naturally connects with
  Kubo--Martin--Schwinger (KMS) conditions on the real line
  and further to standard pairs, as they appear in Algebraic
  Quantum Field Theory. We also exhibit a curious connection between
  Hilbert spaces on the strip and the upper half plane, based on a
  periodization process.
}  

\tableofcontents 

\vspace{1cm}

\section*{Introduction} 
\mlabel{sec:0}

In the complex plane the unit disc
\[ \bD := \{  z \in \C \: |z| < 1 \},\]
the upper half plane
\[ \C_+ := \{ z \in \C \: \Im z > 0\}\]
and the strips
\[ \bS_\beta := \{ z \in \C \: 0 < \Im z < \beta \},\quad \beta > 0 \]
are simply connected proper domains, hence biholomorphically
equivalent by the Riemann Mapping Theorem.  So they all
look alike from the complex analytic perspective and
their automorphism groups are isomorphic to $\PSL_2(\R)$. However, 
each of these domains puts a focus on a specific one-parameter group 
whose action is most natural for the domain: 
\begin{itemize}
\item[\rm(G1)] on the disc $\bD$ it is the rotation group $\T$ fixing $0$,
  generated by an elliptic Lie algebra element.
\item[\rm(G2)] on the upper half plane $\C_+$  it is the translation
  group $\R$, generated by a nilpotent Lie algebra element.
\item[\rm(G3)] on the strip $\bS_\beta$ it is the translation 
  group $\R$,  generated by a hyperbolic Lie algebra element.
\end{itemize}

In the following we write $A$ ($\cong \T$ or $\R$) for the corresponding 
one-parameter group of transformations of the domain.
From this perspective, the domains are rather different, and this is reflected
in various aspects of harmonic analysis.
For $\bD$ and $\C_+$ the
group of biholomorphic maps acts by M\"obius transformations,
but this is different for the strip $\bS_\beta$. Even though
  the translation group acts by M\"obius transformations, this is
  not so for the full automorphism group of $\bS_\beta$.

Each proper simply connected domain $\Omega$ in the complex plane
carries an antiholomorphic involution $\sigma$, which is unique up to
conjugation under automorphisms (\cite[Lemma~B.1]{ANS22}). We pick
$\sigma(z) = \oline z$ for the disc~$\bD$, 
$\sigma(z) = -\oline z$ for the upper half plane $\C_+$, and
$\sigma(z) = \beta i + \oline z$ for the strip $\bS_\beta$.
On $\bD$ and $\C_+$, the antiholomorphic involution  $\sigma$
anticommutes with the specified one-parameter group
from (G1-3), but on $\bS_\beta$ it commutes with the translation group $\R$. 
In any case, we obtain by
\[ f^\sharp(z) := \oline{f(\sigma(z))} \]
on the algebra $H^\infty(\Omega)$ of bounded holomorphic functions on
$\Omega$ the structure of a Banach $*$-algebra.
As it has a unique predual, it carries a natural weak topology. 
We refer to  \cite[App.~B]{ANS22} for a detailed discussion of this
Banach $*$-algebra, including a description of its weakly continuous
states by probability measures on the fixed point set $\Omega^\sigma$.
 All these Banach $*$-algebras are
  isomorphic, but for our three domains, they are topologically
  generated (with respect to the weak topology) by different
  subsemigroups:
  \begin{itemize}
  \item[\rm(S1)]  For the disc it is generated by the $\sharp$-symmetric semigroup
    $(z^n)_{n \in \N_0}$ of monomials. 
  \item[\rm(S2)] For the upper half plane $\C_+$, it is generated
    by the one-parameter semigroup $(e_{it})_{t \geq 0}$, where
    ${e_{it}(z) := e^{itz}}$. Note that this semigroup is $\sharp$-symmetric
      in the sense that $e_{it}^\sharp = e_{it}$ for $t \geq 0$. 
\item[\rm(S3)] For $\bS_\beta$ is its generated by the
  $\sharp$-unitary one-parameter group
   $(g_t)_{t \in \R}$ with $g_t = e^{\beta t/2} e_{it}$
    satisfying $g_t^\sharp = g_{-t}$ for $t \in \R$.    
\end{itemize}
In all three cases, we have a natural homomorphism
  of the character group $\hat A = \Hom(A,\T)$ of the one-parameter group $A$
  introduced in (G1-3)
  to the multiplicative group of $L^\infty(\partial \Omega)$, 
  and a subsemigroup
  $\hat A_+$ (see (S1-3)
mapped into $H^\infty(\Omega)$, whose range spans a
    weakly dense subspace.
For $\Omega = \bD$, the dual
$\hat\T = \Z$ is discrete with
$\hat \bT_+ = \N_0$, and for $\Omega = \C_+, \bS_\beta$,
we have $\hat\R \cong \R$, but $\hat A_+ = \R_+$ for $\C_+$ and
$\hat A_+ = \R$ for $\bS_\beta$. For $\bD$ and $\C_+$ we thus
obtain a so-called symmetric semigroup $(\hat A, \hat A_+, -\id)$.
\begin{footnote}
{A triple $(G,S,\tau)$ is called a 
{\it symmetric semigroup} if $G$ is a group,
$\tau$ is an involutive automorphism of $G$, and 
$S \subeq G$ is a subsemigroup satisfying $\tau(S)^{-1} = S$.}
\end{footnote} 
For the strip $\bS_\beta$, one has to adopt a different perspective
to see a corresponding structure. Here one may consider the triple
$(\T_\beta, \T_{\beta,+}, -\id)$, where $\T_\beta \cong \R/\beta \Z$
is the circle group and $\T_{\beta,+}$ is the image of $[0,\beta/2]$ in $\T_\beta$.
This half circle is not a subsemigroup, but it can be used to specify
so-called reflection positive representations
and reflection positive functions on $\T_\beta$.
We refer to \cite{NO18} for the basics and the physical aspects
  of reflection   positivity.
It is a central goal of the present paper to shed some new light on
how reflection positivity of the triple $(\T_\beta, \T_{\beta,+}, -\id)$
is related to the strip $\bS_\beta$, and to which extent this is analogous
to what we see more directly for the disc $\bD$ and the
upper half plane $\C_+$.\\

The {\bf content of this paper} is as follows.
In Section~\ref{sec:1} we briefly introduce reflection positive 
functions $\phi$ for triples $(G,G_+, \tau)$, as those positive definite
functions $\phi \: G \to \C$ for which the kernel
$(\phi(s\tau(t)^{-1}))_{s,t \in G_+}$ is also positive definite. Here we shall
only deal with the triples
\[ (\Z, \N_0, - \id_\Z), \quad 
(\R, \R_+, - \id_\R) \quad  \mbox{ and } \quad
(\T_\beta, \T_{\beta,+}, - \id_{\T_\beta}).\]
Via a suitable variant of the Gelfand--Naimark--Segal (GNS) construction,
reflection positive functions lead to representations
on {\it reflection positive Hilbert spaces} $(\cE,\cE_+, \theta)$, where
$\cE$ is a Hilbert space, $\theta$ a unitary involution on $\cE$, and $\cE_+$
a closed subspace which is {\it $\theta$-positive} in the sense that
$\la v,v\ra_\theta := \la v,\theta v \ra \geq 0$ for $v \in \cE_+$. We write
$\hat\cE$ for the Hilbert space defined by $\la \cdot,\cdot\ra_\theta$ on
$\cE_+$ by factorization of null vectors and completion
(see \cite{NO18} for more details). 

Section~\ref{sec:3a} connects this concept with simply connected domains
$\Omega \subeq \C$ by showing that any $w \in \Omega^\sigma$ leads to a
reflection positive Hilbert space 
$(L^2(\partial \Omega), H^2(\Omega), \theta_w)$ for which
$\hat\cE$ is one-dimensional, corresponding to a character of the
$*$-algebra $(H^\infty(\Omega),\sharp)$ (Proposition~\ref{prop:2.3}).
In general, reflection positive representations 
associated to the domain $\Omega$ lead to weakly continuous
representations of the $*$-algebra $(H^\infty(\Omega), \sharp)$ on $\hat\cE$
(Proposition~\ref{pos_funct}). 

Sections~\ref{sec:3} and \ref{sec:2} briefly discuss the relevant structures
for the upper half plane and the disc, respectively.
Specializing Proposition~\ref{prop:2.3},
Proposition~\ref{prop:3.1} identifies the weakly continuous characters of
$(H^\infty(\C_+), \sharp)$ with the primitive reflection positive 
functions on $(\R,\R_+, -\id_\R)$, and the positive weakly continuous functionals
on $(H^\infty(\C_+), \sharp)$ with general reflection positive functions.
The same mechanism connects for the disc $\bD$ positive functionals on
$(H^\infty(\bD), \sharp)$ with reflection positive functions on
$(\Z,\Z_+, -\id_\Z)$ (Theorem~\ref{thm:3.2}). 

In Section~\ref{sec:4} we eventually turn to the strip $\bS_\beta$.
We first show that the states of $(H^\infty(\bS_\beta),\sharp)$
are in one-to-one correspondence with normalized positive definite
functions on $\R$ by restricting to the subgroup
$(g_t)_{t \in \R}$ from  (S3). 
The connection to reflection positive functions on
$(\T_\beta, \T_{\beta,+}, -\id)$ is exhibited in Theorem~\ref{thm:5.7}.
A key property is that reflection positive functions
$\phi_\T$ on
$\T_\beta$, considered as functions on the vertical segment
$[0, \beta]i \subeq \C$, extend analytically to holomorphic functions on
$\bS_\beta$. Restricting to $\R \subeq \partial \bS_\beta$ then leads to
positive definite functions $\phi_\R \: \R \to \C$ satisfying the
{\it $\beta$-KMS condition} 
\begin{equation}
  \label{eq:beta-kms}
 \phi_\R(\beta i + t) = \oline{\phi_\R(t)} \quad \mbox{ for } \quad
 t \in \R,
\end{equation}
(\cite{NO19}, \cite{RvD77}). For the measure $\nu$ with Fourier transform
$\hat\nu = \phi_\R$, this corresponds to the
{\it $\beta$-reflection relation} 
\begin{equation}
  \label{eq:nurel-intro} d\nu(-\lambda) = e^{-\beta \lambda} d\nu(\lambda).
\end{equation}
This translates reflection positivity on the circle $\T_\beta$
into the $\beta$-KMS condition on~$\R$.
To describe the natural models for the corresponding unitary representations
of~$\R$, we observe in Theorem~\ref{thm:7.1} that
tempered measures $\nu$ on $\R$, satisfying the
$2\beta$-reflection relation, are in one-to-one correspondence to
reproducing kernel Hilbert spaces of holomorphic functions on $\bS_\beta$
on which $\R$ acts unitarily by translations. The corresponding kernel
is given by $K(z,w) = \hat\nu(z - \oline w)$.
Typical examples of such Hilbert spaces are the
Hardy and  the Bergman space on $\bS_\beta$.

This leads us to another aspect of the strip that is invisible for
the disc and the upper half plane. We recall 
that a closed real subspace $\sV$ of a complex Hilbert space $\cH$ 
is called {\it standard} if 
  $\sV + i \sV$ is dense and $\sV \cap i \sV = \{0\}$.
If $\sV$ is standard, then
the complex conjugation $T_\sV(x+ iy) = x-iy$ on $\sV + i \sV$ has a polar decomposition
$T_\sV = J_\sV \Delta_\sV^{1/2}$, where
$J_\sV$ is a conjugation (an antilinear involutive isometry)
and $\Delta_\sV$ is a positive selfadjoint operator
satisfying $J_\sV \Delta_\sV J_\sV =~\Delta_\sV^{-1}$.
The unitary one-parameter group $(\Delta_\sV^{it})_{t \in \R}$,
the {\it modular group of $\sV$}, preserves the subspace~$\sV$
(see \cite{Lo08} for background on standard subspaces). 
In \cite[Thm.~2.6]{NO19} it is shown that the
$\beta$-KMS condition for a positive definite function
$\psi$ on $\R$ is equivalent to the existence 
of a standard subspace $\sV$ and some $v \in \sV$ such that
\[ \psi(t) = \la v, \Delta^{-it/\beta} v \ra \quad \mbox{ for } \quad t \in \R.\]
So standard subspaces arise naturally by the GNS construction
from positive definite functions satisfying the $\beta$-KMS condition.
This connection is explained in Section~\ref{subsec:5.4}, where
we also discuss some key examples, such as
$\sV = H^2(\bS_\beta)^\sharp \subeq L^2(\R)$ for the translation action.
We conclude Section~\ref{sec:4} with a discussion of the
Cauchy--Szeg\"o kernel and the Poisson kernel of $\bS_\beta$
because they relate naturally with reflection positivity on $\T_\beta$,
resp., $\beta$-KMS conditions on $\R$. 

In Section~\ref{sec:6} we study the passage from functions
on the upper half plane to the strip by periodization, using that
\[ \bS_\beta = \C_+ \cap (\beta i - \C_+).\]
We show in particular that, if $(Q_w^{\C_+})_{w \in \C_+}$ are the evaluation functionals
of the Hardy space $H^2(\C_+)$ and $w \in \bS_\beta$, then
the functions $Q_w^{\C_+}$ on $\C_+$ and $Q_w^{\bS_\beta}$ on $\bS_\beta$ 
  extend meromorphically to $\C$, and in this sense we have
 the identity
  \[  Q_w^{\bS_\beta} = \sum_{k \in \Z} (-1)^k Q_w^{\C_+}(z + 2 ki \beta) \]
  of meromorphic functions (cf.~\eqref{eq:szego-exp2}).
  So  the reproducing kernel $Q^{\bS_\beta}$ of the
Hardy space $H^2(\bS_\beta)$ on the strip has a corresponding
series expansion. A similar result
holds for the Bergman spaces $B^2(\Omega) = L^2(\Omega)\cap\cO(\Omega)$, 
but the Bergman kernel $Q^B(z,w) = Q^B_w(z)$ on $\bS_\beta$ 
is a multiple of a series of the form
$\sum_{k \in \Z} Q_w^B(z + 2 ki \beta)$ without the sign changes. 

In Section~\ref{sec:7} we discuss a unifying perspective
  on normal form results related closely to $\bD$, $\C_+$ and $\bS_\beta$,
  namely for isometries on Hilbert spaces
  (that lead to representations of $H^\infty(\bD)$),
  for one-parameter semigroups of isometries
  (that lead to representations of $H^\infty(\C_+)$),
  and for so-called standard pairs $(U,\sV)$, 
  consisting of a standard subspace $\sV$ of a complex Hilbert space~$\cH$
  and a  unitary one-parameter group $(U_t)_{t \in \R}$ on $\cH$, such that 
$U_t \sV \subeq \sV$ for $t \geq 0$ and $U_t = e^{itH}$ with 
$H \geq 0$. In particular, we explain how the corresponding normal
form results can be derived in a uniform fashion from the
Mackey--Stone--von Neumann Theorem on projective
unitary representations of the product group $A \times \hat A$.

Finally, in the short Section~\ref{sec:8}, we observe that the three
generating subsemigroups in $H^\infty(\C_+)$, $H^\infty(\bD)$ and $H^\infty(\bS_\beta)$
exhibit also some reflection positivity properties in their
parameter. Presently, we do not have any conceptual explanation for this phenomenon.

The following table contains much of the relevant data that we shall
use below. Most of the notation has already been explained in
the introduction, the other items will be explained below.\\[5mm]

\hspace{-24mm}
\begin{tabular}{||l||l|l|l||}\hline
domain $\Omega$ & $\bD$ & $\C_+$ & $\bS_\beta$ \\ 
\hline 
antiholom.\ invol.\ $\sigma$ & $\oline z$  & $-\oline z$ & 
$\beta i + \oline z$ \\ 
$\Omega^\sigma = \Fix(\sigma)$ & $(-1,1)$ & $i \R_+$ & 
$\frac{\beta i}{2} + \R$ 
\phantom{\Big)} \\ 
$(Rf)(x) = f(\sigma(x))$, $x \in \partial\Omega$ 
& $\oline z f(\oline z), z \in \T$ & $f(-x), x \in \R$ & 
$(f(\beta i + x), f(x)), x \in \R$ \\ 
\hline
Szeg\"o kernel $Q(z,w)$ & $\frac{1}{2\pi}\frac{1}{1 - z\oline w}$
 &$\frac{1}{2\pi} \frac{i}{z-\oline w}$  \phantom{$\Big)$}
                                 & 
$\frac{1}{4\beta\cosh(\frac{\pi}{2\beta}(z-\oline w)- \frac{\pi i}{2})}$ 
                                   \phantom{\Big)}\\  
Poisson kernel $P(z,x) = P_z(x)$ & $\frac{1}{2\pi}\frac{1-r^2}{1 - 2r \cos(\theta-t) + r^2}$ & 
$\frac{1}{\pi} \frac{\Im z}{|z-x|^2}$ & 
$\frac{1}{4\beta} \frac{\sin(\frac{\pi}{\beta} \Im z)}
{\sinh^2(\frac{\pi}{2\beta}(\Re z-x))
+ \sin^2(\frac{\pi}{2\beta}\Im z)}\phantom{\Bigg)}$ \\ 
\qquad $z \in \Omega$, $x \in \partial \Omega$ & $z = re^{i\theta}, x=e^{it}$ & $x\in \R$ & $P(z,x + i \beta) = P(i\beta + \oline z,x), x \in\R$\\
\hline
$*$-subsemigroup of $(H^\infty(\Omega),\sharp)$  & 
$c_n(z) = z^n,\ n \in \N_0$ & 
$c_t(z) = e^{itz},\ t \geq 0$ & 
$g_t(z) = e^{t\beta/2} e^{itz}$, $t \in \R$  \phantom{\Big)}\\
                & $c_n^\sharp = c_n$
                        & $c_t^\sharp = c_t$ & $g_t^\sharp = g_{-t}$
                                \phantom{\Big)}\\   
\hline
  boundary group $A$ & $\T$ (elliptic) & $\R$   (unipotent) &
$\R$ (hyperbolic)  \\
  symmetric (semi-)group $\hat A$ for RP & $\N_0 \subeq \Z$,
                        $n^\sharp = n$ & $\R_+\subeq \R$, $t^\sharp = t$ &
$\T_\beta\supeq \T_{\beta,+} = \{ [y] \: 0 \leq y \leq \beta/2\}$ \\
\hline
extremal RP functions  $\phi_\lambda$ 
& $\phi_\lambda(n) = \lambda^{|n|}, \lambda \in [-1,1]$ 
& $\phi_\lambda(t) = e^{-\lambda |t|}, \lambda \geq 0$ 
& $\phi_\lambda(t) = \frac{e^{-\lambda t} + e^{-\lambda(\beta-t)}}{1 + e^{-\lambda\beta}}, \lambda \geq 0$   \phantom{\Big)}\\
&&& on $\T_\beta = \R/\Z \beta$ \\ 
$\hat\phi_\lambda$ on $A$
& $P_\lambda(e^{it}) =\frac{1}{2\pi}\frac{1-\lambda^2}
{1 - 2\lambda\cos(t) + \lambda^2}$
& $P_{i\lambda}(x) = \frac{1}{\pi} \frac{\lambda}{\lambda^2 + x^2}$\phantom{\Big)} & 
  $P_{\frac{i\beta}{2} +\lambda}(x)
= \frac{1}{2\beta} \frac{1} 
  {\cosh(\frac{\pi}{\beta}(\lambda-x))}$
  \phantom{\bigg)}\\
                &\qquad \qquad \qquad $|\lambda| <1$ 
&  &
 \\
 outer fct $F_\lambda = \frac{Q_\lambda}{Q(\lambda,\lambda)^{1/2}}\in H^2(\Omega)$ &
$\frac{1}{\sqrt{2\pi}}\frac {\sqrt{1-\lambda^2}}{1-\lambda z}$
&
$\frac{Q_{i\lambda}(z)}{Q(i\lambda,i\lambda)^{1/2}}
= \frac{\sqrt\lambda}{\sqrt\pi} \frac{i}{i\lambda + z}$
&
$\frac{Q_{\frac{\beta i }{2} + \lambda}(z)}{Q\big(\frac{\beta i}{2} + \lambda,
                       \frac{\beta i}{2} + \lambda\big)^{1/2}}
 = \frac{1}{\cosh(\frac{\pi}{2\beta}(z - \frac{i\beta}{2} - \lambda))}.$
  \\
  
  with  $|F_\lambda^*|^2 = P_\lambda, \lambda \in \Omega^\sigma$
  &&& \\
  $h_\lambda = \frac{F_\lambda^*}{R F_\lambda^*}\phantom{\Big)}
= \frac{F_\lambda^*}{\oline{F_\lambda^*}}$
                & 
$ \frac{1-\lambda e^{-it}}{1-\lambda e^{it}}$ &
    $\frac{i\lambda -x}{i\lambda + x}$
&$\frac{Q_{\frac{\beta i}{2} + \lambda}(x)}{Q_{\frac{\beta i}{2} + \lambda}(\beta i + x)}=
\frac{\cosh(\frac{\pi}{2\beta}(x + \frac{\beta i}{2} - \lambda))}
{\cosh(\frac{\pi}{2\beta}(x - \frac{i\beta}{2} - \lambda))}$.
\phantom{\Bigg(}\\ 
\hline 
\end{tabular} \\[2mm]

\nin {\bf Notation:} 
\begin{itemize}
\item $\R_+ = (0,\infty)$, $\bD = \{ z \in \C \: |z| <1\}$, 
  $\C_+ = \R + i \R_+$ (upper half plane),\\
$\bS_\beta = \{ z \in \C \: 0 < \Im z < \beta\}$ (horizontal strip).
\item We write $\omega \: \bD \to \C_+,  \omega(z) := i \frac{1 + z}{1-z}$ 
for the Cayley transform with 
$\omega^{-1}(w) = \frac{w-i}{w+i}.$ 
\item On the circle $\T = \{ e^{i\theta} \in \C \:  \theta \in \R\}$, we 
use the length measure  of total length~$2\pi$. 
\item We write $e_w(z) := e^{zw}$ for $z,w\in \C$. 
\item For a function $f \:  G \to \C$ on the group $G$, we put 
$f^\vee(g) := f(g^{-1})$. 
\item For a holomorphic function $F$ on $\Omega \subeq \C$, we write 
  $F^*$ for its (non-tangential) boundary values.
\item For a holomorphic function $f$ on $\Omega$, we put
  $f^\sharp(z) = \oline{f(\sigma(z))}$.
\item We use the $L^2$-isometric 
  Fourier transform on the real line: 
  $\hat f(x) =\frac{1}{\sqrt{2\pi}} \int_\R f(y) e^{i xy}\, dy$
and \break $f(y) =\frac{1}{\sqrt{2\pi}} \int_\R \hat f(x) e^{-i xy}\, dx$.
  For a measure $\mu$ on $\R$, we use
  $\hat \mu(x) = \int_\R e^{i xp}\, d\mu(p)$.
\end{itemize}
\vspace{10mm} 

\nin \nin {\bf Acknowledgment:} 
MSA is a Humboldt Research Fellow and gratefully acknowledges support by the Alexander von Humboldt Foundation. Part of this work has been carried out while MSA was a JSPS International Research Fellow, receiving support by the Grant-in-Aid Kakenhi n. 22F21312, and during her stay at the Mathematisches Forschungsinstitut Oberwolfach (MFO) as a Oberwolfach Leibniz Fellow working on projects entitled ``Beurling-Lax type theorems and their connection with standard subspaces in Algebraic QFT'' and ``Reflection positive representations and standard subspaces in algebraic QFT''. MSA was supported by the University of Rome ``Tor Vergata'' funding scheme ``Beyond Borders'' CUP E84I19002200005, and is partially supported by GNAMPA--INdAM. MSA thanks Yoh Tanimoto for useful comments and suggestions. JS acknowledges the support by the UNAM Posdoctoral Program (POSDOC) (estancia posdoctoral realizada gracias al Programa de Becas Posdoctorales en la UNAM (POSDOC)) and the project CONACYT, FORDECYT-PRONACES 429825/2020 (proyecto apoyado por el FORDECYT-PRONACES, PRONACES/429825), recently renamed project CF-2019 /429825. KHN acknowledges support by DFG-grants NE 413/10-1 and NE 413/10-2.

\section{Reflection positive functions} 
\mlabel{sec:1}

This section introduces the main concepts and structures related to
reflection positivity on groups. 

A {\it reflection positive Hilbert space} is a triple $(\cE,\cE_+,\theta)$,  
consisting of a Hilbert space 
$\cE$ with a unitary involution $\theta$ and a closed subspace  
$\cE_+$ satisfying 
\begin{equation}
  \label{eq:thetapos}
\la \xi,\xi\ra_\theta := \la \xi, \theta \xi \ra \geq 0 
\quad \mbox{  for } \quad \xi \in \cE_+.
\end{equation}
This structure immediately leads to a new Hilbert space 
$\hat\cE$ that we obtain from the positive semidefinite form  
$\la \cdot,\cdot\ra_\theta$ on $\cE_+$ by completing the quotient of 
$\cE_+$ by the subspace of null vectors. 
We write $q \: \cE_+ \to \hat \cE, \xi \mapsto \hat\xi$ for the natural map
(\cite{NO18}).

\begin{defn} \mlabel{def:posdeffunc} 
Let $\tau$ be an involutive automorphism of a group $G$ 
and $G_+ \subeq G$ a subset invariant under the involutive map 
$g \mapsto g^\sharp := \tau(g)^{-1}$. Then a $\tau$-invariant 
 function $\phi \:  G \to \C$ is called {\it reflection positive} 
with respect to $(G,G_+,\tau)$ if $\phi$ 
is positive definite, i.e., the 
kernel $(\phi(gh^{-1}))_{g,h \in G}$ is positive definite, and, in addition,  
the kernel $(\phi(st^\sharp))_{s,t \in G_+}$ is also positive definite.

If $\mathbb{C}^G$ denotes all the $\mathbb{C}$-valued functions on $G$, $\varphi$ being reflection positive implies that the reproducing kernel space
$\cE := \cH_K \subeq \C^G$ with the kernel
$K(g,h) := \phi(gh^{-1})$ is reflection
positive with $\cE_+ = \oline{\Spann\{K_h \: h \in G_+\}}$ and
$\theta(f)(g) = f(g^\sharp)$ (cf.\ \cite[\S 3.4]{NO18}). 
\end{defn}

\begin{rem} \mlabel{rem:1.4}
  (a) (Reflection positive matrix coefficients)  
Suppose that $(\cE,\cE_+,\theta)$ is a 
reflection positive Hilbert space and 
$(U,\cE)$ a unitary representation of $G$ on $\cE$ satisfying 
\begin{equation}
  \label{eq:theta-cond}
\theta U(g) \theta = U(\tau(g)) \quad \mbox{  for } \quad g \in G. 
\end{equation}
If $v \in (\cE_+)^\theta = \{ \xi \in \cE_+ \: \theta \xi = \xi\}$ is such that 
$U(G_+)v \subeq \cE_+,$ then the function 
\begin{equation}
  \label{eq:repostate}
 \phi(g) := \la v, U(g) v \ra 
\end{equation}
is clearly positive definite on the group $G$. It is actually reflection positive 
with respect to $(G,G_+,\tau)$ because, for $s,t \in G_+$, we have  
\[ \phi(st^\sharp) 
= \la \theta v, U(st^\sharp) v \ra 
= \la v,\theta U(st^\sharp) v \ra 
= \la v, U(\tau(s)) \theta U(t^\sharp)v \ra 
= \la U(s^\sharp) v, \theta U(t^\sharp) v \ra, \] 
$U(G_+)v \subeq \cE_+$, $s^\sharp, t^\sharp \in G_+$, and $\cE_+$ is $\theta$-positive. 

\nin (b) Reflection positive functions are of particular interest if 
$G_+$ is a subsemigroup of~$G$. 
Accordingly, we call a triple $(G,S,\tau)$ a 
 {\it symmetric semigroup} if $G$ is a group and 
$S \subeq G$ is a subsemigroup satisfying $\tau(S)^{-1} = S$. 
Then $s^\sharp := \tau(s)^{-1}$ defines an involution on $S$. 
In this context one is interested in representations of 
$G$ on reflection positive Hilbert spaces $(\cE,\cE_+,\theta)$ 
satisfying $U(S) \cE_+ \subeq \cE_+$ and $\theta U(g)\theta=U(\tau(g))$ for all $g\in G$. Then any 
$v \in (\cE^+)^\theta$ defines by \eqref{eq:repostate} 
a positive definite function on $G$ which is also positive definite 
on the involutive semigroup $(S,\sharp)$. 
For symmetric 
semigroups $(G,S,\tau)$, the $\tau$-invariant 
reflection positive functions are all of the form 
\eqref{eq:repostate}. For more details on the 
GNS construction in the context of reflection positive functions, 
we refer to \cite[Thm.~3.4.5]{NO18}. 
\end{rem}

\begin{exs} The examples we shall encounter below are: 
  \begin{itemize}
  \item[\rm(a)] $(\Z,\N_0,-\id_\Z)$ a discrete symmetric semigroup. 
  \item[\rm(b)] $(\R,\R_+,-\id_\R)$ a one-dimensional symmetric semigroup. 
  \item[\rm(c)] For $\beta > 0$ the circle $\T_\beta := \R/\beta \Z$ 
of length~$\beta$ 
with the involution $\tau(z) = z^{-1}$ 
provides an important example where $G_+$ is not a semigroup, namely 
\[ \T_{\beta,+} := \{ t + \beta \Z \: 0 < t < \beta/2\}. \] 
  \end{itemize}
\end{exs}

\section{Reflection positivity on simply connected domains}
\mlabel{sec:3a}

In this section we briefly discuss some key structures
relating the Hardy space of a simply connected domain $\Omega$
the $*$-algebra $(H^\infty(\Omega),\sharp)$
and reflection positivity.
In particular, for a domain $\Omega$ with smooth boundary, Proposition~\ref{prop:2.3}
relates the weakly continuous characters of
$(H^\infty(\C_+), \sharp)$ with a reflection positive
Hilbert space for which $(\cE,\cE_+) = (L^2(\partial \Omega), H^2(\Omega))$.

Let $\Omega \subeq \C$ be a proper simply connected domain 
and $\psi \: \Omega \to \bD$ be biholomorphic. The existence of such maps follows 
from the Riemann Mapping Theorem. 
Then we define the Hardy space $H^2(\Omega)$ in such a way that the map 
\begin{equation}
  \label{eq:unih2a}
\Gamma_\psi \: H^2(\bD) \to H^2(\Omega), \quad 
\Gamma_\psi(f) := \sqrt{\psi'} \cdot f \circ \psi 
\end{equation}
is unitary. This is well-defined because the group 
$\Aut(\bD)$ of biholomorphic automorphisms of $\bD$ acts unitarily 
on the Hilbert space $H^2(\bD)$ by \eqref{eq:unih2a}. 

If $\phi \: \Omega_1 \to \Omega_2$ is a biholomorphic map between 
two simply connected proper domains in $\C$, then the map 
\begin{equation}
  \label{eq:unih2}
\Gamma_\phi \: H^2(\Omega_2) \to H^2(\Omega_1), \quad 
\Gamma_\phi(f) := \sqrt{\phi'} \cdot f \circ \phi 
\end{equation}
is unitary up to a positive factor, depending on the normalization 
of the scalar product. 

In particular, we have the following maps 
(see {\rm\cite[p.~200]{Ni19}} for the isomorphisms between $H^2(\bD)$ and $H^2(\C_+)$):
\begin{equation}
  \label{eq:dom-traf}
\begin{tikzcd}[column sep=7pc,row sep=2pc]
\bD \arrow[bend left]{r}{\omega(z)= i \frac{1+z}{1-z}} & \C_+ \arrow[bend left]{l}{\omega^{-1}(z) = \frac{z-i}{z+i}} \arrow[bend left]{r}{\Exp^{-1}(z) =\frac \beta \pi \log\left(z\right)} & \bS_\beta \arrow[bend left]{l}{\Exp(z) 
= \exp(\frac {\pi z} \beta)}
\end{tikzcd}
\end{equation}
These maps induce the following diagram of unitary bijections:
\begin{equation}
  \label{eq:func-traf}
\begin{tikzcd}[column sep=11pc,row sep=2pc]
H^2(\bD) \arrow[bend left]{r}{\Gamma(f)(z) 
= \frac{\sqrt{2i}}{z+i}
f\left(\frac{z-i}{z+i}\right)} 
& H^2(\C_+) \arrow[bend left]{l}{\Gamma^{-1}(f)(z) 
= \frac{\sqrt{2 i}}{1-z}f\left(i \frac{1+z}{1-z}\right)}
\arrow[bend left]{r}{\Phi(f)(z) 
= \sqrt{\frac \pi \beta}\exp\left(\frac {\pi z} {2\beta}\right) 
f\left(\exp\left(\frac {\pi z} \beta\right)\right)}
& H^2(\bS_\beta) \arrow[bend left]{l}{\Phi^{-1}(f)(z) 
= \sqrt{\frac \beta \pi} \frac 1 {\sqrt{z}}f\left(\frac \beta 
\pi \log\left(z\right)\right)}
\end{tikzcd}
\end{equation}

The Hardy space $H^2(\Omega)$ is a reproducing kernel 
Hilbert space, i.e., the point evaluations 
\[  \ev_z \: H^2(\Omega) \to \C, \quad f \mapsto f(z) \]
are continuous linear functionals, hence can be written as 
\[ f(z) = \la Q_z, f \ra \quad \mbox{ for some } \quad Q_z \in H^2(\Omega).\] 
The kernel 
\[ Q \: \Omega \times \Omega \to \C, \quad Q(z,w) := Q_w(z) = \la Q_z, Q_w \ra \] 
is called the {\it Szeg\"o kernel of $\Omega$}. 
If $\Omega$ has smooth boundary (which we assume throughout), we have an isometric 
boundary value map 
\[ H^2(\Omega) \to L^2(\partial \Omega), \quad 
f \mapsto f^*,\] 
then we obtain  the {\it Poisson kernel} of $\Omega$ by 
the {\it Hua formula} (cf.\ \cite{Hu63}, \cite{Ko65}) 
\begin{equation}
  \label{eq:poiss-omega}
P \: \Omega \times \partial \Omega  \to \R, \quad 
 P(z,x) = P_z(x) = \frac{|Q(z,x)|^2}{Q(z,z)} \quad \mbox{ for } \quad 
z \in \Omega, x \in \partial\Omega,
\end{equation}
i.e., 
\[ P_z := P(z,\cdot) = \frac{|Q_z^*|^2}{Q(z,z)} \in L^1(\partial \Omega).\]  
We then have the {\it Poisson formula} 
\begin{equation}
  \label{eq:poissform}
 f(z) = \int_{\partial \Omega} P_z(x) f^*(x)\, dx \quad \mbox{ for } \quad 
f \in H^\infty(\Omega), z \in \Omega.
\end{equation}

\begin{defn} We call a function $F \in H^2(\Omega)$ 
{\it outer} if it is cyclic for the representation 
of the Banach algebra $H^\infty(\Omega)$ on $H^2(\Omega)$ 
(see 
\cite[Thm.~12, Thm.~22]{CHLN16}).
\end{defn}

\begin{lem} \mlabel{lem:a.4} 
The Szeg\"o functions $(Q_w)_{w \in \Omega}$ are outer. 
\end{lem}

\begin{prf} In Appendix B.6 of \cite{ANS22} we have seen that,  
for a biholomorphic map $\phi \:  \Omega_1 \to \Omega_2$ 
between simply connected domains, the Szeg\"o functions transform by 
\begin{equation}
  \label{eq:sz-trafo}
Q^{\Omega_1}_w 
= \oline{\sqrt{\phi'(w)}} \Gamma_\phi(Q^{\Omega_2}_{\phi(w)}) 
\quad \mbox{ for } \quad w \in \Omega_1.
\end{equation}
Further, \eqref{eq:unih2} shows that 
$F \in H^2(\Omega_2)$ is outer if and only if 
$\Gamma_\phi(F)$ is outer in $H^2(\Omega_1)$. 
In view of the Riemann Mapping Theorem, it therefore 
suffices to prove the lemma for the unit 
disc~$\bD$. 

Then $2\pi Q_w(z) = \frac{1}{1 - z \oline w}$. 
For $|z| = 1$, we have 
\[\frac{1}{1+|w|}  \leq 2\pi |Q_w(z)| \leq \frac{1}{1-|w|}.\] 
Therefore $Q_w, Q_w^{-1} \in H^\infty(\bD)$, and now 
$Q_w$  is outer by \cite[Lemma~B.13]{ANS22}.
\end{prf}

Hua's formula \eqref{eq:poiss-omega} now shows that, for every 
$z \in \Omega$, the function 
$F_z := \frac{Q_z}{\sqrt{Q(z,z)}}$ is outer with 
\begin{equation}
  \label{eq:pzout}
 |F_z^*|^2 = P_z 
\end{equation}

Let $\sigma \:  \Omega \to \Omega$ 
be an antiholomorphic involution $\sigma$ on $\Omega$. 
Then 
\[ f^\sharp(z) := \oline{f(\sigma(z))}\quad \mbox{ for } \quad 
z \in \Omega\] 
defines the structure of a 
Banach $*$-algebra on the Banach algebra $H^\infty(\Omega)$ (\cite[App.~B]{ANS22}).  
As the Szeg\"o kernel satisfies 
\begin{equation}
  \label{eq:csprop}
   Q(\sigma z, \sigma w ) = \oline{Q(z,w)} = Q(w,z),
\end{equation}
we have for $Q_w := Q(\cdot, w)$ in $H^2(\Omega)$ the relation 
\begin{equation}
  \label{eq:qwsharp}
 Q_w^\sharp = Q_{\sigma(w)} 
\end{equation}
and in particular $Q_w^\sharp = Q_w$ for $w \in \Omega^\sigma$.
It also follows that $\sharp$ defines a conjugation
(an antilinear involutive isometry) on $H^2(\Omega)$.  

On the boundary space $L^2(\partial \Omega)$, we consider the unitary involution 
\[ (Rf)(x) := f(\sigma(x))\quad\text{for}\quad f\in L^2(\partial\Omega) \] 
and observe that 
\[ RQ_w^* = \oline{Q^*_{\sigma(w)}} = \oline{Q_w^*} \quad \mbox{ for } \quad 
w \in \Omega^\sigma.\] 
Therefore the measurable function 
\begin{equation}
  \label{eq:hw}
 h_w := \frac{Q_w^*}{R Q_w^*} = \frac{Q_w^*}{\oline{Q_w^*}} 
\end{equation}
is unimodular on $\partial \Omega$. It satisfies 
\begin{equation}
  \label{eq:rhw}
 (R h_w)(x) = \oline{h_w(x)} = h_w(x)^{-1} 
\end{equation}
almost everywhere on $\partial \Omega$. Therefore 
\begin{equation}
  \label{eq:thetaw}
  \theta_w := h_w R, \quad \mbox{ i.e.} \quad
  \theta_w(f) = h_w \cdot (f \circ \sigma), 
\end{equation}
defines a unitary involution on $L^2(\partial \Omega)$.
By construction, it satisfies 
\begin{equation}
  \label{eq:thetawq}
  \theta_w Q_w^* = Q_w^*.
\end{equation}
For $F \in H^\infty(\Omega)$ and $f \in H^2(\Omega)$,  we have 
\[  \theta_w(F^* f^*)(x) 
= h_w(x) F^*(\sigma(x)) f^*(\sigma(x))
= F^*(\sigma(x)) \theta_w(f^*)(x).\] 
Writing $F^* \in L^\infty(\partial \Omega)$ for the corresponding 
multiplication operator, this leads to 
\begin{equation}
  \label{eq:thetaw2}
  \theta_w \circ F^* = \oline{F^{\sharp,*}} \circ \theta_w,
\end{equation}

\begin{prop} \mlabel{prop:2.3} 
For any $w \in \Omega^\sigma$, 
the triple 
$(L^2(\partial \Omega), H^2(\Omega), \theta_w)$ is a reflection positive 
Hilbert space with 
\begin{equation}
  \label{eq:flipf1}
 \la f^*, \theta_w f^* \ra = \frac{|f(w)|^2}{Q(w,w)} 
\quad \mbox{ for } \quad f \in H^2(\Omega).
\end{equation}
\end{prop}

\begin{prf} We  first evaluate for $F \in H^\infty(\Omega)$: 
\begin{align*} \label{eq:calc} 
\la F^* Q_w^*, \theta_w F^* Q_w^* \ra 
&=   \la F^* Q_w^*, \oline{F^{\sharp,*}} \theta_w Q_w^* \ra 
=   \la F^* Q_w^*, \oline{F^{\sharp,*}} Q_w^* \ra 
=   \la F^{\sharp,*}  F^* Q_w^*, Q_w^* \ra \\
&=   \la F^\sharp F Q_w, Q_w \ra_{H^2}  
=   (F^\sharp F Q_w)(w) = \oline{F(\sigma(w))} F(w) Q(w,w) \\
&= \oline{F(w)} F(w) Q(w,w) = |F(w)|^2 Q(w,w) = \frac{|(F Q_w)(w)|^2}{Q(w,w)} \geq 0.
\end{align*}
As $Q_w \in H^2(\Omega)$ is outer
by Lemma~\ref{lem:a.4}, the subspace 
$H^\infty(\Omega) Q_w$ is dense in $H^2(\Omega)$, and therefore 
the above calculation proves the proposition. 
\end{prf}

From \eqref{eq:flipf1} we also derive 
that $\la f^*, \theta_w f^*\ra$  vanishes for $f(w) = 0$ and that 
\[ \la Q_w^*, \theta_w Q_w^*\ra = Q(w,w) = \|Q_w\|^2.\] 
Therefore $\hat{\cE_+} \cong \C$ 
and the corresponding quotient map is given by evaluation in $w$: 
\[ q_w(f) = f(w) = \la Q_w, f \ra.\] 

\begin{rem} Polarizing the formula 
\[ \la F^* Q_w^*, \theta_w F^* Q_w^* \ra 
= \oline{F(w)} F(w) Q(w,w)
\quad \mbox{ for } \quad w \in \Omega^\sigma, F \in H^\infty(\Omega),\]
from the proof of Proposition~\ref{prop:2.3}, we obtain 
\[ \la F^* Q_w^*, \theta_w G^* Q_w^* \ra 
= \oline{F(w)} G(w) Q(w,w) \quad \mbox{ for } \quad F, G \in H^\infty(\Omega).\]
This corresponds to the positive functional $\delta_w(F) = F(w)$ on 
the commutative Banach \break {$*$-algebra} $(H^\infty(\Omega),  \sharp)$ 
(\cite[Prop.~B.8]{ANS22}), for which the GNS construction produces  the
one-dimensional representation given by evaluation in~$w$. 
\end{rem}

We now turn to general weakly continuous
  positive functionals on the Banach $*$-algebra $(H^\infty(\Omega), \sharp)$.

\begin{prop}\label{pos_funct} For any positive functional $\eta$ on
  the Banach $*$-algebra $(H^\infty(\Omega),\sharp)$, 
  there exists a finite positive Borel measure $\mu$ on
  $\Omega^\sigma$ and   an outer function \(F \in H^2(\Omega)^\sharp\)  such that
the unitary involution on $L^2(\partial \Omega)$, 
  defined by
  \begin{equation}
    \label{eq:thetaf}
    \theta_F(f) := \frac{F^*}{F^*\circ \sigma} \cdot (f \circ \sigma),\quad\text{for}\quad f\in L^2(\partial\Omega)
  \end{equation}
satisfies
\[\eta(G) = \bracket{F^*}{\theta_F G^*F^*}_{L^2(\partial \Omega)}
  = \bracket{F^*}{G^*F^*}_{L^2(\partial \Omega)}
  = \int_{\Omega^\sigma} G(\lambda)\, d\mu(\lambda) 
  \quad \mbox{ for } \quad G \in H^\infty(\Omega).\]
\end{prop}

\begin{prf} In view of \cite[Lemma~B.1]{ANS22},
  all pairs $(\Omega, \sigma)$ are conjugate under biholomorphic maps.
  We may therefore assume w.l.o.g.\ that \(\Omega = \C_+\)
  and $\sigma(z) = - \oline z$.
  
  First we use \cite[Prop. B.8]{ANS22} to find a finite 
positive Borel measure \(\mu\) on \(\R_+\) such that
\[ \eta(G) = \eta_\mu(G) := \int_0^{\infty} G(i\lambda)\, d\mu(\lambda)
  \quad \mbox{ for } \quad
  G \in H^\infty(\C_+).\]
Then the function
\[ \Psi(p)
  = \frac{1}{\pi}
  \int_0^\infty  \frac{\lambda}{p^2 + \lambda^2}\, d\mu(\lambda) \] 
is $L^1$ and specifies an outer function
\[ F(z) := \Out(\Psi^{1/2})(z)
  := \exp\Big(\frac{1}{2\pi i} \int_\R \Big[ \frac{1}{p-z}
  - \frac{p}{1 + p^2}\Big] \log(\Psi(p))\, dp\Big)\]
in $H^2(\C_+)$ satisfying $F(-\oline z) = \oline{F(z)}$ for $z \in \C_+$, i.e.,
$F = F^\sharp$ 
(see \cite[Lemma~4.2.18]{Sch23} for the existence of the integral
and \cite[Lemma~A.2.6]{Sch23}). The function $F$ is uniquely
  determined up to a multiplicative constant in $\T$ by having the 
  boundary values $|F^*(x)| = |\Out(\Psi^{1/2})(x)|$ for $x \in \R$. 
Now \cite[Cor.~4.3.14]{Sch23} implies that
\[ \la f^*, \theta_F g^* \ra
  = \int_{\R_+} \oline{f(i\lambda)} g(i\lambda)\,
  \frac{ d\mu(\lambda)}{|F(i\lambda)|^2}
  \quad \mbox{ for } \quad f,g \in H^2(\C_+),\]
so that, by $\theta_F(F^*)=F^*$, 
\[ \la F^*, f^* F^* \ra 
 = \la \theta_F F^*, f^* F^* \ra   =\la F^*, \theta_F f^* F^* \ra 
  =  \int_{\R_+} f(i\lambda)\, d\mu(\lambda)=\eta_\mu(f).\qedhere\]
\end{prf}

\begin{prop}
  Let \(\nu\) be a Borel measure on \(\Omega^\sigma\) such that 
  the restriction map $f \mapsto f\res_{\Omega^\sigma}$ defines a
  continuous linear map $r : H^2(\Omega) \to L^2(\Omega^\sigma,\nu)$
  and put $H_\nu := r^* r \in B(H^2(\Omega))$.
  \begin{itemize}
  \item[\rm(a)] If \(\left\lVert H_\nu\right\rVert = \|r\|^2\leq 1\),
then there exists a unitary involution \(\theta \in B(L^2(\partial \Omega))\) such that the triple
  \((L^2(\partial \Omega),H^2(\Omega),\theta)\) is a reflection positive Hilbert space and 
  \[  H_\nu  = P_{H^2(\Omega)} \theta P^*_{H^2(\Omega)}   \]
  holds for the projection $P_{H^2(\Omega)} \: L^2(\partial \Omega) \to
  H^2(\Omega)$ onto the Hardy space. 
\item[\rm(b)] If, in addition, \(\Fix(H_\nu) \neq \{0\}\),
  then 
  any nonzero $F \in \Fix(H_\nu)$ is outer, 
the involution \(\theta = \theta_F\) is unique, and 
   \[ \eta_\nu \: H^\infty(\Omega) \to \C, \quad
     \eta_\nu(f) := \bracket{F^*}{\theta_F f^*F^*} = \bracket{F^*}{f^*F^*}
   = \int_{\Omega^\nu} f(\omega) |F(\omega)|^2 \, d\nu(\omega)\]
defines a positive functional on the Banach
$*$-algebra $(H^\infty(\Omega),\sharp)$. 
\end{itemize}
\end{prop}

\begin{prf} In view of \cite[Lemma~B.1]{ANS22},
  all pairs $(\Omega, \sigma)$ are conjugate under biholomorphic maps.
  We may therefore assume w.l.o.g.\ that \(\Omega = \C_+\)
  and $\sigma(z) = - \oline z$.  
  Assertion (a) now follows from \cite[Thm.~5.3.5]{Sch23} and
  (b) from \cite[Prop.~5.4.3]{Sch23} and with Proposition \ref{pos_funct}
  from 
\[\eta(f^\sharp f)=\bracket{F^*}{\theta_F f^{\sharp,*} f^* F^*} = \bracket{f^*F^*}{\theta_F f^* F^*} \geq 0 \qquad
  \mbox{ for } \quad f \in H^\infty(\C_+).\qedhere\]
\end{prf}

\begin{remark} The situation of Proposition~\ref{prop:2.3}
  corresponds to the point measure \(\nu = \frac 1{Q(w,w)}\delta_w\), which satisfies
  \[|\bracket{f}{H_\nu g}|
    = \frac{|\overline{f(w)}g(w)|}{Q(w,w)}
    = \frac{|\bracket{f}{Q_w}| \cdot |\bracket{Q_w}{g}|}
    {\bracket{Q_w}{Q_w}} \leq \frac{\left\lVert Q_w\right\rVert^2 \cdot \left\lVert f\right\rVert \cdot \left\lVert g\right\rVert}{\left\lVert Q_w\right\rVert^2} = \left\lVert f\right\rVert \cdot \left\lVert g\right\rVert \]
  for $f,g \in H^2(\Omega)$,
so that \(\left\lVert H_\nu\right\rVert \leq 1\). Also
\[\bracket{f}{H_\nu Q_w} = \frac{\overline{f(w)}Q_w(w)}{Q(w,w)} = \overline{f(w)} = \bracket{f}{Q_w} \qquad
  \mbox{ for } \quad f \in H^2(\Omega),\]
so that \(H_\nu Q_w = Q_w\) and therefore \(Q_w \in \Fix(H_\nu) \setminus \{0\}\).
\end{remark}

\section{Reflection positivity on the upper half-plane} 
\mlabel{sec:3}

On the upper half-plane, the group $\R$ acts by translations and the 
$\R$-eigenfunctions in $H^\infty(\C_+)$ are multiples of the 
functions $e_{it}$, $t \geq 0$. 
Here we consider the antiholomorphic involution $\sigma(z) = -\oline z$ fixing 
$i \R_+$ pointwise. As they are real in $i\R_+$, the functions 
$e_{it}$ for $t \geq 0$ satisfy $e_{it}^\sharp = e_{it}$, so that we obtain a 
homomorphism of $*$-semigroups 
\[ (\R_+, \id) \to (H^\infty(\C_+),\sharp), \quad 
t \mapsto e_{it}. \] 

The Hardy space $H^2(\C_+) \subeq L^2(\R)$ is 
the positive spectral part for the translation action 
of $\R$ on $L^2(\R)$ and 
\[ (U_t f)(x) = e^{itx} f(x) = (e_{it}f)(x), \quad t \in \R,\] defines a unitary 
one-parameter group on $L^2(\R)$ for which 
$H^2(\C_+)$ is invariant under $(U_t)_{t \geq 0}$. 

The {\it Szeg\"o kernel} on the upper half-plane is  given by 
\begin{equation}
  \label{eq:cauchyker}
 Q(z,w) = \frac{1}{2\pi} \frac{i}{z - \oline w} 
\quad \mbox{ for } \quad z,w \in \C_+.
\end{equation}
This is an easy consequence of the Residue Theorem, which implies that 
\[ \la Q_z, F \ra = \frac{1}{2\pi i} \int_\R \frac{F^*(x)}{x -z}\, dx = F(z)  
\quad \mbox{ for } \quad 
F \in H^2(\C_+), z \in \C_+\]
(see also  Appendix~\ref{subsec:app-cz-kernel}). 
For the Poisson kernel we obtain with 
Hua's formula \eqref{eq:poiss-omega}
\begin{equation}
  \label{eq:poissker}
 P(z,x) = P_z(x) = \frac{1}{\pi} \frac{\Im z}{|z-x|^2} 
\quad \mbox{ and } \quad 
P_{i\lambda}(x) = \frac{1}{2\pi} \frac{\lambda}{\lambda^2 + x^2}  
\quad \mbox{ for } \quad \lambda > 0.
\end{equation}

General reflection positive functions on $(\R,\R_+,-\id_\R)$ have the integral 
representation 
\begin{equation}
  \label{eq:inrep}
 \phi(x) = \int_{[0,\infty)} e^{-\lambda|x|} \, d\mu(\lambda),
\end{equation}
where $\mu$ is a finite Borel measure on $[0,\infty)$ 
(\cite[Cor.~3.3]{NO14}).
\begin{footnote}{This can also 
be derived from \cite[Cor.~4.4.5]{BCR84}, resp., 
the abstract integral representation of bounded positive definite 
functions on abelian $*$-semigroup \cite[Thm.~4.2.8]{BCR84}.}
\end{footnote}
We now relate this class of functions to the complex geometry
of the upper half plane.

\begin{prop} \mlabel{prop:3.1} 
  For $\lambda > 0$, the following assertions hold:
  \begin{itemize}
\item[\rm(a)] The Fourier transform of the the positive function\footnote{Here and throughout the paper, we will deal with the Poisson kernel as a measure.} $P_{i\lambda}$ 
is the reflection positive function on $(\R,\R_+, - \id_\R)$, given by  
\begin{equation}
  \label{eq:philambdac+}
 \phi_\lambda(t) := \hat{P_{i\lambda}}(t) =  e^{-\lambda|t|} 
\quad \mbox{ for } \quad t\in \R.
\end{equation}
\item[\rm(b)]  The function
\begin{equation}
  \label{eq:Flambda}
 F_\lambda(z) :=  F_\lambda^{\C_+}(z) := \frac{Q(z,i\lambda)}{\sqrt{Q(i\lambda,i\lambda)}} 
=  \sqrt{\frac{\lambda}{\pi}} \frac{1}{z + i\lambda}
\end{equation}
is outer and its boundary values $F_\lambda^*$ on $\R$
satisfy $|F_\lambda^*|^2 = P_{i\lambda}$.
With 
\begin{equation}
  \label{eq:hlambda}
 h_\lambda(x) =  \frac{F_\lambda^*(x)}{F_\lambda^*(-x)} 
= \frac{i\lambda -x}{i\lambda + x}\in \T
\end{equation}
we then obtain on $L^2(\R)$ a unitary involution
\begin{equation}
  \label{eq:theta-inv}
 (\theta_{h_\lambda}f)(x) := h_\lambda(x) f(-x) 
\end{equation}
fixing $F_\lambda^*$, for which 
\begin{equation}
  \label{eq:28}
 \la F_\lambda^{\C_+}, U_t F_\lambda^{\C_+} \ra  =  \phi_\lambda(t)
 \quad \mbox{ for }  \quad t \in \R.
\end{equation}
  \end{itemize}
\end{prop}

Note that (b) implies that the GNS construction, applied
to the positive definite function $\phi_\lambda$ recovers the
representation $(U, L^2(\R))$, that
$(L^2(\R), H^2(\C_+), \theta_{h_\lambda})$ is reflection positive,
and, since the outer function
$F_\lambda^{\C_+}$ is cyclic in $H^2(\C_+)$ for the one-parameter
group $(U_t)_{t \geq  0}$, that the associated space
$\hat{H^2(\C_+)}$ is the one-dimensional representation
of $\R_+$, given by the character $t \mapsto e^{-t\lambda}$.

\begin{prf} (a) First we observe that 
\[  \hat{P_{i\lambda}}(t) 
= \int_\R e_{it}(x) P_{i\lambda}(x)\, dx 
= e_{it}(i\lambda) = e^{-t\lambda} \quad \mbox{ for } \quad \lambda, t > 0. \]
As $P_{i\lambda}$ is symmetric, we thus obtain \eqref{eq:philambdac+}. 
That this function is positive definite follows from the positivity of $P_{i\lambda}$ 
and that it is reflection positive follows from the fact that 
$\phi_\lambda(t) = e^{-\lambda t}$ is a $*$-character of $(\R_+,\id)$. 

\nin (b) 
By Hua's formula, the $L^1$-function $P_{i\lambda}$ 
coincides with the boundary values of $|F_\lambda^*|^2$ for  the outer function 
$F_\lambda^{\C_+} \in H^2(\C_+)$ (cf.\ Lemma~\ref{lem:a.4}).
Further, the relation $h_\lambda(-x) = h_\lambda(x)^{-1}$ implies that 
$\theta_{h_\lambda}$ defines a unitary involution on $L^2(\R)$.
It clearly fixes $F_\lambda^*$. Finally we realize the reflection positive
function $\phi_\lambda$ as the corresponding matrix coefficient:
\[ \la F_\lambda^{\C_+}, U_t F_\lambda^{\C_+} \ra 
= \int_\R |F_\lambda^*(x)|^2 e^{itx}\, dx = \int_\R P_{i\lambda}(x) e^{itx}\, dx 
\ {\buildrel \eqref{eq:philambdac+}\over =}\ \phi_\lambda(t)
\quad \mbox{ for } \quad t \in \R.\qedhere \] 
\end{prf}

\section{Reflection positivity on the disc} 
\mlabel{sec:2}

On the disc, the circle group $\T$ acts by rotations and the 
$\T$-eigenfunctions in $H^\infty(\bD)$ are multiples of the 
functions $e_n(z) := z^n$, $n \in \N_0$. 
We consider the antiholomorphic involution $\sigma(z) = \oline z$ on $\bD$ 
and note that the functions $e_n$ are invariant under the involution 
\[ f^\sharp(z) := \oline{f(\oline z)} \] 
on $H^\infty(\bD)$, so that 
\[ (\N_0, \id) \to (H^\infty(\bD), \sharp), \quad n \mapsto e_n, \] 
is a $*$-homomorphism. 
We identify $H^2(\bD) \subeq L^2(\T)$ with the 
positive spectral part for the $\T$-action on $L^2(\T)$. 
In this context 
\[ (U_n f)(z) = z^n f(z) \quad \mbox{ for } \quad z \in \T,\]
 defines a unitary representation of the pair $(G,S) = (\Z,\N_0)$ 
on $(L^2(\T), H^2(\bD))$. 

The Szeg\"o kernel of the disc is 
\[ Q(z,w) = Q_w(z) = \frac{1}{2\pi} \frac{1}{1 - z \oline w}.\] 
For $f \in H^2(\bD)$, we have 
\[  f(z) = \la Q_z, f \ra 
=  \int_0^{2\pi} \oline{Q_z(e^{i\theta})} f^*(e^{i\theta})\, d\theta 
= \frac{1}{2\pi } \int_{\partial \bD} \frac{f^*(\zeta)}{1 - \oline\zeta z} \, 
\frac{d\zeta}{i\zeta} 
= \frac{1}{2\pi i} \int_{\partial \bD} \frac{f^*(\zeta)}{\zeta - z} \, d\zeta.\]
We thus obtain from Hua's formula \eqref{eq:poiss-omega} 
the Poisson kernel 
\[ P(re^{i\theta}, e^{it}) 
= \frac{1}{2\pi}\frac{1-r^2}{|1 - re^{i(\theta-t)}|^2} 
= \frac{1}{2\pi}\frac{1-r^2}{1 - 2r \cos(\theta-t) + r^2}
\quad \mbox{ for } \quad z = r e^{i\theta}\in \bD, t \in [0,2\pi].\] 
For the circle this leads for each 
point $z = r e^{i\theta}\in \bD$ to a probability 
measure 
\[ P_z(e^{it}) \, dt  = P(z,e^{it}) \,  dt 
=  \frac{1-r^2}{1 - 2r\cos(\theta-t) + r^2}\, \frac{dt}{2\pi}, \] 
and for $\lambda \in (-1,1)$ we have in particular 
\[ P_\lambda(e^{it}) 
= \frac{1}{2\pi}\cdot \frac{1-\lambda^2}{1 - 2\lambda\cos(t) + \lambda^2}.\] 

Recall that a 
function $ \phi \: \Z \to \C$ is {\it reflection positive} 
for $(\Z,\N,-\id_\Z)$ 
if it is positive definite and the kernel 
$(\phi(n+m))_{n,m \geq 0}$ is positive definite (Definition~\ref{def:posdeffunc}).

The Fourier transform of $P_\lambda$ is the positive definite 
function on $\Z$ given by 
\[ \hat{P_\lambda}(n) = \int_0^{2\pi} 
e^{int}  \frac{1-\lambda^2}{1 - 2\lambda\cos(t) + \lambda^2}\, \frac{dt}{2\pi}
= \int_0^{2\pi} 
e^{int}  P(\lambda, e^{it})\, dt 
= \lambda^n \quad \mbox{ for } \quad n \in \N_0. \] 
As $P_\lambda$ is symmetric, we obtain 
\begin{equation}
  \label{eq:disc-philambda}
\phi_\lambda(n) = 
 \hat{P_\lambda}(n)  = \lambda^{|n|} \quad \mbox{ for } \quad 
n \in \Z, \lambda \in (-1,1), 
\end{equation}
which is a positive definite function on $\Z$. 
To see that it is actually reflection positive function on $(\Z,\N,-\id_\Z)$, 
it suffices to observe that its restriction to $\N_0$ is a $*$-homomorphism.

The function 
\[ F_\lambda^\bD(z) := \frac{Q(z,\lambda)}{\sqrt{Q(\lambda,\lambda)}} 
= \frac{1}{\sqrt{2\pi}}\frac {\sqrt{1-\lambda^2}}{1-\lambda z} \] 
is outer in $H^2(\bD)$ (Lemma~\ref{lem:a.4}) 
and by Hua's formula 
$\big|F_\lambda^{\bD,*}\big|^2 = P_\lambda$. 
With the unitary involution $(Rf)(z) =  f(\oline z)$ on $L^2(\T)$, 
we obtain the unimodular function
\[h^{\bD}_\lambda\left(e^{it}\right) := 
\frac{F_\lambda^{\bD,*}\left(e^{it}\right)}{R F_\lambda^{\bD,*}\left(e^{it}\right)} 
= \frac{F_\lambda^{\bD,*}\left(e^{it}\right)}{F_\lambda^{\bD,*}\left(e^{-it}\right)} 
= \frac{1-\lambda e^{-it}}{1-\lambda e^{it}}.\]
The relation $h_\lambda^\bD(z^{-1}) = h_\lambda^\bD(z)^{-1}$ implies that 
\begin{equation}
  \label{eq:theta-invb}
 (\theta_{h^\bD_\lambda}f)(z) := h^\bD_\lambda(z) f(z^{-1}) 
\end{equation}
defines a unitary involution on $L^2(\T)$ fixing $F_\lambda^{\bD,*}$. Note that 
\[ \la F_\lambda^\bD, U_n F^\bD_\lambda \ra 
= \int_\T z^n |F_\lambda^{\bD,*}(z)|^2\, dz 
= \int_\T P_{\lambda}(z) z^n\, dz = \lambda^n
\quad \mbox{ for } \quad n \in \N_0\] 
provides a realization of $\phi_\lambda$ as a matrix coefficient 
as in Remark~\ref{rem:1.4}. 

We conclude our discussion of the disc with a description
of all reflection positive functions on $(\Z,\N_0, -\id)$
in terms of an integral representation. 

\begin{thm} \mlabel{thm:3.2}
{\rm(Reflection positive functions on $\Z$)} 
A function $\phi  \: \Z \to \C$ is reflection positive if and only if  
there exists a finite positive Borel measure $\mu$ on $[-1,1]$ with 
\begin{equation}
  \label{eq:intrep}
\phi(n) = \int_{-1}^1 \lambda^{|n|}\, d\mu(\lambda) \quad \mbox{ for } \quad n \in \Z.
\end{equation}
\end{thm}

\begin{prf} For the point measures in the boundary, we obtain the 
characters 
\[ \chi_+(n) := 1 \quad \mbox{ and } \quad 
\chi_-(n) = (-1)^n,\quad n \in \Z, \]
which are both reflection positive. 
For $|\lambda| < 1$, the function 
$n \mapsto \lambda^{|n|}$ 
is reflection positive by~\eqref{eq:disc-philambda}. 
Therefore every bounded positive Borel measure 
$\mu$ on $[-1,1]$ defines by  \eqref{eq:intrep} a 
reflection positive function for $(\Z,\N_0,-\id_\Z)$. 

If, conversely, $\phi \: \Z \to \C$ is reflection positive,
then it is in particular positive definite, hence bounded. 
Further $\phi\res_{\N_0}$ is positive definite on the involutive 
semigroup $(\N_0, \id)$, so that there exists a 
uniquely determined positive Borel measure $\mu$ on 
$[-1,1]$ with 
\[ \phi(n) = \int_{[-1,1]} x^n\, d\mu(x) 
\quad \mbox{ for } \quad n \in \N_0 \] 
(Hamburger's Theorem; \cite[Thm.~6.2.2]{BCR84}). For $n < 0$ we now have 
$\phi(n) = \oline{\phi(-n)} = \phi(-n)$, 
and this proves \eqref{eq:intrep} for every $n \in \Z$. 
\end{prf}

\section{Reflection positivity on the strip} 
\mlabel{sec:4}

On the strip $\bS_\beta = \{ z \in \C \: 0 < \Im z < \beta\}$,
the group $\R$ acts by translations. It is a
hyperbolic one-parameter group of
$\Aut(\bS_\beta)_e \cong \PSL_2(\R)$. We
consider the  antiholomorphic involution 
\[ \sigma(z) = \beta i + \oline z.\] 
It defines the involution $\sharp$  on $H^\infty(\bS_\beta)$ 
by $f^{\sharp}(z) := \oline{f(\sigma(z))}$, turning it into a
Banach $*$-algebra (cf.\ \cite[App.~B]{ANS22}). 
The $\R$-eigenfunctions for the translation action on the
Banach $*$-algebra $H^\infty(\bS_\beta)$ are multiples of the 
functions $(e_{it})_{t \in \R}$. 
With respect to the involution $\sharp$
on $H^\infty(\bS_\beta)$, they satisfy 
\begin{equation}
  \label{eq:etsharp}
 e_{it}^{\sharp}(z) 
= \oline{e^{it(i\beta + \oline z)}}
= e^{-t\beta} e_{-it}(z),\quad \mbox{ i.e. } \quad 
 e_{it}^{\sharp}
 = e^{-t\beta} e_{-it}.
\end{equation}

\subsection{A unitary subgroup of $H^\infty$}
\mlabel{subsec:hinfty-homo}

The multiplication representation
  of $\R$ on $L^2(\partial \bS_\beta)$ defined by
  $(e_{it})_{t \in \R}$ is {\bf not} unitary on the summand
 corresponding to the upper boundary. The restriction to this
  second summand is $e^{-t \beta}$ times a unitary operator. In particular,
  $\R$ acts by normal operators. 
  Therefore the most natural structure to consider is the homomorphism
  \[ \R \to H^\infty(\bS_\beta), \quad t \mapsto e_{it}. \]
  The following observation shows that the $*$-algebra
  $(H^\infty(\bS_\beta), \sharp)$ plays the role of a ``group algebra''
  of $\R$ in the sense that its states correspond to normalized
  positive definite functions on~$\R$. 

  \begin{prop}  \mlabel{prop:5.3}
    The functions $g_t = e^{t\beta/2} e_{it}$, $t \in \R$,
    define a unitary $1$-parameter group of the Banach $*$-algebra
    $(H^\infty(\bS_\beta),\sharp)$. It has the following properties:
    \begin{itemize}
    \item[\rm(a)] Its  range spans a weakly dense $*$-subalgebra
    and $\|g_t\| = e^{|t|\beta/2}$ for $t \in\R$. 
    \item[\rm(b)] 
    Every weakly continuous state $\omega$ of $H^\infty(\bS_\beta)$
  is   represented by  a probability measure on 
  the fixed point set $\bS_\beta^\sigma = \frac{\beta i}{2} + \R$ via
  \[  \omega_\mu(F)  = \int_{\R} F\Big(\frac{\beta i}{2} + x\Big)\, d\mu(x)
    \quad \mbox{ for } \quad
  F \in H^\infty(\bS_\beta) \] 
  and satisfies
  \[
    \omega_\mu(g_t) = \hat\mu(t)
    = \int_{\R} e^{it\lambda} \, d\mu(\lambda)\quad \mbox{ for } \quad
    t \in \R.\]
  All continuous normalized positive definite
  functions on $\R$ are of this form. 
    \end{itemize}
  \end{prop}

  \begin{prf} (a) First we recall from \cite[Lemma~B.6]{ANS22} that
    the functions $(e_{it})_{t \in \R}$ span a weakly dense
    subspace of $H^\infty(\bS_\beta)$. By $g_t g_s = g_{t+s}$ 
    and \eqref{eq:etsharp} we get 
  \[ g_t^\sharp
    = e^{t\beta/2} e^{-t\beta} e_{-it} 
    = e^{-t\beta/2} e_{-it}  = g_{-t} = g_t^{-1},\]
so that  the span of the $g_t$ is a $*$-subalgebra. Finally, we observe that
  \[ \|g_t\|
    = \sup \{ e^{t\beta/2} e^{-ty} \: 0 < y < \beta\}
    = \max \{ e^{t\beta/2}, e^{-t\beta/2}\} = e^{|t|\beta/2}.\]

\nin (b)  That weakly continuous states of
  $H^\infty(\bS_\beta)$ correspond to probability measures
  $\mu$ on $\frac{\beta i}{2} + \R$ follows from \cite[Prop.~B.8]{ANS22}.
Restricting to the unitary $1$-parameter group $(g_t)_{t \in \R}$,
  the  state $\omega_\mu$ defines a positive
  definite function on $\R$, given by
    \begin{equation}
    \label{eq:phimu}
\omega_\mu(g_t) 
    = \int_{\R} g_t\Big(\frac{i\beta}{2} + \lambda\Big) \, d\mu(\lambda)
    = \int_{\R} e^{t\beta/2} e_{it}\Big(\frac{i\beta}{2} + \lambda\Big) \, d\mu(\lambda)
    = \int_{\R} e^{it\lambda} \, d\mu(\lambda)
    = \hat\mu(t).
  \end{equation}
  According to Bochner's Theorem, we thus obtain all
  continuous positive definite functions on $\R$.
\end{prf}

\subsection{Reflection positive functions on the circle}

The group $\T_\beta = \R/\beta \Z$ 
is a circle of length~$\beta$. 
We write $[t] := t + \beta \Z$ 
for the image of $t$ in $\T_\beta$ and write 
\[ \T_{\beta,+} := \{ [t] \in \T_\beta \: 0 < t < \beta/2\} \] 
for the corresponding semicircle. 
We further fix the involutive automorphism $-\id$, given by 
inversion in $\T_\beta$. 
In the following we identify functions on $\T_\beta$ with $\beta$-periodic functions on $\R$. 

We recall that a function $\phi \: \T_\beta = \R/\beta \Z\to \C$ is 
 reflection positive w.r.t.\ $(\T_\beta, \T_{\beta,+},\tau_\beta)$ 
if it is positive definite and the kernel 
\[ (\phi(t + s))_{0 < t,s < \beta/2} \] 
is positive definite (cf.\ Definition~\ref{def:posdeffunc}).

\begin{ex} \mlabel{ex:5.3} 
(see \cite[Ex.~2.3]{NO15}) Basic examples of reflection positive 
functions on $\T_\beta$ are given on the interval $[0,\beta]$ by  
\begin{equation}
  \label{eq:phi-lambda}
 \phi_\lambda([y]) := c_\lambda(iy) := \frac{e^{-y\lambda} + e^{-(\beta - y)\lambda}}
{1 + e^{-\beta \lambda} }
\quad \mbox{ for } \quad 0 \leq y \leq \beta, \lambda \geq 0,
\end{equation}
(cf.\ Lemma~\ref{lem:polweakdense} for $c_\lambda$). 
A direct calculation shows that the 
Fourier series of the $\beta$-periodic extension  of 
$\phi_\lambda$ to $\R$ (also denoted $\phi_\lambda$) is given by 
\begin{equation}
  \label{eq:fourexp}
\phi_\lambda([y]) = \sum_{n \in\Z} c_n 
\exp\Big(\frac{2\pi i n y}{\beta}\Big) \quad \mbox{ with }  \quad 
c_n = 
\frac{1}{\pi} \frac{(\frac{\beta \lambda}{2\pi})\phantom{\big)}}
{(\frac{\beta\lambda}{2\pi})^2 + n^2}
\cdot \frac{1 - e^{-\beta \lambda}}{1+ e^{-\beta \lambda}},
\end{equation}
The positive definiteness of the 
functions $\phi_\lambda$ 
on the circle group $\T_\beta$ follows immediately 
from the positivity 
of the Fourier coefficients \eqref{eq:fourexp},
which can also be obtained from the Poisson summation formula 
(Lemma~\ref{lem:poisson}). 

The kernel $(\phi_\lambda(t+s))_{0 < t,s < \beta/2}$ is positive definite 
because $\phi_\lambda$ is the Laplace transform of the positive measure 
  $\frac{\delta_\lambda + e^{-\beta \lambda} \delta_{-\lambda}}{1 + e^{-\beta\lambda}}$. 
Therefore $\phi_\lambda$ is reflection positive 
on $(\T_\beta, \T_{\beta,+}, \tau_\beta)$. 
\end{ex}

General continuous 
reflection positive functions $\phi$ on $\T_\beta$ have an
integral representation 
\begin{equation}
  \label{eq:genrptbeta}
 \phi([y])
 = \int_{[0,\infty)} 
\frac{e^{-\lambda y} + e^{-\lambda(\beta-y)}}{1 + e^{-\lambda\beta}}\, d\mu(\lambda) 
 = \int_{[0,\infty)}  \phi_\lambda([y])\, d\mu(\lambda), 
 \qquad 0 \leq y \leq \beta,
\end{equation}
where $\mu$ is a finite positive Borel measure on $[0,\infty)$ 
(\cite[Thm.~2.4]{NO15} and \cite{KL81}). \\

The following lemma is easily verified. It relates measures
$\mu$ on $[0,\infty)$, that parametrize reflection positive functions
on $\bT_\beta$ to measures on $\R$ satisfying the
{\it $\beta$-reflection relation} \eqref{eq:nurel2}. 

  \begin{lem} \mlabel{lem:mu-nu}
    For a finite positive Borel measure $\mu$ on $[0,\infty)$
    and $d\mu^\vee(\lambda):= d\mu(-\lambda)$, 
    consider the following measures on $\R$
\begin{equation}
      \label{eq:gamma-mu}
 \gamma(\mu) 
:= \mu + e_{\beta}\mu^\vee, 
    \end{equation}
   \begin{equation}
      \label{eq:Gamma-mu}
 \Gamma(\mu) 
:= \frac{1}{1 + e_{-\beta}}\mu + \frac{1}{1 + e_{-\beta}}\mu^\vee 
= \frac{\mu + \mu^\vee}{1 + e_{-\beta}}.
    \end{equation}
Then both $\gamma$ and $\Gamma$ define a bijection between finite positive Borel measures
on $[0,\infty)$ and finite positive Borel measures $\nu$ on  $\R$, satisfying
\begin{equation}
  \label{eq:nurel2} d\nu(-\lambda) = e^{-\beta \lambda} d\nu(\lambda).
\end{equation}
Moreover, $\gamma\circ M_\kappa=\Gamma$, where $M_\kappa$ is the multiplication operator by the function $\kappa=\frac1{1+e_{-\beta}}$.
\end{lem}

\begin{prf}
  It is straightforward to check that the maps $\gamma$ and $\Gamma$ are bijections between finite Borel measures on $\mathbb{R}_+$
  and on $\mathbb{R}$ which verify \eqref{eq:nurel2}.
  For a finite Borel measure $\mu$ on $\mathbb{R}_+$, the measure
  $M_\kappa(\mu)$ is also finite. We have
\begin{align*}\gamma\circ M_\kappa(\mu)&=\gamma\left(\frac1{1+e_{-\beta}}\mu\right)=\frac1{1+e_{-\beta}}\mu+e_\beta\left(\frac1{1+e_{-\beta}}\mu\right)^\vee=\frac1{1+e_{-\beta}}\mu+\frac{e_\beta}{1+e_{\beta}}\mu^\vee\\
&=\frac1{1+e_{-\beta}}\mu+\frac1{1+e_{-\beta}}\mu^\vee=\Gamma(\mu). \qedhere
\end{align*}
\end{prf}

\begin{ex} For $\lambda > 0$ and the point measure
    $\mu = \delta_\lambda$, we obtain
    the measure 
    \[ \Gamma(\delta_\lambda) 
    = \frac{1}{1 + e^{-\lambda\beta}} \delta_\lambda
      + \frac{1}{1 + e^{\lambda\beta}} \delta_{-\lambda} 
        = \frac{1}{1 + e^{-\lambda\beta}} (\delta_\lambda
        + e^{-\lambda \beta}  \delta_{-\lambda}) \]
      supported on the $2$-point subset $\{\pm \lambda\} \subeq \R$.
      Its Fourier transform satisfies 
      \[ \Gamma(\delta_\lambda)^\wedge(iy) 
        =  \frac{e^{-\lambda y} +
        e^{-\beta\lambda} e^{\lambda y}}
        {1 + e^{-\lambda\beta}}
        = \phi_\lambda([y]) \quad \mbox{ for } \quad
        0 \leq y < \beta.\]
  \end{ex}

  \begin{thm} \mlabel{thm:5.7}
    {\rm(Reflection positive functions on $\T_\beta$ and
    KMS conditions)} 
  Any reflection 
  positive function $\phi_\T$ on $\T_\beta$ extends analytically
  to a continuous function $\Phi$
  on the closed strip $\oline{\bS_\beta}$ that is
  holomorphic on the interior. We thus obtain a bijection 
  $\phi_\T \mapsto \phi_\R := \Phi\res_\R$ from continuous reflection
  positive functions on $\T_\beta$ to those positive definite functions
  $\phi_\R \: \R \to \C$ satisfying the $\beta$-KMS condition, i.e.,
  $\phi_\R$ extends to a continuous function $\oline{\cS_\beta} \to \C$,
  holomorphic on the interior and satisfying 
  \begin{equation}
    \label{eq:psi-kms}
    \phi_\R(i \beta+t) = \oline{\phi_\R(t)}\quad \mbox{ for } \quad t \in \R.
  \end{equation}
 Any reflection positive function $\varphi_\T$ has an integral
    representation 
  \begin{equation}
  \label{eq:phi-x}
  \phi_\T([y])
  = \int_{[0,\infty)} \phi_\lambda([y])\, d\mu(\lambda) \quad \mbox{ for } \quad
  0 \leq y \leq \beta, 
\end{equation}
where $\mu$ is a finite positive Borel measure on $[0,\infty)$, and the
function $\phi_\R = \hat\nu$ is the Fourier transform of the measure 
$\nu = \Gamma(\mu)$ defined in \eqref{eq:gamma-mu}. 
Vice versa, if $\phi_\R = \hat\nu$ for a finite Borel measure $\nu$ on $\mathbb{R}$, then $\phi_\T([y]):= \hat\nu(iy)$ for
$0 \leq y \leq \beta$ defines a reflection positive function on $\T_\beta$. 
\end{thm}

\begin{prf} By \eqref{eq:genrptbeta},  any reflection
positive function $\phi_\T$ on $\T_\beta$ has an integral representation
\begin{equation}
  \label{eq:phi-xb}
  \phi_\T([y])
  = \int_{[0,\infty)} \phi_\lambda([y])\, d\mu(\lambda)
= \int_{[0,\infty)} c_\lambda(iy)\, d\mu(\lambda)
  \quad \mbox{ for } \quad
  0 \leq y \leq \beta,
\end{equation}
   where $\mu$ is a finite positive Borel measure on $[0,\infty)$.
     As $\|c_\lambda\|_\infty = 1$ for $\lambda \geq 0$
     (Lemma~\ref{lem:ct}), we obtain an analytic continuation
     \[  \Phi \: \oline{\bS_\beta} \to \C, \quad
     \Phi(z)
     = \int_{[0,\infty)} c_\lambda(z)\, d\mu(\lambda) \]
       which is continuous on the closed strip $\oline{\bS_\beta}$ and
       holomorphic on the interior
       (cf.\ \cite[Rem.~2.5]{NO15}).
       The positive definiteness of the kernel
       $(\phi_\T([t + s]))_{0 < t,s < \beta/2}$ now implies that the kernel
       \[  \Phi\Big(\frac{z - \oline w}{2}\Big) \]
       on $\bS_\beta$ is positive definite (\cite[Thm.~A.1]{NO14}).
       Restricting to $\R$, we thus obtain a positive definite function
       $\phi_\R := \Phi\res_\R.$ 
       To represent this function as the Fourier transform
       of a positive measure $\nu$ on~$\R$ (Bochner's Theorem), by Lemma \ref{lem:mu-nu} we have
      \begin{align*}
  \phi_\R(x)
  &= \int_{[0,\infty)} \frac{e^{i\lambda x} + e^{-\beta\lambda} e^{-i\lambda x}}
    {1 + e^{-\beta \lambda}}\, d\mu(\lambda)
    = \int_{[0,\infty)} e^{i\lambda x}\ \frac{d\mu(\lambda)}
      {1 + e^{-\beta \lambda}}
      + \int_{(-\infty,0]} e^{i\lambda x}e^{\beta\lambda}\ \frac{d\mu(-\lambda)}
    {1 + e^{\beta \lambda}}\\
    &= \int_{[0,\infty)} e^{i\lambda x}\ \frac{d\mu(\lambda)}
      {1 + e^{-\beta \lambda}}
      + \int_{(-\infty,0]} e^{i\lambda x}\ \frac{d\mu(-\lambda)}
      {1 + e^{-\beta \lambda}}
      = \int_\R e^{i\lambda x}\ d\Gamma(\mu)(\lambda) = \Gamma(\mu)^\wedge{}(x). 
\end{align*}
Now $\nu := \Gamma(\mu)$ is a finite positive Borel measure on $\R$ satisfying
\begin{equation}
  \label{eq:nurel} d\nu(-\lambda) = e^{-\beta \lambda} d\nu(\lambda)
\quad \mbox{ and }\quad
\cL(\nu)(\beta) = \Phi(i\beta) = \mu([0,\infty))< \infty.
\end{equation} 
If, conversely, $\varphi_\R$ is a positive definite function on $\mathbb{R}$ that verifies the $\beta$-KMS-condition~\eqref{eq:psi-kms}, then there exists a finite positive Borel measure $\nu$ on $\R$ such that $\varphi_\R=\widehat{\nu}$, and $\nu$ satisfies \eqref{eq:nurel} (\cite[Theorem 2.6]{NO19}). 
  Then   $\hat\nu$ extends to continuous function on $\oline{\bS_\beta}$
  which is holomorphic on its interior. Its restriction
  to the vertical segment $[0,\beta]i$ then defines a reflection positive
  function $\phi_\T([y]):= \hat\nu(iy)$ on $\T_\beta$ (cf.\ \cite[Rem.~2.5]{NO15}).
\end{prf}

\begin{rem} ($\beta$-measures as combinations of measures on
  $\R_\pm$)   We have seen in Theorem \ref{thm:5.7} that,
  for a general reflection positive function $\varphi_{\T}$ on the circle $\T$, there exists a finite positive Borel measure on $[0,\infty)$ such that $\varphi_\T([y])=\int_{[0,\infty)}c_{\lambda}(iy)d\mu(\lambda)$.
  As a consequence, $\varphi_{\T}$ extends analytically
  to a function $\Phi$ which is continuous on the closed $\beta$-strip $\overline{\mathbb{S}_\beta}$ and holomorphic on $\mathbb{S}_\beta$. The restriction to $\R$ provides a continuous positive definite function $\varphi_\R=\widehat{\Gamma(\mu)}$, where $\Gamma(\mu)$ is a finite positive Borel measure on $\R$ defined as in Lemma~\ref{lem:mu-nu}. 

Let $\nu$ be the restriction of $\Gamma(\mu)$ to $[0,\infty)$. Then $\nu=\frac{\mu}{1+e_{-\beta}}$ is a finite positive Borel measure on $[0,\infty)$. The function $\varphi_{\mathbb{R}}^+(t):=\int_{[0,\infty)}e^{it\lambda}d\nu(\lambda)$ for $t\in\R$ admits an analytical extension $\Psi_\nu:\overline{\C_+}\to \C$ defined by
$$\Psi_\nu(z):=\int_{[0,\infty)} e^{i\lambda z}d\nu(\lambda).$$ The function $\Psi_\nu$ is continuous on $\overline{\C_+}$ and holomorphic on $\C_+$.
If $\theta$ is the unitary involution defined by $\theta f(w):=f(\beta i-w)$ for a holomorphic function $f$ on $\C_+$, then $\theta \Psi_\nu$ is continuous on $\overline{\beta i-\C_+}$ and holomorphic in $\beta i-\C_+$. Explicitly, $\theta\Psi_\nu$ takes the form 
\[ \theta\Psi_\nu(w)=\Psi_\nu(\beta i-w)=\int_{[0,\infty)} e^{-\beta \lambda}e^{-i\lambda w}d\nu(\lambda)\quad \mbox{ for } \quad w\in\beta i-\C_+.\] 
The function $\Psi_\nu+\theta\Psi_\nu$ is continuous on $\overline{\mathbb{S}_\beta}$ and holomorphic on $\mathbb{S}_\beta$. In particular, for $0\leq y\leq \beta,$
\begin{align*}
\Psi_\nu(iy)+\theta\Psi_\nu(iy)&=\int_{[0,\infty)} \frac{e^{-\lambda y}}{1+e^{-\beta\lambda}}d\mu(\lambda)+\int_{[0,\infty)} \frac{e^{-\beta \lambda}e^{\lambda y}}{1+e^{-\beta\lambda}}d\mu(\lambda)=\int_{[0,\infty)}\frac{e^{-\lambda y}+e^{-\lambda(\beta-y)}}{1+e^{-\beta\lambda}}d\mu(\lambda)=\varphi_\T([y]). 
\end{align*}

On the other hand, if $\varphi_{\R}$ is a continuous positive definite function on $\R$, by Bochner's Theorem, 
there exists a finite positive Borel measure $\rho$ on $\R$ such that $\varphi_\R=\hat{\rho}$. As before, $\Psi_{\rho_+}:\overline{\C_+}\to\C_+$ defined by
$\Psi_{\rho_+}(z):=\widehat{\rho_+}(z)$, where $\rho_+$ is the
restriction of $\rho$ to $[0,\infty)$, is continuous on $\overline{\C_+}$ and holomorphic in $\C_+$. It leads to the following function, which is continuous on $\overline{\mathbb{S}_\beta}$ and holomorphic in $\mathbb{S}_\beta$
\begin{align*}
  \Psi_{\rho_+}(z)+\theta\Psi_{\rho_+}(z)
&=\int_{[0,\infty)} e^{i\lambda z}d\rho_+(\lambda)+\int_{[-\infty,0)} e^{\beta \lambda}e^{i\lambda z}d\rho_+(-\lambda)=\int_{\R}e^{i\lambda z}d(\gamma(\rho_+))(\lambda),
\end{align*}
where $\gamma(\rho_+)$ is the
finite positive Borel measure on $\R$ from
Lemma \ref{lem:mu-nu}. The restriction of $\Psi_{\rho_+}+\theta\Psi_{\rho_+}$ to $\R$ defines a continuous positive definite function that satisfies
the $\beta$-KMS condition. Thus by Theorem \ref{thm:5.7}, $\Psi_{\rho_+}(iy)+\theta\Psi_{\rho_+}(iy)$ for $0\leq y\leq \beta$ is a reflection positive function on $\T$. 
\end{rem}

We have seen above that every reflection positive function on $\T$
decomposes as sum of a holomorphic function on $\C_+$
defined through a continuous positive definite function on $\R$,
and its image under an involution $\theta$, which produces a holomorphic
function on $\beta i-\C_+$. This decomposition is closely related to 
the periodization procedures that appear in Section \ref{sec:6}.

\subsection{Reproducing kernel Hilbert spaces on the strip}
\mlabel{subsec:repro-strip}

In Theorem~\ref{thm:5.7} we have seen how reflection positive
functions on $\T_\beta$ lead by analytic continuation
to holomorphic functions $\hat\nu$ on the strip $\bS_\beta$.
Then the functions $\hat\nu(z-\oline w)$ define positive
definite kernels on $\bS_{\beta/2}$. The following theorem
clarifies the representation theoretic significance of
such kernels for Hilbert spaces of holomorphic functions on
$\bS_\beta$ on which translations define unitary operators.

\begin{thm} \mlabel{thm:7.1}
{\rm(a)}  If $\nu$ is a positive Borel measure on 
$\R$ whose Laplace transform \break
$\cL(\nu)(y) = \int_\R e^{-\lambda y}\,
d\nu(\lambda)$ is finite on the interval
  $(0,2\beta)$, then
  \[ K(z,w) := \hat\nu(z-\oline w) = \int_\R e^{i\lambda(z - \oline w)} \,
  d\nu(\lambda)\]
  is a translation invariant sesquiholomorphic positive
  definite kernel on
  the strip $\bS_\beta$, and all translation invariant
  sesquiholomorphic kernels on $\bS_\beta$ are of this  form, 
  for a uniquely determined measure~$\nu$.

  \nin{\rm(b)}  The measure $\nu$ satisfies 
  \begin{equation}
    \label{eq:nu-beta-x}
    d\nu(-\lambda) = e^{-2\beta \lambda} d\nu(\lambda)
  \end{equation}
if and only if
$(\theta F)(z) := F(\beta i - z)$ 
defines a unitary involution on the corresponding
reproducing kernel space~$\cH_K$.

\nin {\rm (c)} We have a unitary map 
\begin{equation}
  \label{eq:psi-b.7}
 \Psi \: L^2(\R, \nu) \to \cH_K,\quad
  \Psi(f)(z) = \la e_{-i\oline z}, f \ra \quad \mbox{ for} \quad
  z \in \cS_\beta, f \in L^2(\R, \nu).
\end{equation}
\end{thm}

\begin{prf} (a) If $\nu$ is given, then 
$\hat\nu(z)
= \int_\R e^{iz\lambda}\, d\nu(\lambda)$
defines a holomorphic function on $\bS_{2\beta}$. With 
$e_z(\lambda) := e^{\lambda z}$, we see that 
\[ K(z,w) := \hat\nu(z - \oline w)
= \int_{\R} e^{i\lambda z}e^{-i\lambda \oline w}\, d\nu(\lambda)
= \la e_{iw}, e_{iz} \ra_{L^2(\R,\nu)} \]
is a positive definite translation invariant
sesquiholomorphic kernel on $\bS_\beta$.

Suppose, conversely, that $K$ is a
translation invariant sesquiholomorphic positive
  definite kernel on $\bS_\beta$, i.e.,
\[ K(z,w) = K(z + t, w + t) \quad \mbox{ for } \quad t \in \R,
z,w \in \bS_\beta.\]
A direct analytic continuation argument then shows that
$K(z,w) = \kappa(z - \oline w)$ 
for a holomorphic
function $\kappa \: \bS_{2\beta} \to \C.$ 
We obtain on $\cH_K$ a unitary one-parameter group by
\[ (U_t F)(z) := F(z + t) \quad \mbox{ for } \quad t \in \R, z \in \bS_\beta.\]
We then have
\[ (U_t K_w)(z) = K_w(z + t) = K_{w-t}(z) \quad \mbox{ for }\quad
  K_w(z) := K(z,w),\]
so that we obtain on $\R$ a positive definite function 
\begin{equation}
  \label{eq:po-de-phi}
 \phi(t) :=
\la K_{\beta i/2}, U_t K_{\beta i/2} \ra 
= \la K_{\beta i/2}, K_{-t + \beta i/2} \ra 
= K\Big(\frac{\beta i}{2}, \frac{\beta i}{2} - t\Big) 
= \kappa(\beta i + t).
\end{equation}
Now Bochner's Theorem implies that there exists a finite positive
Borel measure $\mu$ on $\R$ with
\[ \hat\mu(t) = \kappa(\beta i + t).\]
We consider the measure
\[ d\nu(\lambda) := e^{\beta\lambda}\, d\mu(\lambda) \]
on $\R$.
As the Fourier transform $\hat\mu$ extends to a
holomorphic function on 
$\{ z \in \C \: |\Im z| < \beta \}$, 
we obtain 
\begin{equation}
  \label{eq:Phi}
  \kappa(z) = \hat\mu(z - \beta i) = \hat\nu(z) \quad \mbox{ for }
  \quad z \in \bS_{2\beta}.
\end{equation}
In particular, we find that
\[ \kappa(iy) = \cL(\mu)(y - \beta) = \cL(\nu)(y) \quad \mbox{ for } \quad
  0 < y < 2 \beta,\]
so that the Laplace transform $\cL(\nu)$ is finite 
on the interval $(0,2\beta)$.

\nin (b) If \eqref{eq:nu-beta-x} is satisfied, then
$\hat\nu(z) = \hat\nu(2\beta i - z)$.
We thus obtain
\[
  K(\beta i - z, \beta i - w)
= \hat\nu(2 \beta i - (z - \oline w)) 
= \hat\nu(z - \oline w) = K(z,w),\]
so that $\theta$ defines a unitary operator on~$\cH_K$.
If, conversely, this is the case, then
\[ \hat\nu(z - \oline w)
= K(z,w) = K(\beta i - z, \beta i - w)
= \hat\nu(2 \beta i - (z - \oline w))\]
implies $\hat\nu(z) =   \hat\nu(2\beta i -z)$ for
$z \in \bS_{2\beta},$ 
and this implies \eqref{eq:nu-beta-x}.

\nin (c) For $z \in \bS_\beta$, the function
$e_{-i\oline z}$ is contained in $L^2(\R,\nu)$ and these functions form a total
subset. We thus obtain an injective map
\begin{equation}
 \Psi \: L^2(\R, \nu) \to \cO(\bS_\beta),\quad
  \Psi(f)(z) = \la e_{-i\oline z}, f \ra \quad \mbox{ for} \quad
  z \in \cS_\beta, f \in L^2(\R, \nu)
\end{equation}
whose range is a reproducing kernel Hilbert space $\cH_K$.
Its kernel is given by 
\[ K(z,w) =  \la e_{-i\oline z}, e_{-i\oline w} \ra_{L^2(\R, \nu)}
  =  \hat\nu(z-\oline w). \qedhere\]
\end{prf}

If the measure $\nu$ is finite, then $\hat\nu$ is a continuous
function on $\R$, and if it is merely tempered, then
$\hat\nu$ defines a tempered distribution on $\R$, given by
the boundary values of the holomorphic function
$\hat\nu$ on $\bS_{2\beta}$ on the real line are possibly singular.
On $\R$ they define the distribution kernel
\[ K(x,y)
= \int_\R e^{i\lambda(x-y)}\, d\nu(\lambda)
= \hat\nu(x-y).\]

\begin{exs} \mlabel{ex:qs}
  From \cite[\S 6.1]{NOO21} and \cite[Thm.~VII.3.1]{FK94}, we recall the
  Riesz measures on $\R_+$, defined by 
  \begin{equation}
    \label{eq:riesz1}
 d\mu_s(p) = \Gamma(s)^{-1} p^{s-1}\, dp \quad  \mbox{ for } \quad s > 0,
  \end{equation}
where $\Gamma$ is the $\Gamma$-function. Their Fourier transform satisfies 
\[ \hat\mu_s(z) = \int_{\R_+} e^{izp}\, d\mu_s(p)
  = (-iz)^{-s} = \Big(\frac{i}{z}\Big)^s\quad \mbox{  for } \quad \Im z > 0.\] 
The boundary values on $\R$ define a tempered distribution
$\hat\mu_s$ by
  \[  \hat\mu_s(\phi) = \lim_{\eps \to 0^+}
    \int_\R \phi(x) \hat\mu_s(x + i \eps)\, dx.\]
  Here the singularity in $x=0$ can be treated in terms
  of principal values (valeurs finis).
The distribution $\hat\mu_s$
is represented  on the open subset $\R^\times$ by the function 
\begin{equation}
  \label{eq:1d-impart}
\hat\mu_s(x) 
=  e^{\sgn(x) s \frac{\pi i}{2}} |x|^{-s} 
\quad \mbox{ with } \quad  \Im\hat\mu_s(x) 
=  \sgn(x) \sin\Big(s \frac{\pi}{2}\Big) |x|^{-s}.
\end{equation}
In particular, the distribution $\Im\hat\mu_s$ vanishes on $\R^\times$  
if and only if $s \in 2 \Z$, and then  it is 
supported in~$\{0\}$, cf. \cite[Lemma 6.6]{NOO21}.

On $\C_+$ we obtain for $s > 0$ the positive definite kernels
\[ Q_s(z,w) := \Big(\frac{i}{z - \oline w}\Big)^s
  = \hat\mu_s(z - \oline w), \quad s > 0.\]
For $s = 1$ this is a multiple of the Szeg\"o kernel,
and for $s = 2$ a multiple of the Bergman kernel. The corresponding Riesz 
measure is
\[ d\mu_2(p) = p\, dp \quad \mbox{ on }\quad \R_+.\]

The reproducing kernel spaces 
  $\cH_s(\C_+) := \cH_{Q_s}$ with kernel $Q_s$
  transform naturally to $\bS_\beta$ via the exponential
  map $\Exp \: \bS_\beta\to \C_+, z \mapsto e^{\pi z/\beta}$. 
  The map
  \[ \Phi_s \: \cH_s(\C_+) \to \cH_s(\bS_\beta), \quad
    \Phi_s(f)(z) = (\Exp'(z))^{s/2} f(\Exp(z))
    = \Big(\frac{\pi}{\beta}\Big)^{s/2} \exp\Big(\frac{s\pi z}{2\beta}\Big)
    f\Big(\exp\Big(\frac{\pi z}{\beta}\Big)\Big) \]
  defines a unitary operator to a reproducing kernel Hilbert space~$\cH_s(\bS_\beta)$.   With the formula \eqref{eq:cauchyker} for the Szeg\"o kernel
  $\Q^{\C_+}$ on $\C_+$ and by using \eqref{eq:sz-trafo}, we find for this kernel the formula 
  \begin{align*}
    Q_s(z,w)
    &= (\Exp'(z))^{s/2} Q_s^{\C_+}(\Exp(z),\Exp(w)) (\oline{\Exp'(w)})^{s/2}    \\
    &= \Big(\frac{\pi}{\beta}\Big)^s
      \exp\Big(\frac{s\pi (z+ \oline w)}{2\beta}\Big)
      \frac{1}{(2\pi)^s}
      \Bigg(\frac{i}{\exp(\pi z/\beta) - \exp(\pi\oline w/\beta)}\Bigg)^s \\
    &= \Big(\frac{1}{2\beta}\Big)^s
      \Bigg(\frac{i \exp\big(\frac{\pi (z+ \oline w)}{2\beta}\big)}
      {\exp(\pi z/\beta) - \exp(\pi\oline w/\beta)}\Bigg)^s 
=       \Bigg(\frac{i}
      {4\beta\sinh\big(\frac{\pi}{2\beta}(z- \oline w)\big)}\Bigg)^s \\
&=       \Bigg(\frac{1}
      {4\beta\cosh\big(\frac{\pi}{2\beta}(z- \oline w)- \frac{\pi i}{2}\big)}\Bigg)^s. 
  \end{align*}
  So $Q_s(z,w) = Q(z,w)^s$, where $Q$ is the Szeg\"o kernel
  on $\bS_\beta$ (see \eqref{eq:cskernel-strip} in
  Appendix~\ref{subsec:app-cz-kernel}).
\end{exs}

\begin{rem} Suppose that $\psi$  is a  positive definite
distribution $\psi$ on $\R$ satisfying the $2\beta$-KMS condition
\[ \psi(t + 2\beta i) = \oline{\psi(t)} \]
in the sense of boundary values of the holomorphic extension to the
strip $\bS_{2\beta}$. By the Bochner--Schwartz Theorem
  (\cite[Thm.~XVIII, \S VII.9]{Schw73}), there exists a tempered
measure $\nu$ on $\R$
with $\psi = \hat\nu$ that satisfies
\[ d\nu(-p) = e^{-2\beta p}\, d\nu(p)\]
(cf.\ Theorem~\ref{thm:5.7} for the finite case). 
For the imaginary part 
\[ \kappa(t) := 2 i \Im \psi(t) = \psi(t) - \psi(-t)\]
we then obtain 
\begin{align*}
 \kappa(t)
  &  = \int_\R e^{itp}\, d\nu(p) - \int_\R e^{-itp}\, d\nu(p)
    = \int_\R e^{itp}\, d\nu(p) - \int_\R e^{itp}\, d\nu(-p)\\  
  &  = \int_\R e^{itp}\, d\nu(p) - \int_\R e^{itp - 2\beta p}\, d\nu(p)
    = \int_\R e^{itp}(1 - e^{-2\beta p})\, d\nu(p).
\end{align*}
Hence $\nu(\{0\})$ does not contribute to $\kappa$,
but as the Fourier transform is bijective on tempered distributions, 
it follows that $\nu\res_{\R^\times}$ is uniquely determined by $\kappa$
(cf.\ \cite[p.~619]{BY99}).   
\end{rem}

\begin{ex} Typical examples that arise naturally in Quantum
  Field Theory   are the measures
  \[ d\nu(p) = \frac{p \gamma(p^2)dp}{1 - e^{-2\beta p}}, \]
  where $\gamma$ is a polynomial that is $\geq 0$ on $\R_+$
  (\cite{Yn94}, \cite{BY99}).
  For $\gamma(x) = x^n$, $n \in \N_0$, we thus obtain
  measures with the non-singular density
  $\frac{p^{2n+1}dp}{1 - e^{-2\beta p}}$.
  These restrict on $\R_+$ to the measures 
  \begin{equation}
    \label{eq:mu2n+2}
    \frac{(2n+1)!\; d\mu_{2n+2}(p)}{1 - e^{-2\beta p}}.
  \end{equation}
\end{ex}

\begin{ex}
For $s > 0$, consider on $\R$ the positive measure
\[ d\nu_s(p) := \frac{1}{\Gamma(s)} \frac{|p|^{s-2} p\, dp}{1 - e^{-2\beta p}},
  \quad \mbox{ i.e.,} \quad
  \nu_s =
  \frac{\mu_s}{1 - e_{-2\beta}} + 
  \frac{\mu_s^\vee}{e_{-2\beta}-1} \]
with $\mu_s$ as in \eqref{eq:riesz1} in Example~\ref{ex:qs}.  
It satisfies the $2\beta$-KMS relation 
\begin{align*}
 d\nu_s(-p)
&  = \frac{1}{\Gamma(s)} \frac{|p|^{s-2} p\, dp}{e^{2\beta p}-1}
  = e^{-2\beta p} \frac{1}{\Gamma(s)} \frac{|p|^{s-2} p\, dp}{1-e^{-2\beta p}}
=  e^{-2\beta p}\, d\nu_s(p).
\end{align*}
Further, 
\begin{align*}
 \kappa(t)
  &  = \int_\R e^{itp} (1 - e^{-2\beta p})\, d\nu_s(p)
    = \frac{1}{\Gamma(s)} \int_{\R_+} e^{itp}|p|^{s-1}\, dp 
  - \frac{1}{\Gamma(s)} \int_{\R_-} e^{itp} |p|^{s-1}\, dp \\ 
  &= \int_{\R_+} e^{itp}\, d\mu_s(p)
    -  \int_{\R_+} e^{-itp}\, d\mu_s(p) = \hat\mu_s(t) - \hat\mu_s(-t)
    = 2 i \Im \hat\mu_s(t),
\end{align*}
in the sense of tempered distributions.
By \eqref{eq:1d-impart}, this distribution is given on $\R^\times$
by the function
\[ 2 i \sgn(t) \frac{\sin(\pi s/2)}{|t|^s}.\]
In particular, it vanishes on $\R^\times$ for
$s  \in 2 \N$ and then the distribution $\kappa$ is supported in~$\{0\}$.

To evaluate these distributions, we recall from
\cite[Ex.~6, \S V.3]{RS75} that
\[\lim_{\eps \to 0^+} \frac{1}{x+i\eps}
  = \lim_{\eps \to 0^+} \frac{x- i \eps}{x^2 + \eps^2} 
  = \lim_{\eps \to 0^+} \frac{x}{x^2 + \eps^2} - \frac{i\eps}{x^2 + \eps^2} 
  = \cP\Big(\frac{1}{x}\Big) - \pi i \delta_0,\]
where $\cP$ denotes the {\it principal value}. 
We conclude that
\[ \lim_{\eps \to 0^+} \Im \frac{1}{x+i\eps} = - \pi \delta_0,\]
and by taking derivatives in the distributional sense, that
\begin{equation}
  \label{eq:ims}
  \lim_{\eps \to 0^+}   \Im \frac{1}{(x+i\eps)^n}
  = \frac{(-1)^n\pi}{(n-1)!} \delta_0^{(n-1)}.
\end{equation}

For $s = 2 + 2n$, we have
\[ \hat\mu_{2 + 2n}(z) = \Big(\frac{i}{z}\Big)^{2 + 2n} =
  \frac{(-1)^{n+1}}{z^{2 + 2n}},\]
  so that the boundary values on $\R$ can be written as
  \[ \hat\mu_{2 +2n}(x) = \lim_{\eps \to 0^+} \frac{(-1)^{n+1}}{(x + i \eps)^{2 + 2n}},\]
  which  leads with \eqref{eq:ims} to
  \[ \Im\hat\mu_{2 + 2n} =  \frac{\pi}{(2n+1)!} 
    \delta_0^{(2n+1)}.\] 
An explicit description of the Fourier transform
  of $\nu_{2 + 2n}$ can be derived from the special case $n = 0$ as follows.
  We have
  \[ d\nu_{2n +2}(p)
    = \frac{1}{(2n+1)!} \frac{p^{1 + 2n}}{1 - e^{-2\beta p}} \, dp
    = \frac{p^{2n}}{(2n+1)!} \, d\nu_2(p)
    \quad \mbox{ with } \quad
    d\nu_2(p) = \frac{pdp}{1 - e^{-2\beta p}}.\]
  Therefore the $2n$th derivative of $\hat\nu_2$ given by $\hat\nu_2^{(2n)}
    = (2n+1)! (-1)^n \hat\nu_{2n+2}$ 
    leads to
    \[ \hat\nu_{2n+2}  = \frac{(-1)^n}{(2n+1)!}\hat\nu_2^{(2n)}. \] 

  The measure $(1 + e_{-2\beta}) \nu_s$ is symmetric.
  For $s = 2$, it takes the form
  \begin{equation}
    \label{eq:theta-form}
    \frac{1 + e^{-2\beta p}}{1 - e^{-2\beta p}} p\, dp
    = \frac{\cosh(\beta p)}{\sinh(\beta p)} p\, dp.
  \end{equation}
\end{ex}

\subsection{Relating to standard subspaces}
\mlabel{subsec:5.4}

We assume the context of Theorem~\ref{thm:5.7},
where $\nu$ is a finite positive measure on $\R$ 
    satisfying 
    \begin{equation}
      \label{eq:beta-rel-x}
      d\nu(-\lambda)    = e^{-\beta\lambda} d\nu(\lambda).
    \end{equation}
    In Theorem~\ref{thm:5.7} we have seen that \eqref{eq:beta-rel-x}
    corresponds to the $\beta$-KMS relation
    \[ \hat\nu(\beta i + x) =  \oline{\hat\nu(x)} \quad \mbox{ for } \quad
      x \in \R, \]
and  this in turn relates to standard subspaces via
\eqref{eq:fac2} below. This motivates our discussion of standard
subspaces in the present context. 

\begin{defn} A closed real subspace $\sV$ of the complex Hilbert space
  is called {\it standard} if
  $\sV + i \sV$ is dense and $\sV \cap i \sV = \{0\}$.
If $\sV$ is standard, then
the complex conjugation $T_\sV$ on $\sV + i \sV$ has a polar decomposition
$T_\sV = J_\sV \Delta_\sV^{1/2}$, where
$J_\sV$ is a conjugation (an antilinear involutive isometry)
and $\Delta_\sV$ is a positive selfadjoint operator
satisfying $J_\sV \Delta_\sV J_\sV =~\Delta_\sV^{-1}$.
The unitary one-parameter group $(\Delta_\sV^{it})_{t \in \R}$,
the {\it modular group of $\sV$}, preserves the subspace~$\sV$.
We also note that, for any pair $(\Delta, J)$ of a positive selfadjoint
operator $\Delta$ and a conjugation $J$, satisfying
$J \Delta J = \Delta^{-1}$, the subspace $\sV = \Fix(J \Delta^{1/2})$
is standard with $T_\sV = J \Delta^{1/2}$. 
We refer to \cite{Lo08, NO17} for more on standard subspaces.
\end{defn}

\begin{proposition}\mlabel{prop:5.11}
  {\rm(Standard subspaces for multiplication representations)} 
    Let $\nu$ be a tempered positive measure on $\R$ satisfying 
    $d\nu(-\lambda)    = e^{-\beta\lambda} d\nu(\lambda)$.
  On $L^2(\R,\nu)$ we consider the positive selfadjoint operator
and the conjugation defined by 
\[ (\Delta f)(\lambda) := e^{-\beta \lambda} f(\lambda)
\quad \mbox{ and } \quad
(Jf)(\lambda) := e^{-\beta\lambda/2} \oline{f(-\lambda)}. \]
  Then $J \Delta J = \Delta^{-1}$ and the standard subspace  
  $\sV := \Fix(J \Delta^{1/2})\subeq L^2(\R,\nu)$ is 
  \[ \sV     = \{ f \in L^2(\R,\nu) \:  f(\lambda)  = \oline{f(-\lambda)}\}.\]
\end{proposition}

\begin{prf} That $J$ is a conjugation on $L^2(\R,\nu)$
  follows from 
  \begin{align*}
   \int_\R |(Jf)(\lambda)|^2\, d\nu(\lambda)
&=     \int_\R e^{-\lambda\beta} |f(-\lambda)|^2\, d\nu(\lambda)
=     \int_\R e^{\lambda\beta} |f(\lambda)|^2\, d\nu(-\lambda)\\
&=     \int_\R e^{\lambda\beta} |f(\lambda)|^2 e^{-\beta\lambda}\, d\nu(\lambda)
    =     \int_\R |f(\lambda)|^2 \, d\nu(\lambda).
  \end{align*}
  For $\Delta$ we then have $J \Delta J = \Delta^{-1},$ 
  so that $\sV := \Fix(J \Delta^{1/2})$
  is a standard subspace. 
  Its elements are characterized by the relation
 $e^{-\beta\lambda/2} f(\lambda)  = e^{-\beta \lambda/2}    \oline{f(-\lambda)}$ 
  which is equivalent to $f(\lambda)  = \oline{f(-\lambda)}$.
\end{prf}

\begin{rem} If, in the context of Proposition~\ref{prop:5.11}, 
  the measure $\nu$ is finite, then
  the constant function $1$ is contained in the standard subspace
  $\sV$ and for the
  representation of $\R$ on $L^2(\R,\nu)$, given by 
  \[ (U_t f)(\lambda) := (\Delta^{-it/\beta}f)(\lambda)
    = e^{it\lambda} f(\lambda), 
    \quad \mbox{ we have } \quad
    \hat\nu(t) = \la 1, U_t 1 \ra \quad \mbox{ for } \quad t \in \R\]
  (cf.\ \cite[Thm.~2.6]{NO19})
  and $(\Delta f)(\lambda) = e^{-\beta \lambda} f(\lambda)$.
\end{rem}

\begin{ex} \mlabel{ex:5.12} 
  (a)  (Standard subspaces for the translation representation)
  We consider $\cH = L^2(\R)$ and
  the standard subspace $\sV \subeq L^2(\R)$, specified by
  \[ Jf = \oline f \quad \mbox{ and }  \quad
    (\Delta^{-it/2\beta} f)(x) = f(x + t), \quad x,t \in \R.\]
For $f \in \cD(\Delta^{1/2})$ we have 
    \[ (\Delta^{1/2}f)(x) = f(x + i \beta),  \]
    so that $f$ is fixed by $J \Delta^{1/2}$ if and only if
    $f^\sharp = f$, where
    \[ f^\sharp(x) := \oline{f(x + \beta i)}
      \quad \mbox{ for }  \quad x \in \R.\]
    This shows that 
  \begin{equation}
    \label{eq:5.12}
    \sV = \{ f \in \cD(\Delta^{1/2}) \: f^\sharp = f\}. 
\end{equation}
Endowed with the graph topology, we have
$\cD(\Delta^{1/2}) \cong \Gamma(\Delta^{1/2})$, and further
\[ \Gamma(\Delta^{1/2}) \cong H^2(\bS_\beta)
  \subeq L^2(\R)^{\oplus 2}.\] 
In this picture, the Tomita involution of $\sV$ corresponds to the
involution on $H^2(\bS_\beta)$, given by
\[ f^\sharp(z) = \oline{f(\beta i + \oline z)}
  \quad \mbox{ for } \quad z \in \bS_\beta\]
and the lower boundary value map thus induces an isometry
\begin{equation}
  \label{eq:h2sharp}
 H^2(\bS_\beta)^\sharp := \{ f \in H^2(\bS_\beta) \: f^\sharp = f\}
 \to \sV, \quad  f \mapsto f\res_\R
\end{equation}
 (cf.\ \cite[Ex.~3.16]{NO17}).
 On the pairs $(f_1, f_2) = (f, \Delta^{1/2}f)
 \in \Gamma(\Delta^{1/2}) \subeq L^2(\R)^{\oplus 2}$
of boundary values of elements of $H^2(\bS_\beta)$, we have
    \[ (f_1, f_2)^\sharp = (\oline{f_2}, \oline{f_1}).\] 

    \nin (b) For any element $v$ in a general standard subspace $\sV$,
    the function
  \begin{equation}
    \label{eq:fac2}
 \psi \: \R \to \C, \quad    \psi(t) := \la v, \Delta^{-it/\beta} v \ra
  \end{equation}
  on $\R$ is positive definite and satisfies a $\beta$-KMS condition
  such as \eqref{eq:beta-kms} 
  (cf.\ \cite[Thm.~2.6]{NO19}). 
Now take $\sV = H^2(\bS_\beta)^\sharp$ as in (a) above,
 and consider the element 
  $v = Q_{\frac{\beta i}{2}}^* \in \sV \subeq L^2(\R)$.
From
\[ Q(z,w) = \frac{1}{2\pi} \hat\nu(z - \oline w)\quad \mbox{ for } \quad
  d\nu(\lambda) =  \frac{d\lambda}{1 + e^{-2\beta \lambda}} 
\]
(see\footnote{Here and throughout the paper, we will use $\lambda$ to denote both a variable  and the Lebesgue measure $d\lambda$.} \eqref{eq:cskernel-strip} 
in Appendix~\ref{subsec:app-cz-kernel}),
we get
\[   Q_{\beta i/2}(z)=Q\left(z,\frac{\beta i}2\right)
  = \frac{1}{2\pi}\hat\nu\Big(z + \frac{\beta i}{2}\Big) 
  = \frac{1}{2\pi}
  \int_\R \frac{e^{i\lambda z} e^{-\beta \lambda/2}}{1 + e^{-2\beta \lambda}} \, d\lambda.\]
So $2\pi Q_{\beta i/2} = \widehat{\nu'}$ for a new measure
\[ d\nu'(\lambda)
  := e^{-\beta \lambda/2} d\nu(\lambda).\] 
Note that 
\begin{equation}
  \label{eq:nu'beta}
 d\nu'(-\lambda)
  = e^{\beta \lambda/2} d\nu(-\lambda) 
  = e^{\beta \lambda/2} e^{-2\beta\lambda}\, d\nu(\lambda) 
  = e^{\beta \lambda} e^{-2\beta\lambda}\, d\nu'(\lambda) 
  = e^{-\beta\lambda}\, d\nu'(\lambda).
\end{equation}
We thus derive from the isometry of the Fourier transform $L^2(\R) \to L^2(\R)$
the relation 
  \begin{align*}
    \psi(t)
    &= \la Q_{\frac{\beta i}{2}}, 
      \Delta^{-it/\beta}   Q_{\frac{\beta i}{2}} \ra
= \int_\R \oline{Q_{\frac{\beta i}{2}}(x)} 
    Q_{\frac{\beta i}{2}}(x + 2t) \, dx \\ 
    &=\frac{1}{2\pi} \int_\R
      \frac{e^{-\beta \lambda/2}}{1 + e^{-2\beta\lambda}}
      e^{2it\lambda} \frac{e^{-\beta \lambda/2}}{1 + e^{-2\beta\lambda}} \, d\lambda\\
    &=\frac{1}{2\pi} \int_\R
      \frac{e^{-\beta \lambda}}{(1 + e^{-2\beta\lambda})^2}
      e^{2it\lambda} \, d\lambda 
=\frac{1}{4\pi} \int_\R
      \frac{e^{-\beta \lambda/2}}{(1 + e^{-\beta\lambda})^2}
      e^{it\lambda} \, d\lambda.
  \end{align*}
\end{ex}

\section{From the half plane to the strip by periodization}
\mlabel{sec:6}

In this section we study the passage from functions
on the upper half plane to the strip by periodization, using that
\[ \bS_\beta = \C_+ \cap (\beta i - \C_+).\]
We show in
particular that the evaluation functionals $(Q_w^{\C_+})_{w \in \C_+}$
of the Hardy space $H^2(\C_+)$ extend meromorphically to $\C$ and 
that the Szeg\"o kernel functions
$(Q_w^{\bS_\beta})_{w \in \bS_\beta}$ of the strip $\bS_\beta$  
are represented (up to a positive factor) by the series
\[  \sum_{k \in \Z} (-1)^k Q_w^{\C_+}(z + 2 ki \beta). \]
We also provide an abstract perspective on this result which
applies to translation invariant Hilbert space on $\bS_\beta$.
In particular, a similar result
holds for the Bergman space, but in this case one has to consider the series
$\sum_{k \in \Z} Q_w(z + 2 ki \beta)$, where $Q$ is the kernel of the
Bergman space $\cO(\C_+) \cap L^2(\C_+)$.

\subsection*{Partial fractions  expansion}

The meromorphic function
\[ q(z) := \frac{i}{4\beta} \frac{1}{\sinh\big(\frac{\pi z}{2\beta}\big)} \]
represents the Szeg\"o kernel of the strip $\bS_\beta$ by
\[ Q(z,w)
=  \frac{i}{4\beta}  \frac{1}{\sinh\big(\frac{\pi(z-\oline w)}{2\beta}\big)}
= q(z- \oline w) \]
(cf.\ \eqref{eq:cskernel-strip} in Appendix~\ref{subsec:app-cz-kernel}).
The following lemma provides a partial fraction expansion
of~$q$.

\begin{lem} \mlabel{lem:q-exp} We have the partial fractions expansion
  
\begin{equation}
  \label{eq:sinh-sum}
 \frac{\pi}{2\beta} \frac{1}{\sinh\big(\frac{\pi z}{2\beta}\big)}
  =  \sum_{k \in \Z}  \frac{(-1)^k}{z + 2k i \beta}.
\end{equation}
\end{lem}
\nin Note that the series on the right  converges because
\[ \frac{1}{z + 2k i \beta} - \frac{1}{z + 2(k+1) i \beta}
  =\frac{2i\beta}   {(z + 2ki\beta)(z + 2(k+1)i\beta)} \]
has a denominator that is quadratic in~$k$.

\begin{prf} We start from the well-known partial fractions expansion
\[ \pi \cot(\pi z) 
  = \frac{1}{z} + \sum_{n \not=0} \Big(\frac{1}{z-n} + \frac{1}{n}\Big)
\quad \mbox{ that yields } \quad  \pi \cot\Big(\frac{\pi}{2}z\Big)
  = \frac{2}{z} + \sum_{n \not=0} \Big(\frac{2}{z-2n} + \frac{1}{n}\Big) \]
and
\[  \pi \tan\Big(\frac{\pi}{2}z\Big)
= \pi \cot\Big(\frac{\pi}{2}(1-z)\Big)
  = \frac{2}{1-z} + \sum_{n \not=0} \Big(\frac{2}{1-z-2n} + \frac{1}{n}\Big)
  = -\frac{2}{z-1} - \sum_{n \not=0} \Big(\frac{2}{z - (1+2n)}
  +  \frac{1}{n}\Big). \]
Using the relation $\cot(z) + \tan(z) = \frac{2}{\sin(2z)}$,
the sum of both series now provides a partial fractions 
expansion of $\frac{1}{\sin}$: 
\begin{align} \label{eq:align1}
 \frac{\pi }{\sin(\pi z)}
  &=  \frac{1}{2}\Big(\pi \cot\Big(\frac{\pi}{2}z\Big)
    +    \pi \tan\Big(\frac{\pi}{2}z\Big)\Big)\notag\\
  &=  \frac{1}{z} + \sum_{n \not=0} \Big(\frac{1}{z-2n} + \frac{1}{2n}\Big) 
    -\frac{1}{z-1} - \sum_{n \not=0} \Big(\frac{1}{z - (1+2n)}
    +  \frac{1}{2n}\Big)\notag \\
  &=  \sum_{n \in \Z} \frac{1}{z-2n} - \frac{1}{z-(2n+1)} 
  =  \sum_{n \in \Z} \frac{(-1)^n}{z-n}.
\end{align}
Rescaling now leads to the expansion 
\begin{equation}
  \label{eq:szego-sum}
 \frac{\pi}{2\beta} \frac{1}{\sinh(\frac{\pi}{2\beta} z)}
   = \frac{\pi}{2\beta} \frac{i}{\sin(\frac{\pi i}{2\beta} z)}
=  \sum_{n \in \Z} \frac{(-1)^n i}{i z-2\beta n}
=  \sum_{n \in \Z} \frac{(-1)^n} {z+ 2n\beta i}.
\qedhere\end{equation}
\end{prf}

For the Szeg\"o kernel $Q$ of $\bS_\beta$,
we thus obtain  the series expansion
\begin{equation}
  \label{eq:szego-exp} Q(z,w) =
  \frac{i}{2\pi} \sum_{n \in \Z}
  \frac{(-1)^n} {z-\oline w+ 2n\beta i}.
\end{equation}
For the meromorphic functions
\[ q_w :=   \frac{1}{2\pi}  \frac{i}{z-\oline w}, \quad w \in \C, \]
representing for $w \in \C_+$ the Szeg\"o kernel function $Q_w^{\C_+}$, 
this leads on $\bS_\beta$ to the pointwise representation 
\begin{equation}
  \label{eq:szego-exp2}
  Q_w = \sum_{n \in \Z} (-1)^nq_{w + 2n\beta i}. 
\end{equation}
This sum splits naturally as
\begin{equation}
  \label{eq:szego-exp3}
  Q_w = Q_w^+ + Q_w^- = \sum_{n \geq 0} (-1)^nq_{w + 2n\beta i}
  +   \sum_{n < 0} (-1)^nq_{w + 2n\beta i}. 
\end{equation}
We shall see below that this splitting corresponds to the decomposition
of $H^2(\bS_\beta)$ into positive and negative spectral subspaces
for the translation group and that the series on the right
converge in the Hilbert space norm.
\medskip

As $H^2(\bS_\beta)$ is generated by the functions
$(Q_w)_{w \in \bS_\beta}$, 
the preceding observation suggests to consider functions
in $H^2(\bS_\beta)$ as functions on a larger domain in $\C$ satisfying
$f(z + 2 \beta i) = - f(z)$.
This relation determines the values of
$f$ on $\bS_\beta + \Z 2\beta i$, and on the open subset 
$\bS_\beta + (2\Z + 1) \beta i$, these functions will typically
have singularities, as the functions $(Q_w)_{w \in \bS_\beta}$ show.

  On the Hardy space $H^2(\bS_\beta)$ we consider the unitary operators with 
  \[ (U_t f)(z)  = f(z + t)
    \quad \mbox{ and }  \quad \theta(f)(z) = f(\beta i - z) \quad
    \mbox{ for }\quad t \in \R, z \in \bS_\beta.\]
  Note that
  \[     \theta U_t \theta= U_{-t}. \]
  Writing $U_t = e^{itP}$, the positive selfadjoint operator
  $\Delta := e^{-2\beta P}$ satisfies
  \[ (\Delta^{1/2} f)(x) = f(\beta i + x), \ x \in \R, \quad \mbox{ for }\quad
  f \in \cD(\Delta^{1/2}) = H^2(\bS_\beta)\res_\R.\] 
So we obtain a graph realization of  the Hardy space as 
  \[ H^2(\bS_\beta) \cong \Gamma(\Delta^{1/2})
    \subeq \cD(\Delta^{1/2}) \oplus \cD(\Delta^{-1/2})
    \subeq L^2(\R)^{\oplus 2}.\]
The closed subspaces $H^2(\C_\pm) \subeq L^2(\R)$ (via non-tangential
  boundary values) are the spectral subspaces 
  of $P$ for $\pm [0,\infty)$. This implies that the
  positive spectral subspace for the operator
  $P \oplus P$ on $H^2(\bS_\beta)$ is
  \[ H^2(\bS_\beta)_+
    = H^2(\bS_\beta) \cap L^2(\R)^{\oplus 2}_+ 
    = H^2(\bS_\beta) \cap (H^2(\C_+) \oplus H^2(\C_+))
    = \Gamma(\Delta^{1/2}\res_{H^2(\C_+)}).\]
As $\Delta^{1/2}\res_{H^2(\C_+)}$ is a symmetric contraction,
we thus obtain a  closed embedding of
$H^2(\C_+)$ with the range $H^2(\bS_\beta)_+$. 
  We likewise obtain an isomorphism
  \[ H^2(\bS_\beta)_-
    = H^2(\bS_\beta) \cap (H^2(\C_-) \oplus H^2(\C_-))
    = \{ (\Delta^{-1/2}f, f) \: f \in H^2(\C_-)\}
    \cong H^2(\beta i - \C_+).\]
So we identify $H^2(\C_+)$ and $H^2(\beta i - \C_+)$
  with closed subspace of $H^2(\bS_\beta)$ satisfying
  \[ H^2(\bS_\beta) \cong 
H^2(\bS_\beta)_+ \oplus H^2(\bS_\beta)_- \cong 
H^2(\C_+) \oplus H^2(\beta i - \C_+).\]
We also note that
\[ \cD(\Delta^{1/2})
  = H^2(\C_+) \oplus \cD(\Delta^{1/2}\res_{H^2(\C_-)}) \]
leads to
\begin{equation}
  \label{eq:h2split}
 H^2(\bS_\beta) = \Gamma(\Delta^{1/2})
  = H^2(\C_+) \oplus \Gamma(\Delta^{1/2}\res_{H^2(\C_-)})
  \cong H^2(\C_+) \oplus H^2(\beta i - \C_-).
\end{equation}
Note that
\begin{equation}
  \label{eq:thetastrip}
 \theta(H^2(\C_+)) = H^2(i\beta - \C_+) \quad \mbox{ and } \quad
 \theta \Delta  = \Delta ^{-1} \theta
\end{equation}
holds for the restriction of $\Delta $ to $H^2(\C_+)$ and the restriction
of $\Delta ^{-1}$ to $H^2(\beta i - \C_+)$,
respectively.

As $\Delta $ is a symmetric contraction on $H^2(\C_+)$, the series
\[ \sum_{n = 0}^\infty (-1)^n \Delta ^n = (\1 + \Delta )^{-1} \]
converges strongly on $H^2(\C_+)$ to a positive contraction. Likewise
$\sum_{n = 1}^\infty (-1)^n \Delta ^{-n}$ 
converges to a bounded operator on $H^2(\beta i -\C_+)$. 

For $w \in \bS_\beta$, we have
$q_w \in H^2(\C_+)$ and
$q_{w-2\beta i} \in H^2(\beta i - \C_+)$.
Therefore the series expansion \eqref{eq:szego-exp2} takes the form 
\begin{equation}
  \label{eq:qseriesonstrip}
  Q_w
  = \sum_{n \geq 0} (-1)^n \Delta ^n q_w
  - \sum_{n \leq 0} (-1)^n \Delta ^{n} q_{w-2\beta i} 
  = \sum_{n \geq 0} (-1)^n \Delta ^n q_w
  +  \sum_{n < 0} (-1)^n \Delta ^{n+1} q_{w-2\beta i}. 
\end{equation}
The first series converges
in the Hilbert space $H^2(\C_+)$ and the second one in
$H^2(\beta i - \C_+)$. Hence the decomposition
\eqref{eq:h2split} shows that the series \eqref{eq:qseriesonstrip}
converges in $H^2(\bS_\beta)$. 

\begin{rem} Here is a more direct argument that leads to this
  series representation of the Szeg\"o kernel of the strip.
First we use the orthogonal decomposition 
\[ H^2(\bS_\beta) =H^2(\bS_\beta)_+ \oplus H^2(\bS_\beta)_- \]
from \eqref{eq:h2split} to write the Szeg\"o kernel functions
as
\[ Q_w = Q_w^+ + Q_w^-. \]

The evaluation functional $\ev_w \:  H^2(\bS_\beta) \to \C,
f \mapsto f(w)$, is represented on the closed subspace
$H^2(\bS_\beta)_+ = \Gamma(\Delta^{1/2}\res_{H^2(\C_+)})$ by $Q_w^+$.
For $f \in H^2(\C_+)$, identified with a closed subspace of $H^2(\bS_\beta)$
and $Q_w^+ = (r^+_w, \Delta^{1/2} r^+_w)$ for $r^+_w \in H^2(\C_+)$, we now have
\begin{align*}
  \la Q_w^{\C_+}, f \ra_{H^2(\C_+)}
  &= f(w) = \la Q_w^+, f \ra_{H^2(\bS_\beta)}
  = \la (r^+_w, \Delta^{1/2} r^+_w), (f, \Delta^{1/2} f) \ra\\
&  = \la r^+_w, f \ra_{H^2(\C_+)} + \la \Delta r^+_w, f \ra_{H^2(\C_+)}
  = \la (\1 + \Delta) r^+_w, f \ra_{H^2(\C_+)}.
\end{align*}
This implies that
\[ r^+_w
  = (\1 + \Delta)^{-1} Q_w^{\C_+} 
  = \sum_{k \geq 0} (-1)^k \Delta^k Q_w^{\C_+}
  = \sum_{k \geq 0} (-1)^k Q_{w + 2 k \beta i}^{\C_+},\]
where the series on the right converges in $H^2(\C_+)$.
This is the above series expansion of the kernel~$Q_w^+$.
For $Q_w^- = (r_w^-, \Delta^{1/2} r_w^-)$, we first observe that
\eqref{eq:thetastrip}, combined with $\theta(Q_w) = Q_{\beta i -w}$,
implies 
\[ Q_w^- = \theta(Q_{\beta i - w}^+).\]
So 
\begin{align*}
 r_w^-(z)
& = \theta(r_{\beta i - w}^+)(z)
  = r_{\beta i - w}^+(\beta i - z)
 = \sum_{k \geq 0}  \frac{(-1)^k}{2\pi}
    \frac{i}{\beta i - z + 2k \beta i + \beta i + \oline w} \\
  & = \sum_{k \geq 0} \frac{(-1)^k}{2\pi}
    \frac{i}{2(k+1)\beta i - z + \oline w}
    = \sum_{k \geq 0} \frac{(-1)^{k+1} }{2\pi}
    \frac{i}{z - \oline w - 2(k+1)\beta i} \\
  & = \sum_{k \geq 0} (-1)^{k+1} q_{w - 2(k+1)\beta i}
    = \sum_{k < 0} (-1)^k q_{w + 2k\beta i}.
\end{align*}
Now \eqref{eq:qseriesonstrip} is expressed as a series
converging in the Hilbert norm by writing
\[ Q_w = Q_w^+ + Q_w^- = (r_w^+ + r_w^-, \Delta^{1/2}(r_w^+ + r_w^-)).\] 
\end{rem}

\subsubsection*{An abstract setting}

To put this construction into a general framework, we
start with an abstract complex Hilbert space $\cH$, a 
unitary involution $\theta_0$ and a unitary
$1$-parameter group $(U_t)_{t \in \R}$ satisfying the relation
$\theta_0 U_t \theta_0 = U_{-t}$ for $t \in \R$.
We write for $\beta>0$
\[ U_t = e^{itP} \quad \mbox{ and } \quad \Delta := e^{-2\beta P}\]
and consider the graph 
\[ \cH^\beta := \Gamma(\Delta^{1/2}) \subeq \cH^{\oplus 2}.\]
We assume that $\ker(P) = \{0\}$ and write
$\cH_\pm$ for the positive and negative spectral subspace of $P$ in~$\cH$,
so that $\cH^\beta_\pm = \cH^\beta \cap  \cH_\pm^{\oplus 2}$ are the corresponding
spectral subspaces of $\cH^\beta$ for $P \oplus P$. Note that
\[ \theta(f_1, f_2) := (\theta_0 f_2, \theta_0 f_1) \]
defines a unitary involution on $\cH^\beta$ because
$f_2 = \Delta^{1/2} f_1$ implies
\[ \Delta^{1/2} \theta_0 f_2 = \theta_0 \Delta^{-1/2} f_2  = \theta_0 f_1.\]
On the domain of $\Delta^{1/2}$, 
the involution $\theta$ thus takes the form
\begin{equation}
  \label{eq:theta-form2}
 f \mapsto \Delta^{-1/2}\theta_0(f)\quad \mbox{ for }
 \quad f \in \cD(\Delta^{1/2}).
\end{equation}
This involution also satisfies 
$\theta U_t^\beta \theta = U_{-t}^\beta$ for the unitary
representation $(U_t^\beta)_{t \in \R}$ on $\cH^\beta$,
defined by $U^\beta_t(f_1, f_2) = (U_t f_1, U_t f_2)$ for $t \in \R$.

\begin{ex} Clearly, this generalizes the situation where
  \[ \cH = L^2(\R), \ \ 
    \cH_\pm = H^2(\C_\pm), \quad
    \cH^\beta = H^2(\bS_\beta) \quad \mbox{ and } \quad 
    (U_tf)(x)= f(x+ t), \ \ (\theta_0 f)(x) = f(-x).\] 
\end{ex}
\medskip

In the general context, the situation is most easily analyzed
for cyclic representations, which means that
$\cH_+$ is realized as a reproducing
kernel Hilbert space $\cH_{Q^{\C_+}} \subeq \cO(\C_+)$ on which
$(U_t F)(z) = F(z + t)$. 
As the projection
$p_+ \:  \cH^\beta_+ \to \cH_+$ is injective, we may also  consider
$\cH^\beta_+$ as a space of holomorphic functions on $\C_+$
for which point evaluations are continuous.
But this space has a different reproducing kernel~$Q^+$.
For the corresponding reproducing kernel spaces we then have
\[ \cH_{Q^+} \cong \cH^\beta_+ \into \cH_{Q^{\C_+}} \cong \cH_+.\]

We now determine $Q^+$ in terms of the kernel $Q^{\C_+}$.
For $f \in \cH^\beta_+$, considered as a holomorphic function
on $\C_+$ and $Q_w^+ = (r^+_w, \Delta^{1/2} r^+_w)$
for $r^+_w \in \cH_+$, we obtain as
for the Hardy space of the strip: 
\begin{align*}
  \la Q_w^{\C_+}, f \ra_{\cH_+} 
  &= f(w) = \la Q_w^+, f \ra_{\cH^\beta} 
    = \la (r^+_w, \Delta^{1/2} r^+_w), (f, \Delta^{1/2} f) \ra
    = \la (\1 + \Delta) r^+_w, f \ra_{\cH_+}. 
\end{align*}
This implies that
\[ r^+_w
  = (\1 + \Delta)^{-1} Q_w^{\C_+} 
  = \sum_{k \geq 0} (-1)^k \Delta^k Q_w^{\C_+}
  = \sum_{k \geq 0} (-1)^k Q_{w + 2 k \beta i}^{\C_+},\]
where the series on the right converges in $\cH_+$.
So we also obtain a series expansion
\[ Q^+(z,w)   = \sum_{k \geq 0} (-1)^k Q^{\C_+}(z, w + 2 k \beta i)\]
of $Q^+$, resembling
the one for the Hardy space.

To push the analogy further,
we observe that the unitary involution $\theta$ on $\cH_\beta$
induces a unitary operator $\cH^\beta_+ \to \cH^\beta_-$.
We also realize $\cH_-$ as a reproducing kernel Hilbert space
$\cH_{Q^{\C_-}}$ on the lower half plane, where
\[ \theta_0 \: \cH_+ \cong \cH_{Q^{\C_+}} \to \cH_{Q^{\C_-}} \cong \cH_-, \quad
  \theta_0(f)(z) := f(-z).\]
This means in particular that
\[ Q^{\C_-}(z,w) = Q^{\C_+}(-z,-w).\] 
In this picture \eqref{eq:theta-form2} implies that the unitary operator
$\theta \: \cH_{Q^+} \to \cH_{Q^-}$ is implemented by
\[ \theta \: \cO(\C_+) \supeq \cH_{Q^+} \cong \cH^\beta_+  \to
  \cH^\beta_- \cong \cH_{Q^-}
  \subeq \cO(\beta i - \C_+),
  \quad \theta(F)(z) = F(-(z - \beta i)) = F(\beta i - z).\]
We thus obtain a realization of $\cH^\beta$ as a
reproducing kernel Hilbert space on the strip $\bS_\beta$,
whose kernel has a series expansion in terms of
the functions
$Q_{w + 2 k \beta i}^{\C_+}$ on $\C_+$ and the functions
$\theta(Q_{w + 2 k \beta i}^{\C_+})$ on $\beta i - \C_+$. 

\begin{ex} \mlabel{ex:series}
  It is worthwhile to take a closer look at these
  expansions in the Fourier transformed picture.
  We consider a positive symmetric Borel measure\footnote{A Borel measure $\mu$ is symmetric when $\mu(-A)=\mu(A)$ for every Borel set $A\subseteq\R$}
  $\mu$ on $\R$ with $\mu(\{0\}) = 0$. 
On  $\cH := L^2(\R,\mu)$ we then have the following operators: 
    \[ (U_t f)(\lambda) = e^{it\lambda} f(\lambda), \quad
    (Hf)(\lambda) = \lambda f(\lambda), \quad
    (\theta_0 f)(\lambda) = f(-\lambda), \quad
    (\Delta f)(\lambda) = e^{-2\beta \lambda} f(\lambda).\]
  We then have $\cH_\pm = L^2(\R_\pm,\mu_\pm)$, where
  $\mu_\pm$ is the restriction of $\mu$ to $\R_\pm$. 

  \nin (a)  We assume that the measure 
  $\mu_+ := \mu\res_{[0,\infty)}$ has a
  Laplace transform 
  $\cL(\mu_+)$ which is finite on $(0,\infty)$. By symmetry of the measure $\mu$, the measure $\mu_- := \mu\res_{[-\infty,0)}$ has finite Laplace transform. 
  We consider the family
  \begin{equation}
    \label{eq:qws}
 q_w = \chi_{\R_\pm} e_{-i \oline w} \quad \mbox{ for } \quad
 \pm \Im w > 0,
  \end{equation}
 Then $q_w \in \cH_\pm$ for $w \in \C_\pm$, and this is equivalent to the finiteness condition of the Laplace transform. We
 thus obtain a positive definite kernel on $\C_+$ by
 \[ Q_+(z,w) = \la q_z, q_w \ra
 = \int_0^\infty e^{i\lambda(z - \oline w)}\, d\mu(\lambda)
 = \hat\mu_+(z - \oline w).\]
Next we observe that, for $w \in \C_+$, 
 \begin{equation}
   \label{eq:deltaq}
 (\Delta q_w)(\lambda)
 = e^{-2\beta \lambda} \chi_{\R_+}(\lambda)   e^{-i\lambda \oline w}
 = (\chi_{\R_+} e_{-i\oline w-2\beta})(\lambda)   
 = q_{w + 2 i \beta}(\lambda),
 \end{equation}
 i.e., that
 \begin{equation}
   \label{eq:deltaq2}
   \Delta q_w = q_{w + 2 i \beta}
   \quad \mbox{ for }  \quad w \in \C_+.
 \end{equation}
 We likewise obtain
 \begin{equation}
   \label{eq:deltaq3}
   \Delta^{-1} q_w = q_{w - 2 i \beta}\quad \mbox{ for }\quad w \in \C_-.
 \end{equation}
 Note that, for $w \in \bS_\beta$, 
 \[ w +  k 2i \beta \in
  \begin{cases}
    \C_+ & \text{ for } k \geq 0 \\
    \C_- & \text{ for } k <  0.
  \end{cases}
\]
With \eqref{eq:deltaq2} and \eqref{eq:deltaq3} we now get 
 \begin{align*}
   & \sum_{k \geq  0} (-1)^k  \Delta^k q_w
   + \sum_{k \leq  0} (-1)^k  \Delta^k q_{w - 2i \beta}
 = (\1 + \Delta)^{-1} q_w + (\1 + \Delta^{-1})^{-1} q_{w - 2i\beta}\\
&     = \chi_{\R_+} \frac{ e_{-i\oline w}} {1 + e_{-2\beta}} 
+ \chi_{\R_-} \frac{ e_{-i(\oline w + 2 i \beta)}} {1 + e_{2\beta}} 
   = \chi_{\R_+} \frac{ e_{-i\oline w}} {1 + e_{-2\beta}}
+ \chi_{\R_-} \frac{ e_{-i\oline w}} {1 + e_{-2\beta}}
= \frac{ e_{-i\oline w}} {1 + e_{-2\beta}}.
 \end{align*}

For the measure
 \[ d\nu(\lambda) := \frac{d\mu(\lambda)}{1 + e^{-2\beta\lambda}} \]
      we have
      $\cL(\nu)(y) < \infty$ for $0 < y < 2 \beta$. Therefore, by Theorem \ref{thm:7.1},
      the Fourier transform of $\nu$ defines the positive definite
      sesquiholomorphic kernel 
      $K(z,w) := \hat\nu(z-\oline w)$ on~$\bS_\beta$.
We write $\nu = \nu_+ + \nu_-$ with
\begin{equation}
  \label{eq:dnu+}
 d\nu_+(\lambda)
        = \frac{d\mu_+(\lambda)}{1 + e^{-2\beta \lambda}}
        = \sum_{k = 0}^\infty (-1)^k e^{-2k \beta\lambda}\, d\mu_+(\lambda)
\end{equation}
and 
\begin{equation}
  \label{eq:dnu-}
 d\nu_-(\lambda)
        = \frac{d\mu_-(\lambda)}{1 + e^{-2\beta \lambda}}
        = e^{2\beta\lambda}\frac{d\mu_-(\lambda)}{1 + e^{2\beta \lambda}}
        = \sum_{k = 0}^\infty (-1)^k e^{2(1 + k)\beta\lambda} \, d\mu_-(\lambda).   
\end{equation}
      Then the above calculation provides the series expansion
      \begin{align*}
        \hat\nu(z) 
&=\hat\nu_+(z) + \hat\nu_-(z) \\
&=\sum_{k \geq 0} (-1)^k \hat\mu_+(z + 2 k i \beta)
+ \sum_{k \geq 0} (-1)^k \hat\mu_-(z -  2 (k+1) i \beta) \\
&=\sum_{k \geq 0} (-1)^k \hat\mu_+(z + 2 k i \beta)
+ \sum_{k \geq 0} (-1)^k \hat\mu_+(-z +  2 (k+1) i \beta) \\
&=\sum_{k \geq 0} (-1)^k \hat\mu_+(z + 2 k i \beta)
+ \sum_{k \geq 0} (-1)^k \hat\mu_+(2\beta i -z + 2 k i \beta) \\
&=\hat\nu_+(z) + \hat\nu_+(2i \beta - z).
      \end{align*}
      So we also see that $\hat\nu(2i\beta - z) = \hat\nu(z)$,
      which corresponds to the relation
      $d\nu(-\lambda) = e^{-2\beta \lambda}\, d\nu(\lambda)$ from
      Lemma~\ref{lem:FTcosh}.

\nin (b) There is an interesting variant of this construction.
Under the same assumptions on $\mu$ as above,
we consider for $w \in \bS_\beta$ the series 
 \begin{align} \label{eq:bseries} 
   \sum_{k \in \Z} q_{w + 2 k i \beta}
   &     =
\sum_{k \geq  0}  \Delta^k q_w + 
     \sum_{k \leq  0}  \Delta^k q_{w - 2i \beta}
     = (\1 - \Delta)^{-1} q_w + (\1 - \Delta^{-1})^{-1} q_{w - 2i\beta} \notag \\
&  = \chi_{\R_+} \frac{ e_{-i\oline w}} {1 - e_{-2\beta}}
                                  + \chi_{\R_-} \frac{ e_{-i(\oline w + 2 i \beta)}} {1 - e_{2\beta}}
                                  = \chi_{\R_+} \frac{ e_{-i\oline w}} {1 - e_{-2\beta}}
+ \chi_{\R_-} \frac{ e_{-i\oline w}} {e_{-2\beta}-1}. 
 \end{align}

For the positive measure
\[ d\nu(\lambda)
  := \frac{d\mu_+(\lambda)}{1 - e^{-2\beta\lambda}}
  +  \frac{d\mu_-(\lambda)}{e^{-2\beta\lambda}-1}
  = \frac{\sgn(\lambda) d\mu(\lambda)}{1 - e^{-2\beta\lambda}} \]
we then find with $\nu = \nu_+ + \nu_-$ the series expansion
      \begin{align*}
        \hat\nu(z) 
        &=\hat\nu_+(z) + \hat\nu_-(z)
          =\sum_{k \geq 0}  \hat\mu_+(z + 2 k i \beta)
+  \sum_{k \geq 0}  \hat\mu_+(2i\beta -z + 2 k i \beta)
      \end{align*}
      and also in this case $\hat\nu(2i\beta - z) = \hat\nu(z)$.
\end{ex}

\begin{ex} \mlabel{ex:b.12}
(a) (Hardy space)   If $d\lambda$ denotes Lebesgue measure on $\R$ and
  $\mu = \frac{1}{{2\pi}}\, d \lambda,$ 
   then
\[ \hat\mu_+(z)
 = \int_0^\infty e^{i\lambda z}\, d\mu_+(\lambda)
 = \frac{1}{{2\pi}} \int_0^\infty e^{i\lambda z}\, d\lambda
 = \frac{1}{{2\pi}} \cdot \frac{1}{-iz} 
 = \frac{i }{{2\pi}\cdot z}\]
is the Szeg\"o kernel on $\C_+$, and on $\bS_{2\beta}$, the Fourier
transform $\hat\nu$ defines the Szeg\"o kernel of
$\bS_\beta$~\eqref{eq:cskernel-strip}.
The expansion in \eqref{eq:dnu+} and \eqref{eq:dnu-} thus spezializes to 
      \begin{align*}
\hat\nu(z) 
&=\sum_{k \geq 0} (-1)^k \hat\mu_+(z + 2 k i \beta)
 + \sum_{k \geq 0} (-1)^k \hat\mu_+(-z +  2 (k+1) i \beta) \\
&=\sum_{k \geq 0} (-1)^k \frac{i}{2\pi(z + 2 k i \beta)} 
 + \sum_{k \geq 0} (-1)^k \frac{i}{2\pi(-z + 2 (k+1) i \beta)}\\
&=\sum_{k \geq 0} (-1)^k \frac{i}{2\pi(z + 2 k i \beta)} 
 + \sum_{k \geq 0} (-1)^{k+1} \frac{i}{2\pi(z - 2 (k+1) i \beta)}\\
&=\sum_{k \geq 0} (-1)^k \frac{i}{2\pi(z + 2 k i \beta)} 
 + \sum_{k < 0}(-1)^{k} \frac{i}{2\pi(z + 2 k i \beta)}
 =\sum_{k \in \Z} (-1)^k \frac{i}{2\pi(z + 2 k i \beta)}.
      \end{align*}
This is our series expansion of the Szeg\"o kernel of $\cS_\beta$:
\[ Q(z,w) = \hat\nu(z - \oline w)
  =  \frac{1}{2\pi} \sum_{k \in \Z} (-1)^k \frac{i}{z -\oline w+ 2i k \beta}.\] 

\nin (b) (Bergman space)
Taking instead the symmetric measure
  $\mu = \frac{1}{(2\pi)^2}\, |\lambda| d \lambda$  on $\R$, 
the relation 
\[ \hat\mu_+(z)
= \int_0^\infty e^{i\lambda z}\, d\mu_+(\lambda)
= -i \frac{d}{dz} \int_0^\infty \frac{e^{i\lambda z}}{\lambda} d\mu_+(\lambda)
= -i \frac{d}{dz} \frac{1}{4\pi^2} \frac{i}{z}
=  -\frac{1}{4\pi^2} \frac{1}{z^2} \]
implies that $\hat\mu_+(z - \oline w)$ is  the Bergman
kernel on $\C_+$ (cf.\ Example~\ref{ex:qs}).

For the measure
\[ d\nu(\lambda) :=
  \frac{d\mu_+(\lambda)}{1 - e^{-2\beta\lambda}}
+  \frac{d\mu_-(\lambda)}{e^{-2\beta\lambda}-1}
= \frac{1}{(2\pi)^2}\frac{\lambda d\lambda}{1 - e^{-2\beta\lambda}} \]
whose Fourier transform defines the Bergman kernel by \eqref{eq:bergman},
the expansion in Example~\ref{ex:series}(b)
now leads to the series expansion 
\begin{align*}
\hat\nu(z) 
 &=\sum_{k \geq 0}  \hat\mu_+(z + 2 k i \beta)
 +  \sum_{k \geq 0}  \hat\mu_+(2i\beta -z + 2 k i \beta)
\\&= \sum_{k \geq 0}  -\frac{1}{4\pi^2} \frac{1}{(z + 2 k i \beta)^2}
  +  \sum_{k \geq 0}  -\frac{1}{4\pi^2} \frac{1}{(2i\beta -z + 2 k i \beta)^2}
  =-\frac{1}{4\pi^2} \sum_{k \in \Z} \frac{1}{(z + 2 k i \beta)^2}.
\end{align*}
So we obtain a series expansion of the Bergman kernel
of the strip $\bS_\beta$ in terms of translations and reflections
of the Bergman kernel of~$\C_+$.
\end{ex}

\section{A uniform perspective on some normal forms} 
\mlabel{sec:7}

In this section we discuss a unifying perspective
  on normal form results related closely to Hilbert spaces
  of holomorphic functions on $\bD$, $\C_+$ and $\bS_\beta$, respectively, 
  namely for isometries on Hilbert spaces,
    for one-parameter semigroups of isometries, 
  and for standard pairs $(U,\sV)$,  
  consisting of a standard subspace $\sV$ of a complex Hilbert space $\cH$
  and a  unitary one-parameter group $(U_t)_{t \in \R}$ on $\cH$ such that 
$U_t \sV \subeq \sV$ for $t \geq 0$ and $U_t = e^{itH}$ with 
$H \geq 0$. In particular, we explain how the corresponding classical
normal form results can be derived in a uniform fashion from the
Mackey--Stone--von Neumann Theorem on projective 
unitary representations of the product group $A \times \hat A$,
where $A \cong \T$ or $\R$ in our case.

More specifically, we consider the following structures:
\begin{itemize}
\item[$(\bD)$] A single isometry $S \: \cE_+ \to \cE_+$ on a Hilbert space
  $\cE_+$ (which can always be enlarged to a unitary operator 
$S \:  \cE \to \cE$ leaving the subspace $\cE_+$ invariant. 
If $(S,\cE_+)$ is regular in the sense that $\bigcap_n S^n \cE_+ = \{0\}$, 
then the Wold decomposition leads to an equivalence of $S$ with the 
unilateral shift on $\ell^2(\N_0,\cK)$, where $\cK = \cE_+ \cap
(S\cE_+)^\bot$ is a multiplicity space. 
Then we can realize $S$ as a multiplication operator 
on $H^2(\bD,\cK)$, given by $(Sf)(z) = z f(z)$,
and $(\N_0, \id) \to (H^\infty(\bD),\sharp),
n \mapsto z^n$, is a morphism of involutive
semigroups. 

\item[$(\C_+)$] A one-parameter semigroup of isometries 
$(S_t)_{t \geq 0}$ on a complex Hilbert space $\cE_+$. 
Here the Lax--Phillips Theorem (\cite[\S 4.4]{NO18}, \cite[Thm.~1]{LP64})
provides, under the regularity 
assumption  $\bigcap_{t > 0} S_t \cE_+ = \{0\}$, an equivalence 
with the unilateral shifts on $L^2(\R_+,\cK)$.
Applying the Fourier transform leads to
the realization  on the Hardy space
$H^2(\C_+,\cK)$ by  the multiplication operators 
$(S_t f)(z) = e^{itz} f(z)$, which define a morphism of
involutive semigroups $(\R_+, \id) \to (H^\infty(\C_+),\sharp)$. 

\item[$(\bS_\beta)$] {\it Standard pairs} $(U,\sV)$ consist of a standard subspace
$\sV$ of a complex Hilbert space $\cH$ and a 
unitary one-parameter group $(U_t)_{t \in \R}$ on $\cH$ such that 
$U_t \sV \subeq \sV$ for $t \geq 0$ and $U_t = e^{itH}$ with $H \geq 0$. 
By the Borchers--Wiesbrock Theorem
(\cite[\S 3.4]{NO17}, \cite{Bo92}, \cite{Wi93}), any standard pair 
defines an antiunitary positive energy representation of 
$\Aff(\R) \cong \R \rtimes \R^\times$ by 
\begin{equation}  \label{eq:affrep1}
  U(b,e^s) := U_b \Delta_\sV^{-is/2\pi} \quad  \mbox{ and }
  \quad 
U(0,-1) := J_\sV.
\end{equation}
If the regularity condition $\bigcap_{b >0} U_b \sV = \{0\}$ is satisfied, 
then this representation can be realized on $L^2(\R,\cK)$ via 
\[ (U(b,e^s) f)(x) = e^{ib e^{\frac{\pi x}{\beta}}} f(x + s) 
\quad  \mbox{ and } \quad  
 U(0,-1) f = J_\cK f,\] 
 where $J_\cK$ is a conjugation on $\cK$.
 Then $\sV$ corresponds to the standard subspace
 \[ H^2(\bS_\beta,\cK)^\sharp = \{ f \in H^\infty(\bS_\beta,\cK) \:
 f^\sharp = f\}, \quad f^\sharp(z) := J_\cK f(\beta i + \oline z),
 \quad z \in \bS_\beta. \] 
\end{itemize}

These normal form results can be derived from a common
source, which is the Stone--von Neumann--Mackey Theorem 
for locally compact abelian groups $A$ (cf.\ \cite[Thm.~A.VIII.6]{Ne99})
which classifies the irreducible projective unitary representations
$(U,\cH)$ of $A \times \hat A$ 
satisfying the canonical commutation relations:
\[ U_a U_\chi U_a^{-1} = \chi(a) U_\chi \quad \mbox{ for } \quad
  (a,\chi) \in A \times \hat A,\]
hence define a unitary representation of the group
\[ \Heis(A) := \T \times \hat A  \times A, \quad
  (z,\chi, a)  (z',\chi', a') =  (zz' \chi'(a), z + z', \chi + \chi').\]
\begin{itemize}
\item[($\bD$)] For the disc $\bD$,
  we consider the rotation action
  $(U_z f)(w) = f(zw)$, $z,w \in \T$, of $A := \T$ on $L^2(\T)$,
  which is complemented by the representation of $\Z\cong \hat A$ 
  on $L^2(\T)$ by the multiplication operators $(U_n f)(z) = z^n f(z)$.
  Both do not commute, but satisfy the commutation relations
  $U_z U_n U_z^{-1} = z^n U_n$, hence define a
  representation of the central extension 
$\Heis(\Z) \cong \T \times (\Z \times \hat\Z)$.
Starting with a unitary operator $U$ on a complex Hilbert space $\cE$, 
and a subspace $\cE_+$ invariant under~$U$, the regularity conditions 
$\bigcap_{n \in \N} U^n\cE_+~=~\{0\}$ and $\bigcap_{n \in \N} U^{-n}\cE_+^\bot = \{0\}$ 
imply the existence of a spectral measure
$P$ on $\Z \cong \hat \T$ satisfying 
$P(\{n\}) = U^n(\cE_+ \cap (U\cE_+)^\bot)$, hence also
$P(\N_0) = \cE_+$.
We thus obtain a unitary representation
of $\T$ satisfying the canonical commutation relations
and the Stone--von Neumann--Mackey Theorem applies.

\item[$(\C_+)$] For the upper half plane, the multiplication operators 
$(S_t)_{t \in \R}$ define a unitary 
representation of $\R$ on $L^2(\R)$, which is complemented by the 
translation action of $\R$ to a projective
representation of $\R \times \hat\R \cong \R^2$.
Starting with a unitary one-parameter group
$(U_t)_{t \in \R}$ on $\cE$
a subspace $\cE_+$ invariant under $(U_t)_{t \geq 0}$, the regularity condition 
implies the existence of a spectral measure on $\R \cong \hat \R$ 
assigning $\cE_+$ to $[0,\infty)$. We thus obtain a unitary representation
of $\R$ satisfying the canonical commutation relations
and the classical Stone--von Neumann Theorem on the representations
of the Heisenberg group $\Heis(\R)$ applies
(see \cite[Thm.~4.4.1]{NO18} for details).
\item[$(\bS_\beta$)] For the strip we obtain an antiunitary representation of 
  $\Aff(\R)$ on $L^2(\R)$ for which the inclusion of
  \[ H^2(\bS_\beta)^\sharp := \{ f \in H^2(\bS_\beta) \: f^\sharp = f\},
  \quad f^\sharp(z) = \oline{f(\beta i + \oline z)} \]
  by non-tangential boundary values on $\R$ yields a standard subspace
  whose modular group acts by translations
  \[ (\Delta_\sV^{-is/2\beta} f)(x) = f(x + s) \quad \mbox{ and } \quad
  J_\sV f = \oline f.\]
  With $(U_b f)(x) := e^{ib e^\frac{\pi x}\beta} f(x)$ we thus obtain via
  \eqref{eq:affrep1} an antiunitary positive energy representation of
  $\Aff(\R)$ on $L^2(\R)$.

  If, conversely, $(U,\cH)$ is an antiunitary positive energy
  representation of $\Aff(\R)$ on $\cH$, then there exists a unique
  standard subspace $\sV$ with $U(0,-1) = J_\sV$ and
  $U(0,e^s) = \Delta_\sV^{-is/2\beta}$ for $s \in \R$.
  Positive energy now implies that $U_s \sV \subeq \sV$ for $s \geq 0$,
  and if the regularity condition is satisfied, then
  $A := \partial U(1,0)$ is strictly positive, so that
  $Q := \log A$ makes sense and
  $\pi(t,s) := e^{itQ} U(0,e^s)$ defines a projective unitary representation
  of $\R \times \R$, satisfying the canonical commutation relations.
  Hence we obtain the asserted  normal form situation 
  with a representation on $L^2(\R_+,\cK)$ (cf.\ \cite[Prop.~2.38]{NO17}),
  respectively $L^2(\R,\cK) \supeq H^2(\bS_\beta)^\sharp$ in the spectral
  picture (see \eqref{eq:h2sharp} and also \cite{ANS23} for more details). 
\end{itemize}

\begin{ex} In the same setting as before, the operator $R = e^{-\beta P}$ on $L^2(\R)$, where $Pf=if'$ for $f\in L^2(\R)$, satisfies 
with the unitary operators
\[ (U_t f)(x) = e^{itx} f(x) 
  \quad \mbox{ and } \quad
  (V_s f)(x) = (e^{isP}f)(x) = f(x + s) \quad \mbox{ for } \quad s,t \in \R,
  x \in \R,\]
the relation 
\[   V_s U_t = e^{its} U_t V_s \quad \mbox{ for } \quad s,t \in \R.\]
This leads in the form
$U_t V_s U_t^{-1} = e^{-its} V_s$ with measurable functional calculus to
$U_t R U_t^{-1} = e^{t\beta}R$, i.e., to 
\begin{equation}
  \label{eq:R-rel}
  R U_t = e^{-t \beta} U_t R \quad \mbox{ for } \quad t \in \R.
\end{equation}
\end{ex}

More in general, we have the following

\begin{prop}
Let $(U_t)_{t \in \R}$ be a unitary one-parameter group
on the complex Hilbert space $\cH$ and
$R = R^*$ be a selfadjoint operator satisfying
\begin{equation}
  \label{eq:r-rel}
 R U_t = e^{-\beta t} U_t R \quad \mbox{ for } \quad t \in \R.   
\end{equation}
Then the following assertions hold:
\begin{itemize}
\item[\rm(a)] The one-parameter group $(U_t)_{t \in \R}$
  preserves the spectral subspaces $\cH_-, \cH_0, \cH_+$ of $R$ corresponding to
  $(-\infty,0), \{0\}, (0,\infty)$.
\item[\rm(b)] 
  On $\cH_+$ the operator $R$ is positive, so that the
    unitary one-parameter group
  $(R^{is})_{s \in \R}$ is defined and satisfies the canonical
commutation relations 
\begin{equation}
  \label{eq:r-rel3}
 U_t R^{is} U_t^{-1} = e^{i\beta st} R^{is} \quad \mbox{ for } \quad s,t \in \R.    
\end{equation}
Thus, there exists a unitary operator $\Phi_+ \: \cH_+ \to L^2(\R, \cK_+)$ intertwining
$U_t$ and $R$ with 
\begin{equation}
  \label{eq:dag2}
 (\tilde U_t f)(x) = e^{itx} f(x)  \quad \mbox{ and } \quad
 (\tilde R^{is} f)(x) = f(x - s \beta).
\end{equation}
In a similar way, on $\mathcal{H}_-$ the operator $R_-:=-R$ is positive. Thus the unitary one-parameter group $(R_-^{is})_{s\in\mathbb{R}}$ is defined and satisfies the canonical commutation relations as in \eqref{eq:r-rel3}. Therefore, there exists a unitary operator $\Phi_- \: \cH_- \to L^2(\R, \cK_-)$ that verify similar properties as in \eqref{eq:dag2}.\\
If $\mathrm{ker}(R)=\{0\}$, then $\Phi:=\Phi_+\oplus\Phi_-:\mathcal{H}\to L^2(\mathbb{R},\mathcal{K}_+)\oplus L^2(\mathbb{R},\mathcal{K}_-)$ is a unitary
intertwining operator.
\end{itemize}
\end{prop}

\begin{prf} Rewriting \eqref{eq:r-rel}, we obtain
\begin{equation}
  \label{eq:r-rel2}
 U_t R U_t^{-1} = e^{\beta t} R \quad \mbox{ for } \quad t \in \R.   
\end{equation}
This relation implies (a), so that
$\cH$ decomposes into the $U_\R$-invariant subspaces
\[ \cH = \cH_+ \oplus \cH_0 \oplus \cH_-, \]
where $\cH_0 = \ker(R)$ and $\cH_\pm$ are the positive/negative spectral
subspaces of $R$.

On $\cH_+$ the operators $(R^{is})_{s \in \R}$ satisfy the canonical
commutation relations \eqref{eq:r-rel3}.
Therefore the Stone--von Neumann Theorem implies the existence
of a Hilbert space $\cK_+$ and of
a unitary operator
\[ \Phi_+ \: \cH_+ \to L^2(\R, \cK_+), \]
intertwining
$U_t$ and $R$ with \eqref{eq:dag2}.
On $\cH_-$, the operator $R$ is negative.
Thus $R_-=-R$ is a positive operator on $\cH_-$ and, as before, $(R_-^{is})_{s\in\R}$ verifies \eqref{eq:r-rel3}. Thus, there exists a unitary operator $\Phi_- \: \cH_- \to L^2(\R, \cK_-)$ with the same properties as $\Phi_+$.
 If $\mathrm{ker}(R)=\{0\}$, then the direct sum $\Phi_+\oplus\Phi_-$
is again a unitary operator. 
\end{prf}
  
\begin{rem} (A one-parameter group preserving $\Gamma(R)$) 
  In the context of the preceding proposition,
  the graph $\Gamma(R)$ is a closed subspace of $\cH \oplus \cH$.
It is invariant under the non-unitary one-parameter group
\[ W_t := U_t \oplus e^{-\beta t} U_t \]
because
\[ W_t (f, Rf)
  = (U_t f, e^{-\beta t} U_t R f) 
  = (U_t f, R U_t f) \]
for $f \in \cD(R)$. Moreover, \eqref{eq:r-rel2} implies that
$U_t \cD(R) = \cD(R)$ for $t \in \R$. Further, if $R$ is a positive operator, then $\cH=\cH_+$, $\cK:=\cK_+$ and $\Phi=\Phi_+:\cH\to L^2(\R,\cK)$, thus the isometry  
\[ \Phi \oplus \Phi \: \cH \oplus \cH \to L^2(\R,\cK) \oplus L^2(\R,\cK),\]
maps the subspace $\Gamma(R)$ to the Hardy space
$H^2(\bS_\beta, \cK)$ (cf.\ Example~\ref{ex:5.12}).
\end{rem}

\section{Reflection positivity in the parameter} 
  \mlabel{sec:8}

  In this section we discuss the interesting observation that
  the three generating subsemigroups in $H^2(\C_+)$, $H^2(\bD)$ and $H^2(\bS_\beta)$
exhibit also some reflection positivity properties in their
parameter. Presently, we do not have any conceptual explanation for this phenomenon.

\begin{itemize}
\item On $\C_+$ we have two semigroups isomorphic to 
$\R_+$, one acts by multiplications and the other $i\R_+$ by 
translations. Evaluation of the exponentials yields 
\[ \R_+ \times i \R_+ \to \R^\times, \quad 
(t,i\lambda) \mapsto e_{it}(i\lambda) = e^{-\lambda t} \] 
which is bihomomorphic, i.e., a semigroup homomorphism in both
  arguments. This map extends naturally to 
\[ \R_+ \times \C_+ \to \oline{\bD}, \quad 
(t,z) \mapsto e_{it}(z) = e^{itz} \] 
which defines a $*$-homomorphism $(\R_+,\id) \to (H^\infty(\C_+), \sharp)$ for which 
\[ \C_+ = \{ z \in \C \: (\forall t > 0)\ |e^{itz}| < 1 \}.\]  
\item On $\bD$ we have two semigroups: the semigroup $(z^n)_{n \in \N_0}$ 
acts by multiplications on $H^2(\bD)$ and the semigroup 
$S := \{r \in \R^\times  \: |r| < 1\} = \{\pm 1\} \exp(-\R_+) \subeq \T_\C = \C^\times$ acts 
geometrically. Again, the map 
\[(\N_0,+) \times ((-1,1),\cdot) \to 
(\R_+,\cdot), \quad  (n,\lambda) \mapsto \lambda^n  \] 
is a bihomomorphism. It extends to a map 
\[ \N_0 \times \bD \to \oline{\bD}, \quad 
(n,z) \mapsto z^n \] 
which defines a $*$-homomorphism $(\N_0, \id) \to (H^\infty(\bD),\sharp)$ 
and 
\[ \bD = \{ z \in \C \: (\forall n \in \N) |z^n| < 1 \}.\]  

For each $n \in \N$, the power function $p(\lambda) := \lambda^n$ 
on the subsemigroup $S$ 
extends to a reflection positive function on 
$(\R^\times, \tau)$, where $\tau(a) = a^{-1}$. In fact, identifying 
$(\R^\times, \tau)$ with $(\R \times \{\pm 1\}, -\id)$ via the exponential 
function, the reflection positive extension is given by 
\[ p_n(t,\eps) = \eps^n e^{-n|t|}.\] 

\item For $\bS_\beta$ the semigroup acting geometrically has to 
be replaced by the domain 
\[ \T_{\beta,+} = \{ [iy] \: 0 < y < \beta/2\} \subeq \T_\beta,\] 
which only ``acts'' on the lower half of the strip. We also have an  
evaluation map 
\[ \R_+ \times \T_{\beta,+} \to \R,\quad 
(t, iy) \mapsto  c_t(iy) = 
\frac{e^{-t y} + e^{-t(\beta -y)}}{1 + e^{-t\beta}} 
= \frac{\cosh\big(t\big(\frac{\beta}{2}-y\big)\big)}
{\cosh(t\beta/2)}, \] 
but this map is not homomorphic in any argument. It 
extends to a map 
\[ \R \times \bS_\beta \to \C, \quad 
(t, z) \mapsto  c_t(z) = 
\frac{e^{it z} + e^{-t \beta} e^{-it z}}{1 + e^{-t\beta}} \] 
that defines a map 
$\R \to H^\infty(\bS_\beta)$. In view of Lemma~\ref{lem:ct}, we have 
\[ \bS_\beta = \{ z \in \C \: (\forall t> 0) \ |c_t(z)| < 1 \}. \]
Note that the functions $c_t(z)$ are not very far from being multiplicative 
in $z$: for each $\beta$ it is a linear combination of two 
homomorphisms.
\end{itemize}

\appendix 

\section{Appendix}

\subsection{Total families in $H^\infty$} 

In this subsection we briefly discuss the generating families
of $H^\infty$ for $\bD, \C_+$ and $\bS_\beta$ with an emphasis on
the strip. 

\begin{lem} \mlabel{lem:polweakdense}
  \begin{itemize}
  \item[\rm(a)] The subspace $\C[z] \subeq H^\infty(\bD)$ 
of polynomials is dense with respect to the weak topology. 
  \item[\rm(b)] The one-parameter semigroup 
$(e_{it})_{t > 0}$ spans a weakly dense subspace of $H^\infty(\C_+)$. 
\item[\rm(c)]   The functions $(c_t)_{t \geq 0}$, defined by 
\[ c_t = \frac{e_{it} + e^{-\beta t} e_{-it}}{1 + e^{-\beta t}} 
, \quad t \geq  0 \] 
 span a 
weakly dense subspace of $H^\infty(\bS_\beta)^\tau$, 
the fixed points under the holomorphic involution $\tau(z) = \beta i - z$. 
  \end{itemize}
\end{lem}

\begin{prf} For (a) and (b) we refer to \cite[Lemma~B.6]{ANS22}. 
To verify (c), we observe that, 
for $t \in \R$ and $0 < y < \beta$, we have 
\[  e_{it}(\tau(iy)) = e_{it}(i(\beta - y)) 
= e^{-t\beta} e_{-it}(iy),\quad \mbox{ hence } \quad 
e_{it} \circ \tau = e^{-t\beta} e_{-it}.\] 
Therefore the projection 
\[ p_\tau \: H^\infty(\bS_\beta) \to H^\infty(\bS_\beta)^\tau, \quad 
p_\tau(f) := \frac{f + f \circ \tau}{2} \] 
maps $e_{it}$ to a positive multiple of $c_t$ 
and 
\[ p_\tau(e_{-it}) 
= \frac{e_{-it} + e^{t\beta} e_{it}}{2} 
= e^{t\beta} \frac{e^{-t\beta} e_{-it} + e_{it}}{2} 
= e^{t\beta} p_\tau(e_{it}).\] 
As $p_\tau$ is surjective 
and weakly continuous, it  
suffices to recall from \cite[Lemma~B.6(c)]{ANS22} 
that the functions  $(e_{it})_{t \in \R}$ 
span a dense subspace of $H^\infty(\bS_\beta)$. 
\end{prf}

\begin{lem}\mlabel{lem:ct}
For $t \in \R$, the function  
\[ c_t = \frac{e_{it} + e^{-\beta t} e_{-it}}{1 + e^{-\beta t}} 
\quad \mbox{\ on\ $\bS_\beta$ satisfies  } \quad 
c_t^{\sharp} = c_t \] 
and defines an  element of $H^\infty(\bS_\beta)$ with 
$\|c_t\|_\infty = 1$.
These functions  satisfy 
  \begin{equation}
    \label{eq:stripcahr}
 \bS_\beta = \{ z \in \C \: (\forall t > 0)\ |c_t(z)| < 1 \}.
\end{equation}
\end{lem}

\begin{prf} As  $c_t(0)= 1$ and 
$|c_t(x)|, |c_t(\beta i + x)| \leq 1$ for $x \in \R$, 
we derive  $\|c_t\|_\infty  = 1$ from Hadamard's Three Lines Theorem 
(\cite[Thm.~12.8]{Ru86}). 
Therefore, by the Maximum Modulus Theorem, we have 
$|c_t(z)|< 1$ for $z \in \bS_\beta$, so setting
\[M := \{ z \in \C \: (\forall t > 0)\ |c_t(z)| < 1 \},\]
we have \(\bS_\beta \subeq M\).
To verify that \(M \subeq \bS_\beta\), we observe that, for $z = x + i y$, we have  
\[ c_t\Big(z + \frac{i\beta}{2}\Big) 
= \frac{\cosh(i tx - ty)}{\cosh\big(\frac{t\beta}{2}\big)}.\] 
Next we observe that 
\[ \cosh(x + i y) 
= \cosh(x) \cos(y)  + i \sinh(x) \sin(y), \] 
leads to 
\begin{align*}
 |\cosh(x + i y)|^2
&= \cosh^2(x) \cos^2(y) + \sinh^2(x) \sin^2(y)\\
&= (1 + \sinh^2(x)) \cos^2(y) +\sinh^2(x) \sin^2(y)\\
&= \sinh^2(x) +\cos^2(y) = \cosh^2(x) - \sin^2(y) \leq \cosh^2(x). 
\end{align*}
We thus obtain 
\[ 
\Big|c_t\Big(z + \frac{i\beta}{2}\Big)\Big|^2 
=  \frac{|\cosh(i tx - ty)|^2}{\cosh\big(\frac{t\beta}{2}\big)^2}
=  \frac{\cosh^2(ty) - \sin^2(tx)}{\cosh\big(\frac{t\beta}{2}\big)^2}.\]
For $t \to \infty$, the asymptotics of this expression is 
\[ \frac{\cosh^2(ty) - \sin^2(tx)}{\cosh\big(\frac{t\beta}{2}\big)^2} 
\approx \frac{\cosh^2(ty)}{\cosh\big(\frac{t\beta}{2}\big)^2} 
\approx \frac{e^{2t|y|}}{e^{t\beta}} = e^{t(2|y| - \beta)}.\] 
We conclude that, for $|y| > \beta/2$, there exists a $t > 0$ 
for which this expression is $> 1$. 

It remains to check the cases \(y = \pm \frac \beta 2\). For \(y = -\frac \beta 2\) we have \(z + \frac{i\beta}{2} = x \in \R\), which implies
\[c_t\Big(z+ \frac{i\beta}{2}\Big) = c_t(x) = \frac{e^{itx}+e^{-\beta t}e^{-itx}}{1+e^{-\beta t}}.\]
For \(x=0\) this is 
\[c_t\Big(z+ \frac{i\beta}{2}\Big) = c_t(x) = 1 \quad \mbox{ for }  \quad
  t \in \R\]
and for \(x \neq 0\), setting \(t := \frac{2\pi}{|x|}\), we get
\[c_t\Big(z+ \frac{i\beta}{2}\Big) = c_t(x) = \frac{e^{2\pi i \sgn(x)}+e^{-\beta t}e^{-2\pi i \sgn(x)}}{1+e^{-\beta t}} = \frac{1+e^{-\beta t}}{1+e^{-\beta t}} = 1.\]
For \(y = \frac \beta 2\) we use the relation
$c_t^\sharp = c_t$ to obtain the corresponding conclusion.
This completes the proof.
\end{prf}

\subsection{An application of Poisson summation} 
\mlabel{app:c.2}

\begin{lem} \mlabel{lem:poisson} 
For $\lambda, \beta > 0$ the function 
$\psi_\lambda \:  \R \to \R, 
\psi_\lambda(x) := \frac{1}{\pi}\frac{\lambda}{\lambda^2 + x^2}$ 
satisfies 
\[ 
\sum_{k \in \Z} \psi_{\frac{\beta\lambda}{2\pi}}(k) e^{i k x \frac{2\pi}{\beta}}= \frac{ e^{-\lambda x} + e^{-\lambda (\beta-x)}}
{1 - e^{-\lambda \beta}} \quad \mbox{ for }\quad 0 \leq x \leq \beta. \]
\end{lem}

The left hand side of this equation is a $\beta$-periodic function,
so that this identity describes the Fourier expansion of the
$\beta$-periodic extension of the function on the right. 

\begin{prf} {\bf Step 1:} For any 
function $f \in L^1(\R)$ whose Fourier transform is summable 
on $\Z$, the Poisson summation formula (\cite[Ex.~3.6.4]{DE09}) 
asserts that we have almost 
everywhere on $\R$ the identity
\[ \sum_{k \in \Z} f(x + k) = \sum_{k \in \Z} \tilde f(k) e^{2\pi i kx}, 
  \quad \mbox{ where } \quad \tilde f(y) := \int_\R f(x) e^{-2\pi i yx}\, dx.\]
Below we shall need the rescaled formula, which leads 
for $f_\beta(x) := f(\beta x)$ to $\tilde f_\beta(y) = \frac{1}{\beta} \tilde f(y/\beta)$, so that 
\begin{equation}
  \label{eq:beta-poisson}
\sum_{k \in \Z} f(x + \beta k) 
= \sum_{k \in \Z} f_\beta\Big(\frac{x}{\beta} + k\Big) 
= \sum_{k \in \Z} \tilde f_\beta(k) e^{\frac{2\pi}{\beta} i kx}
= \frac{1}{\beta} \sum_{k \in \Z} \tilde f\Big(\frac{k}{\beta}\Big)
e^{\frac{2\pi}{\beta} i kx}.
\end{equation}

\nin {\bf Step 2:} For any $s > 0$, the functions $\psi_\lambda$ satisfy
\begin{equation}
  \label{eq:rescale}
 \psi_{\lambda/s}(x) = s \psi_\lambda(sx) \quad \mbox{ for } \quad x \in \R.
\end{equation}
From Proposition~\ref{prop:3.1}, we recall that, for $\lambda>0$,  
\[\phi_\lambda(y) = e^{-\lambda|y|} = \int_\R e^{-iyx} P_{i\lambda}(x)\, dx = \int_\R e^{-iyx} \psi_{\lambda}(x)\, dx\]
so that 
\[\phi_\lambda(2\pi y) = \tilde{\psi_\lambda}(y), \qquad y \in \R.\] 
This leads with \eqref{eq:rescale} to 
\begin{equation}
  \label{eq:fourel}
\widetilde{\phi}_\lambda(x) 
= \int_\R e^{-2\pi i xy} \tilde\psi_\lambda\big(\frac{y}{2\pi}\big)\, dy
= 2\pi \int_\R e^{-4\pi^2 i xy} \tilde\psi_\lambda(y)\, dy
= 2\pi \psi_\lambda(2\pi x) = \psi_{\frac{\lambda}{2\pi}}(x).
\end{equation}

\nin {\bf Step 3:} We can now apply the Poisson summation formula to the 
functions $\phi_\lambda$ as follows. 
By \eqref{eq:fourel}, we have 
\[ \frac{1}{\beta} \widetilde{\phi}_\lambda\Big(\frac{k}{\beta}\Big)
= \frac{1}{\beta} \psi_{\frac{\lambda}{2\pi}}\Big(\frac{k}{\beta}\Big)
=  \psi_{\frac{\lambda\beta}{2\pi}}(k).\] 
We thus obtain with the Poisson summation formula \eqref{eq:beta-poisson} 
\begin{align*}
&\sum_{k \in \Z} \psi_{\frac{\beta\lambda}{2\pi}}(k) e^{i k x \frac{2\pi}{\beta}} 
 = \sum_{k \in \Z} \frac{1}{\beta} \widetilde{\phi}_\lambda\Big(\frac{k}{\beta}\Big)
e^{i k x \frac{2\pi}{\beta}} 
= \sum_{k \in \Z} \phi_\lambda(x + k \beta)
= \sum_{k \in \Z} e^{-\lambda |x+k\beta|}.
\end{align*}
For $0 \leq x < \beta$, the right hand side can be evaluated to 
\begin{align*}
 \sum_{k \in \Z} e^{-\lambda |x+k\beta|}
&= \sum_{k = 0}^\infty e^{-\lambda (x+k\beta)} + \sum_{k = 1}^\infty e^{\lambda (x-k\beta)}
= e^{-\lambda x} \sum_{k = 0}^\infty (e^{-\lambda \beta})^k 
 + e^{\lambda x} \sum_{k = 1}^\infty (e^{-\lambda \beta})^k \\
&= e^{-\lambda x} \frac{1}{1 -e^{-\lambda \beta}} 
 + e^{\lambda x} \frac{e^{-\lambda \beta}}{1 - e^{-\lambda \beta}} 
= \frac{ e^{-\lambda x} + e^{-\lambda (\beta-x)}}
{1 - e^{-\lambda \beta}}. 
\end{align*}
Eventually, we get 
\[ 
\sum_{k \in \Z} \psi_{\frac{\beta\lambda}{2\pi}}(k) e^{i k x \frac{2\pi}{\beta}}= \frac{ e^{-\lambda x} + e^{-\lambda (\beta-x)}}
{1 - e^{-\lambda \beta}}. \qedhere\]
\end{prf}

\subsection{A direct derivation of the Szeg\"o kernel}
\mlabel{subsec:app-cz-kernel}

In this subsection we calculate the Szeg\"o kernels directly
in terms of scalar products of $L^2$-functions and Fourier transforms. 

\subsubsection*{The disc $\bD$}

For the circle $\T = \partial \bD$, the Hardy space
$H^2(\bD)$ can be identified with the space of all Fourier series 
\[ \hat f(z) = \frac{1}{\sqrt{2\pi}}
\sum_{n = 0}^\infty f_n z^n \quad \mbox{ with } \quad
\sum_{n=0}^\infty |f_n|^2 < \infty.\]
Note that the map $\ell^2(\N_0) \to H^2(\bD), f = (f_n) \mapsto \hat f$
is isometric. Here
\[ \hat f(e^{i\theta}) =\frac{1}{\sqrt{2\pi}}
\sum_{n=0}^\infty f_n e^{in\theta} \quad \mbox{ with } \quad
f_n = \frac{1}{\sqrt{2\pi}}\la e_n, \hat f \ra
\quad \mbox{ and }\quad
e_n(\theta) = e^{in\theta}.\]
Therefore
\[ \hat f(z) = \la \hat{q_z}, \hat f \ra \quad \mbox{ for } \quad
\hat{q_z} = \frac{1}{2\pi} \sum_{n = 0}^\infty  \oline z^n e_n,
\quad q_z \in L^2(\T).\]
We conclude that the Szeg\"o kernel of $\bD$ is given with
$(q_z)_n = \frac{1}{\sqrt{2\pi}} \oline z^n$ by
\[ Q(z,w) =  \la \hat{q_z}, \hat{q_w} \ra = \hat{q_w}(z)
= \frac{1}{2\pi} \sum_{n = 0}^\infty  \oline w^n z^n
= \frac{1}{2\pi}\frac{1}{(1 - z\oline w)}.\] 

\subsubsection*{The upper half plane $\C_+$}

For the upper half plane $\C_+$, the Hardy space 
$H^2(\C_+)$ consists of all Fourier transforms
$\hat f$ of $f \in L^2(\R_+)$, where
\[ \hat f(z) =   \frac{1}{\sqrt{2\pi}}
\int_\R e^{iz \lambda} f(\lambda)\, d \lambda
= \frac{1}{\sqrt{2\pi}}\la e_{-i\oline z}, f \ra_{L^2(\R_+)}. \]

If $\hat{q_z}$ is the Szeg\"o kernel of $\C_+$, by definition we have
$$\hat f(z)=\langle \hat{q_z}, \hat f\rangle_{H^2(\C_+)}=\langle q_z, f\rangle_{L^2(\R_+)}=\int_0^\infty \overline{q_z(\lambda)}f(\lambda)d\lambda$$

Therefore, comparing the two previous formulas, we have
\[ q_z(\lambda)=\frac{e^{-i\bar z \lambda}}{\sqrt{2\pi}}, \quad
\mbox{ i.e.,} \quad  q_z = \frac{e_{-i\oline z}}{2\pi}.  \]

For the measure $d\nu(\lambda) := \frac{1}{2\pi}d\lambda$ on $\R_+$,
we have by Example~\ref{ex:b.12}(a),
\[ \hat\nu(z) = \frac{1}{2\pi} \frac{i}{z} \quad \mbox{ for } \quad
  z \in \C_+,\]
so that
\[ \hat{q}(z,w)=\hat{q_w}(z)=\frac1{\sqrt{2\pi}}\int_0^\infty e^{i\lambda z}q_w(\lambda)d\lambda=\int_0^\infty e^{i(z - \oline w)\lambda}\, \frac{d \lambda}{2\pi}
= \hat\nu(z - \oline w) = \frac{i}{2\pi(z - \oline w)}.\]
Therefore the Szeg\"o kernel of $\C_+$ is 
\[ Q(z,w) =\la \hat{q_z},\hat{q_w}\ra_{H^2(\C_+)} 
= \oline{q_w}(z) =  \frac{1}{2\pi}\frac{i}{z - \oline w}.\]

We observe that $2\pi\nu$ is the Lebesgue measure $d\lambda$, which implements the isometric isomorphism between $H^2(\C_+)$ and $L^2(\mathbb{R}_+)$ given by the Fourier transform $L^2(\R_+)\ni f\mapsto \hat{f} \in H^2(\C_+)$.

\subsubsection*{The strip $\bS_\beta$}

To derive the Szeg\"o kernel of the strip
$\bS_\beta$, we first investigate the Hardy space $H^2(\mathbb{S}_\beta)$. 
\smallskip

For the strip $\bS_\beta$, we realize $H^2(\bS_\beta) \subeq
L^2(\R)^{\oplus 2}$ as the graph of the operator
$(Rf)(x) = f(x + \beta i)$ which maps the lower boundary values of a function $F\in H^2(\bS_\beta)$ to its upper boundary
values (cf.~Example~\ref{ex:5.12}(a) and
\cite[App.~A]{Ta15}). 
If $\mathcal{K}$ is defined by
$$\mathcal{K}:=\{f\in L^2(\mathbb{R}):e_{-\beta}f\in L^2(\mathbb{R})\}\subset L^2(\mathbb{R}),$$
by \cite[Prop~5.1]{Go69}, 
$H^2(\mathbb{S}_{\beta})$ consists of the Fourier transforms $\widehat{f}$ of functions $f\in\mathcal{K}$ corresponding to lower boundary values of functions $F$ in the Hardy space $H^2(\mathbb{S}_{\beta})$. 

Let $M_{e_{-\beta}}$ be the multiplication operator by $e_{-\beta}$  on $\mathcal{K}=\mathrm{Dom}(M_{e_{-\beta}})$ \cite[App.~A]{Ta15}. Then $\mathcal{K}$ can be
realized as a graph of the operator $M_{e_{-\beta}}$ in the following way 
\begin{equation}
  \label{eq:graphreal} 
  \mathcal{K}\to\Gamma(M_{e_{-\beta}}),\quad f\mapsto (f, e^{-\beta x}f)\subseteq L^2(\mathbb{R})^{\oplus 2}
\end{equation}
Its Fourier transform $\widehat{M_{e_{-\beta}}}=R$ acts as a vertical translation by $i\beta$. Therefore, $H^2(\mathbb{S}_{\beta})$ is isomorphic to the graph of the operator $\widehat{M_{e_{-\beta}}}$ in the following way
\[ H^2(\mathbb{S}_{\beta})\to\Gamma(\widehat{M_{e_{-\beta}}}),\quad
  \widehat{f}\mapsto (\widehat{f}, \widehat{M_{e_{-\beta}}}\widehat{f}), \quad
F \mapsto (F\res_\R, F\res_{\R + \beta i}). \]

For $\widehat{f}\in H^2(\mathbb{S}_{\beta})$, $z\in\mathbb{S}_{\beta}$, we have
\begin{equation}\label{eq:b.6.3a}
  \widehat{f}(z)=\frac1{\sqrt{2\pi}}
  \int_{\mathbb{R}}e^{i\lambda z}f(\lambda)\,d\lambda
  =\frac1{\sqrt{2\pi}}\langle e_{-i\overline{z}},f\rangle_{L^2(\mathbb{R})}
\end{equation}
for some $f\in\mathcal{K}$. Our aim is to compute the Szeg\"{o} kernel $\widehat{q}(z,w)=\widehat{q_w}(z)$ of $H^2(\mathbb{S}_{\beta})$ for $q_w\in\mathcal{K}$. By definition, for $z\in\mathbb{S}_{\beta}$ we have
$$\widehat{f}(z)=\langle \widehat{q_z},\widehat{f}\rangle_{H^2(\mathbb{S}_{\beta})}.$$
Using our characterization of $H^2(\mathbb{S}_{\beta})$ as a graph of $\widehat{M_{e_{-\beta}}}$, we have
\begin{equation}\label{eq:b.6.3b}
  \widehat{f}(z)=\langle \widehat{q_z},\widehat{f}\rangle_{H^2(\mathbb{S}_{\beta})}=\langle (1+\widehat{M_{e_{-2\beta}}})
  \widehat{q_z}, \widehat{f}\rangle_{L^2(\mathbb{R})}
  =\int_{\mathbb{R}}(1+e^{-2\beta\lambda})\overline{q_z(\lambda)}
  f(\lambda)\,d\lambda.
\end{equation}
Therefore, comparing \eqref{eq:b.6.3a} and \eqref{eq:b.6.3b} 
we find 
\[ q_z(\lambda)=\frac1{\sqrt{2\pi}}\frac{e^{-i\overline{z}\lambda}}{1+e^{-2\beta\lambda}},\quad\lambda\in\mathbb{R}. \] 
Using the Fourier transform, we obtain 
\begin{align*}
\widehat{q}(z,w)&=\widehat{q_w}(z)=\frac1{\sqrt{2\pi}}\int_{\mathbb{R}}e^{i\lambda z}q_w(\lambda)d\lambda=\frac1{2\pi}\int_\mathbb{R}\frac{e^{i\lambda(z-\overline{w})}}{1+e^{-2\lambda\beta}}d\lambda=\int_\mathbb{R}e^{i\lambda(z-\overline{w})}d\nu(\lambda)\\
&=\la e_{-i\oline z},e_{-i\oline w}\ra_{L^2(\R,\nu)},
\end{align*}
where
\begin{equation}
  \label{eq:sharp1}
  d\nu(\lambda):=\frac1{2\pi}\frac{d\lambda}{1+e^{-2\lambda\beta}}.
\end{equation}

\begin{lem} \mlabel{lem:FTcosh}
  The measure  $\nu$
on $\R$ is tempered, satisfies
$d\nu(-\lambda) = e^{-2\beta\lambda}d\nu(\lambda)$,
and its Fourier transform
defines a holomorphic function $\hat\nu$ on the strip $\bS_{2\beta}$, given by
\[ \hat \nu(z) = \frac{1}{2\pi}
  \int_\R \frac{ e^{iz\lambda}}{1 + e^{-2\beta\lambda}}\, d\lambda
    = \frac{1}{4\beta}
\frac{i}{\sinh\big(\frac{\pi z}{2\beta}\big)}
  = \frac{1}{4\beta}
    \frac{1}{\cosh\big(\frac{\pi z}{2\beta}-\frac{\pi i}{2}\big)}
\quad \mbox{ for } \quad z \in \bS_{2\beta}.\] 
  \end{lem}

 \begin{prf}
Starting from the well-known formula for the Fourier transform
of $\frac{1}{\cosh}$:
\begin{equation}
  \label{eq:cosh-fou}
 \int_\R \frac{e^{ix\xi}}{\cosh x}\, dx
 = \frac{\pi}{\cosh\big(\frac{\pi}{2}\xi\big)},
\end{equation}
we obtain for
$z \in \bS_{2\beta}$  
\begin{align*}
2\pi\hat\nu(z) &=  \int_\R \frac{e^{i\lambda z}}{1 + e^{-2\beta \lambda}}\, d\lambda
=  \int_\R \frac{e^{i\lambda z + \lambda\beta}}
      {e^{\beta\lambda} + e^{-\beta \lambda}}\, d\lambda
      = \frac{1}{2} \int_\R \frac{e^{i\lambda (z - \beta i)}}
      {\cosh(\beta\lambda)} \, d\lambda
      = \frac{1}{2\beta} \int_\R \frac{e^{i \frac{\lambda}{\beta}(z - \beta i)}}
    {\cosh(\lambda)} \, d\lambda\\
&    = \frac{\pi}{2\beta}
    \frac{1}{\cosh\big(\frac{\pi}{2}
      \big(\frac{1}{\beta}(-\beta i + z)\big)\big)} 
    = \frac{\pi}{2\beta}
    \frac{1}{\cosh\big(-\frac{\pi i}{2} + \frac{\pi z}{2\beta}\big)}
    = \frac{\pi}{2\beta}
    \frac{i}{\sinh\big(\frac{\pi z}{2\beta}\big)}.   
\qedhere\end{align*}
  \end{prf}

  By Lemma~\ref{lem:FTcosh},  the Szeg\"{o} kernel of
  $H^2(\mathbb{S}_{\beta})$ is given by 
\begin{equation}\label{eq:cskernel-strip} 
  \widehat{q}(z,w)=\widehat{q_w}(z)
  =\widehat{\nu}(z-\overline{w}) 
  =\frac{i}{4\beta}  \frac{1}{\sinh\big(\frac{\pi(z-\oline w)}{2\beta}\big)}
  =\frac{1}{4\beta}
    \frac{1}{\cosh\big(\frac{\pi (z-\oline w)}{2\beta}-\frac{\pi i}{2}\big)}.
\end{equation}

\begin{rem} Observe that $L^2(\R) \subeq L^2(\R,\nu)$, so that
  $\Gamma(M_{e_{-\beta}}) \subeq L^2(\R)^{\oplus 2}$.
  For the corresponding graph norm on $\cK$, we obain
  for $f,g\in\mathcal{K}$ and their images $(f,M_{e_{\beta}}f)$, $(g,M_{e_{\beta}}g)$ in $\Gamma(M_{e_{-\beta}})$ 
\begin{align*}
&\langle f,g\rangle_{L^2(\R,\nu)}+\langle M_{e_{-\beta}}f,M_{e_{-\beta}}g\rangle_{L^2(\R,\nu)}=\left\langle(1+M_{e_{-2\beta}})f,g\right\rangle_{L^2(\R,\nu)}\\
  &=\int_\mathbb{R}(1+e^{-2\beta\lambda})\overline{f(\lambda)}g(\lambda)d\nu
    =\frac{1}{2\pi} \int_\mathbb{R}(1+e^{-2\beta\lambda})\overline{f(\lambda)}g(\lambda)\frac{d\lambda}{1+e^{-2\lambda\beta}}=\frac{1}{2\pi} \langle f,g\rangle.
\end{align*}
\end{rem}

\subsection{The Poisson kernel of the strip} 

From \eqref{eq:cskernel-strip} in Appendix~\ref{subsec:app-cz-kernel},
we recall the Szeg\"o kernel of $\bS_\beta$: 
\[ 
  Q(z,w)
 = \frac{i}{4\beta}  \frac{1}{\sinh\big(\frac{\pi(z-\oline w)}{2\beta}\big)}
\quad \mbox{ and }  \quad
 Q(z,z) 
  =\frac{1}{4\beta\cos\big(\frac{\pi \Im z}{\beta}-
    \frac{\pi}{2}\big)}. \]

For the involution $\sigma(z) = \beta i + \oline z$, the fixed point set 
in $\bS_\beta$ is the line $\frac{\beta i}{2} + \R$, and on this line the 
kernel satisfies 
\[ Q\Big(x + \frac{\beta i}{2}, y + \frac{\beta i}{2}\Big)
= \frac{i}{4\beta} \frac{1}{\sinh\big(\frac{\pi(x-y)}{2\beta} 
+ \frac{\pi i}{2}\big)}
= \frac{1}{4\beta} \frac{1}{\cosh\big(\frac{\pi(x-y)}{2\beta}\big)}.\] 
On the imaginary axis, the kernel specializes for 
$0 < \lambda,\mu < \beta$ to 
\[ Q(i\lambda,i\mu) 
=  \frac{i}{4\beta}  \frac{1}{\sinh\big(\frac{i\pi(\lambda + \mu)}{2\beta}\big)}
=  \frac{1}{4\beta}  \frac{1}{\sin\big(\frac{\pi(\lambda + \mu)}{2\beta}\big)},
\] 
which is the kernel of this reproducing kernel Hilbert space 
obtained by restricting the Hardy space to $(0,\beta)i = \bS_\beta^{\sigma_2}$, where $\sigma_2(z):=-\bar{z}$. 

\begin{lem} The Poisson kernel on $\bS_\beta$ is given by 
\begin{equation}
  \label{eq:5.2.1}
 P_z(x) 
=\frac{1}{4\beta} \frac{\sin(\frac{\pi}{\beta} \Im z)}
{\sinh^2\big(\frac{\pi}{2\beta}(\Re z-x)\big)
  + \sin^2\big(\frac{\pi}{2\beta}\Im z\big)} \quad \mbox{ for } \quad x \in \R, 
\end{equation}
and 
\begin{equation}
  \label{eq:5.2.2}
 P_z(x + \beta i) 
=\frac{1}{4\beta} \frac{\sin(\frac{\pi}{\beta} \Im z)}
{\sinh^2\big(\frac{\pi}{2\beta}(\Re z-x)\big)
  + \cos^2\big(\frac{\pi}{2\beta}\Im z\big)}
\quad \mbox{ for } \quad x \in \R. 
\end{equation}
\end{lem}

\begin{prf}
  \begin{footnote}
{See also \cite{BK07} for the Poisson kernel on the strip 
$\{ z \in \C \: |\Im z| < 1 \}$.}
  \end{footnote}
For the computation of the Poisson kernel of $\bS_\beta$, we  recall that 
\[ \sinh(x + i y) 
= \sinh(x) \cos(y) + i \cosh(x) \sin(y), \] 
which leads to 
\begin{align}\label{eq:dagg}
 |\sinh(x + i y)|^2
&= \sinh^2(x) \cos^2(y) + \cosh^2(x) \sin^2(y) \notag \\
&= \sinh^2(x) \cos^2(y) +(1 +  \sinh^2(x)) \sin^2(y)\notag \\
&= \sinh^2(x) +\sin^2(y). 
\end{align}

For the Poisson kernel, we thus obtain for $x \in \R$ with 
\eqref{eq:cskernel-strip} and Hua's formula 
\eqref{eq:poiss-omega} 
\begin{align*}
P_z(x) 
&= \frac{|Q(z,x)|^2}{Q(z,z)} 
=\frac{1}{4\beta} \frac{\sin(\frac{\pi}{\beta} \Im z)}
{|\sinh\big(\frac{\pi}{2\beta}(z-x)\big)|^2}
=\frac{1}{4\beta} \frac{\sin(\frac{\pi}{\beta} \Im z)}
{\sinh^2\big(\frac{\pi}{2\beta}(\Re z-x)\big)
+ \sin^2\big(\frac{\pi}{2\beta}\Im z\big)}, 
\end{align*}
and for $\beta i + x$ in the upper boundary 
\begin{align*}
P_z(x + \beta i) 
= P(\beta i + \oline z,x) 
&=\frac{1}{4\beta} \frac{\sin(\frac{\pi}{\beta} \Im z)}
{\sinh^2\big(\frac{\pi}{2\beta}(\Re z-x)\big)
+ \cos^2\big(\frac{\pi}{2\beta}\Im z\big)}
\qedhere\end{align*}
\end{prf}

For $z = \frac{\beta i}{2} + \lambda$, Equation~\eqref{eq:5.2.1}
specializes to 
\begin{align}  \label{eq:hor-poisson}
P_{\lambda + \frac{\beta i}{2}}(x) 
&=\frac{1}{4\beta} \frac{1} 
{\sinh^2\big(\frac{\pi}{2\beta}(\lambda-x)\big)
  + \sin^2\big(\frac{\pi}{4}\big)}
=\frac{1}{4\beta} \frac{1} 
{\sinh^2\big(\frac{\pi}{2\beta}(\lambda-x)\big) + \frac{1}{2}} \\
  &=  P_{\lambda + \frac{\beta i}{2}}(x + \beta i).\notag
\end{align}
With 
  \[ \frac{\cosh(2x)}{2} = \frac{\cosh(x)^2 + \sinh(x)^2}{2} 
  = \frac{1}{2} +  \sinh(x)^2 \]
\eqref{eq:hor-poisson} can be rewritten as 
  \begin{equation}
    \label{eq:new-poisson} P_{\lambda + \frac{\beta i}{2}}(x) 
=\frac{1}{4\beta} \frac{1} 
{\sinh^2\big(\frac{\pi}{2\beta}(\lambda-x)\big) + \frac{1}{2}} 
=\frac{1}{2\beta} \frac{1} 
{\cosh\big(\frac{\pi}{\beta}(\lambda-x)\big)}.
  \end{equation}

\begin{rem} \mlabel{rem:cos2}
With the Poisson representation we obtain for $t \in \R$ 
\begin{align*}
 e_{it}\Big(\frac{\beta i}{2} + \lambda\Big)
&  = e^{-t\beta/2} e^{it\lambda} 
= \int_\R P_{\frac{\beta i}{2} + \lambda}(x) e^{itx}\, dx 
 +  \int_\R P_{\frac{\beta i}{2} + \lambda}(\beta i + x) e^{itx} e^{-t\beta}\, dx\\
& = (1 + e^{-t\beta}) \int_\R P_{\frac{\beta i}{2} + \lambda}(x) e^{itx}\, dx,
\end{align*}
and therefore
\begin{equation}
  \label{eq:cosh}
 \int_\R P_{\frac{\beta i}{2} + \lambda}(x) e^{itx}\, dx
  = \frac{e^{-t\beta/2} e^{it\lambda}}{1 + e^{-t\beta}}
  = \frac{e^{it\lambda}}{e^{t\beta/2} + e^{-t\beta/2}}
  = \frac{e^{it\lambda}}{2\cosh(t\beta/2)}.
\end{equation}
For $\lambda = 0$ this formula  also describes the Fourier transform of
$\frac{1}{\cosh}$ in terms of $P_{\beta i/2}$; see also \eqref{eq:cosh-fou}. 
\end{rem}

\subsection{More Fourier transforms}

We start with the Fourier transform of the square of 
$\sech = \cosh^{-1}$:
\[ \frac{1}{\sqrt{2\pi}} \int_\R
  \frac{e^{ix\lambda}}{\cosh(x)^2}\, dx
  = \frac{\sqrt{\pi}}{\sqrt 2}
  \frac{\lambda}{\sinh\big(\frac{\pi\lambda}{2}\big)}
  = \sqrt{2\pi} 
  \frac{\lambda}{e^{\pi \lambda/2} - e^{-\pi \lambda/2}}.\]
One can compute such Fourier transform by noticing that the derivative of $\tanh$ is $\frac1{\cosh^2}$ and compute the distributional
Fourier transform of $\tanh$ instead.
It is immediately seen that $\tanh$ has countable singularities in the imaginary axis. Using the Residue Theorem \cite[Thm.~10.42]{Ru86}, respectively on the upper and lower half-plane, one gets the desired conclusion\footnote{For similar computations in the case of the Fourier transform of $\frac1{\cosh}$, see\\\url{https://math.stackexchange.com/questions/799616/fourier-transform-of-1-cosh}}. 
\smallskip

Inverting the Fourier transform leads to
\begin{equation}
  \label{eq:invfout}
  \frac1{\sqrt{2\pi}}\int_\R   \frac{e^{-i\lambda x}\lambda}{e^{\pi \lambda/2} - e^{-\pi \lambda/2}}\,
  d\lambda = \frac{1}{\cosh(x)^2}.
\end{equation}

The Bergman kernel on the strip $\bS_\beta$ is the square of the
Szeg\"o kernel $Q_1(z,w)$, hence
\[ Q_2(z,w)
  = \frac{1}{(4\beta)^2}
  \frac{1}{\cosh^2
    \big(\frac{\pi (z-\oline w)}{2\beta}-\frac{\pi i}{2}\big)}.\]
For $q_2(z) := Q_2(z,0)$, this leads with \eqref{eq:invfout} to 
\[ q_2(z)
  = \frac{1}{(4\beta)^2}
  \int_\R   \frac{e^{i\lambda(\frac{\pi z}{2\beta} - \frac{\pi i}{2})}
    \lambda}{e^{\pi \lambda/2} - e^{-\pi \lambda/2}}\,  d\lambda 
  =  \frac{1}{(4\beta)^2} 
  \int_\R e^{i\lambda(\frac{\pi z}{2\beta})}
  \frac{\lambda\, d\lambda}{1 - e^{-\pi \lambda}}
  = \frac{1}{(2\pi)^2} 
  \int_\R e^{i\lambda z} 
  \frac{\lambda\, d\lambda}{1 - e^{-2\beta \lambda}}
  = \hat\nu(z)\]
for
\[ d\nu(\lambda) = \frac{1}{(2\pi)^2}
  \frac{\lambda\, d\lambda}{1 - e^{-2\beta \lambda}}.\]
For the Bergman kernel we thus obtain
\begin{equation}
  \label{eq:bergman}
  Q_2(z,w) = \hat\nu(z- \oline w).
\end{equation}

\begin{rem} (Higher order kernels) 
An elementary calculation shows that for $n\geq1$ 
\[ \frac{1}{\cosh(x)^{n+2}}
  = \frac{1}{n(n+1)} \Big(n^2 - \frac{d^2}{dx^2}\Big) \frac{1}{\cosh(x)^n}.\]
For the Fourier transform
$\cF(f) = \hat f$, this leads to
\[ \cF(\cosh^{-n-2})(p) 
  = \frac{n^2 + p^2}{n(n+1)} \cF(\cosh^{-n})(p).\] 
We thus obtain inductively if $n=2(k-1)$ for $k\geq2$
\[   \cF(\cosh^{-2k})(p) 
  =\frac{(2^2 + p^2)\cdots ((2k-2)^2 + p^2)}{(2k-1)!}
  \cF(\cosh^{-2})(p).\]
One likewise obtains a formula for the Fourier transform
of $\cosh^{-(2k+1)}$. 
\end{rem}

\end{document}